\documentclass[11pt]{article}
\usepackage[utf8]{inputenc}
\usepackage[T1]{fontenc}
\usepackage{latexsym,amssymb,amsmath,amsfonts,amsthm}
\usepackage{graphics}
\usepackage[section]{placeins}
\usepackage{graphicx}
\usepackage{mathrsfs}
\usepackage{subfigure}
\usepackage{extarrows}
\usepackage{float}
\usepackage{color}
\usepackage{framed}
\usepackage{tikz}
\usetikzlibrary{calc}
\usepackage{epstopdf}
\newcommand{\R}{{\mathbb R}}

\newcommand{\N}{{\mathbb N}}

\newcommand{\be}{\begin{eqnarray}}
\newcommand{\ben}{\begin{eqnarray*}}
\newcommand{\en}{\end{eqnarray}}
\newcommand{\enn}{\end{eqnarray*}}

\newcommand{\real}{{\rm Re\,}}
\newcommand{\ima}{{\rm Im\,}}

\newcommand{\s}{\mathbb{S}}
\newcommand{\G}{\Gamma}

\newcommand{\range}{{\rm Range}}

\newtheorem{theorem}{Theorem}[section]

\newtheorem{remark}[theorem]{Remark}

\newtheorem{lemma}[theorem]{Lemma}
\newtheorem{corollary}[theorem]{Corollary}

\definecolor{rot}{rgb}{0,0,0}
\definecolor{hw}{rgb}{0,0,0}
\definecolor{rot1}{rgb}{0,0,0}
\newcommand{\tcr}{\textcolor{rot}}
\newcommand{\rot}{\textcolor{rot1}}
\newcommand{\tcb}{\textcolor{hw}}

\begin{document}
\renewcommand{\theequation}{\arabic{section}.\arabic{equation}}
\begin{titlepage}
\title{\bf Inverse wave-number-dependent source problems for the Helmholtz equation with partial information on radiating period}

\author{
 Hongxia Guo\thanks{School of Mathematical Sciences and LPMC, Nankai University, 300071 Tianjin, China. ({\tt hxguo@nankai.edu.cn})}  \and Guanghui Hu\thanks{School of Mathematical Sciences and LPMC, Nankai University, 300071 Tianjin, China. ({\tt ghhu@nankai.edu.cn})} \and
 Guanqiu Ma\thanks{Corresponding author: School of Mathematical Sciences and LPMC, Nankai University, 300071 Tianjin, China. ({\tt gqma@nankai.edu.cn})}
 }
\date{}
\end{titlepage}
\maketitle


\begin{abstract}
This paper addresses a factorization method for imaging the support of a wave-number-dependent source function from  multi-frequency data measured at a finite pair of symmetric receivers in opposite directions. The source function is given by the inverse Fourier transform of a compactly supported time-dependent source whose initial moment or terminal moment for radiating is unknown.  Using the multi-frequency far-field data at two opposite observation directions, we provide a computational criterion for characterizing the smallest strip containing the support and perpendicular to the directions.
A new parameter is incorporated into the design of test functions for indicating the unknown moment.
The data from a finite pair of opposite directions can be used to recover the $\Theta$-convex polygon of the support.
Uniqueness in recovering the convex hull of the support is obtained as a by-product of our analysis using all observation directions.
Similar results are also discussed with the multi-frequency near-field data from a finite pair of  observation positions in three dimensions.
We further comment on possible extensions to source functions with two disconnected supports. Extensive numerical tests in both two and three dimensions are implemented to show effectiveness and feasibility of the approach. The theoretical framework explored here should be seen as the frequency-domain analysis for inverse source problems in the time domain.

\vspace{.2in} {\bf Keywords: Inverse source problem, Helmholtz equation, wave-number-dependent sources, multi-frequency data, factorization method.}
\end{abstract}

\section{Introduction}
\subsection{Problem formulation}
Consider the inhomogeneous Helmholtz equation with a wave-number-dependent source term
 	\begin{equation}\label{Helmholtz}
 		\Delta w(x,k) + k^2 w(x,k) = - f(x,k), \qquad x\in \R^{3},
 	\end{equation}
where $k>0$ is called the wave-number. In this paper  the source function $f(x,k)$ is supposed to be the inverse Fourier transform of a time-dependent source $S(x,t)$ with compact support in both temporal and spatial variables. More precisely, we make the following assumptions:
\begin{itemize}
\item $\mbox{supp}\, S(x,t)=\overline{D} \times [t_{\min}, t_{\max}]\subset B_R\times \R_+$ with some $t_{\max}>t_{\min}\geq 0$ and some $R>0$. Here $B_R=\{x\in \R^3: |x|<R\}$.
\item $D\subset \R^3$ is a bounded Lipschitz domain such that $\R^3\backslash\overline{D}$ is connected.
\item $S(x,t)\in C([t_{\min}, t_{\max}], L^\infty(D))$ is a \tcb{real-valued function} fulfilling the positivity constraint
\begin{equation}\label{F}
S(x, t)\geq c_0>0\qquad\mbox{a.e.}\; x\in \overline{D},\quad t\in [t_{\min}, t_{\max}].
\end{equation}
\end{itemize}
Hence, the source function in the frequency-domain takes the integral form
\begin{equation}\label{fxt}
f(x,k)=\frac{1}{\sqrt{2\pi}}\int_{\R} S(x,t)e^{{\rm i} k t}dt=\frac{1}{\sqrt{2\pi}}\int_{t_{\min}}^{t_{\max}} S(x,t)e^{{\rm i} k t}dt.
\end{equation}
This implies that the function $k\mapsto f(x, k)$ is real-analytic for every $x\in \R^3$.
The Helmholtz equation \eqref{Helmholtz} arises from the inverse Fourier transform of the time-dependent acoustic wave equation with the function $S(x, t)$ as the forcing term and  satisfying the homogeneous initial conditions (see \eqref{TU} in the Appendix).
The unique solution to \eqref{Helmholtz} satisfies the Sommerfeld radiation condition
 	\be\label{SRC}
	\lim\limits_{r \to \infty} r (\partial_r w - ikw) = 0,\quad r = |x|,\en
which holds uniformly in all directions $x/|x|$. In fact, the radiating behavior of $w$ at a fixed frequency can be  derived from the inverse Fourier transform of the outgoing solution to the time-dependent wave equation in $\R^3$ (\cite{GGH2022}).
Moreover, $w$ is of  a convolution form in spatial variables:
\begin{equation}\label{expression-w}
 	w(x,k)= \int_{\mathbb{R}^3} \Phi(x-y;k) f(y,k)\,dy,
\end{equation}
where $\Phi(x;k)$ is the fundamental solution to the Helmholtz equation $(\Delta + k^2)w = 0$, given by
\begin{equation*}
		\Phi(x;k) = \frac{e^{ik|x|}}{4 \pi |x|},\quad x\in \mathbb{R}^3,\, |x| \neq 0.
	\end{equation*}
The Sommerfeld radiation conditions \eqref{SRC} for $w$ and $\Phi$ give rise to the following asymptotic behavior at infinity:
\begin{equation}\label{far-field}
w(x,k)=\frac{e^{{\rm i} k |x|}}{4\pi|x|}\left\{w^\infty(\hat{x},k)+O\left(\frac{1}{|x|}\right)\right\}\quad\mbox{as}\quad|x|\rightarrow\infty,
\end{equation}
where $w^\infty(\cdot, k)\in C^\infty(\s^2)$ is referred to as the far-field pattern (or scattering amplitude) of $w$. It is well known that the function  $\hat{x}\mapsto w^\infty(\hat{x}, k)$ is real analytic on $\s^2$, where $\hat{x}\in \s^2$ is usually called the observation direction. In the appendix we shall prove that $w^\infty(\hat{x}, k)$ coincides with the inverse Fourier transform of the time-dependent far-field data  in terms of the time variable. By \eqref{expression-w}, the far-field pattern $w^\infty$ can be expressed as
\be\label{u-infty}
w^\infty(\hat{x}, k)=\int_{D} e^{-{\rm i}k\hat{x}\cdot y} f(y, k)\,dy,\quad \hat{x}\in \s^2,\quad k>0.
\en
Since the time-dependent source $S$ is real valued, we have $f(x, -k)=\overline{f(x,k)}$  and thus $w^\infty(\hat{x}, -k)=\overline{w^\infty(\hat{x}, k)}$ for all $k\in \R$. The far-field pattern
$w^\infty(\hat{x}, k)$ depends analytically on $k\in \R$, because $D$ is bounded and $f$ is real-analytic with respect to $k\in \R$.

In our previous work \cite{GGH2022} we have studied the inverse source problem of identifying $\partial D$ from the multi-frequency data detected at one or several (but not necessarily symmetric) observation directions/points, provided the source radiating period $[t_{\min}, t_{\max}]$
is completely known. We note that both the initial moment $t_{\min}$ and the terminal moment $t_{\max}$ are essentially required in \cite{GGH2022}. The primary task of this paper is to relax this condition by assuming that one of $t_{\min}$ and $t_{\max}$ is not given.
Let $0\leq k_{\min}<k_{\max}$ and denote by $[k_{\min}, k_{\max}]$ the bandwidth of available wave numbers of the Helmholtz equation. The inverse source problems to be considered within this paper are described as follows:
\begin{framed}
(ISPs): Extract information on the position and shape of the support $D$ of $S(\cdot, t)$ from knowledge of the multi-frequency far-field patterns
$$\left\{w^\infty(\pm\hat{x}_j, k): k\in[k_{\min}, k_{\max}], \, j=1,2,\cdots, J \right\}, $$
or from the multi-frequency near-field data
$$\left\{w(\pm{x}_j, k): k\in[k_{\min}, k_{\max}], \,|x_j|=R,\; j=1,2,\cdots, J \right\}. $$
Here, $2J\in \N$ denotes the number of opposite observation directions/positions and either the initial or the terminal moment is not given.

\end{framed}

\subsection{Literature review and comments}

If the time-dependent source function is of the form $S(x,t)=s(x)\delta(t-t_0)$ (which corresponds to the critical case that $t_{\min}=t_{\max}=t_0$), the source function in the frequency domain can be expressed as $f(x,k)=s(x) e^{ikt_0}/\sqrt{2\pi}$. In the special case that $t_0=0$, the source term $f(x,k)=s(x)/\sqrt{2\pi}$ is independent of the wave-number. If the impulse moment $t_0$ is given, the function $u(x,k)=w(x,k)e^{-ikt_0}\sqrt{2\pi}$ turns out to be the radiation solution of the Helmholtz equation $\Delta u+k^2u=s$ whose right hand side is again independent of the wave-number. The far-field pattern of $u$ is nothing  else but the Fourier transform of the space-dependent source term $s$ at the Fourier variable  $\xi=k\hat{x}\in \R^3$ multiplied by some constant.
A wide range of literatures is devoted to such inverse wave-number-independent source problems with multi-frequency data, for example,  uniqueness proofs and increasing stability analysis with near-field measurements \cite{BLLT, BLT10, CIL, EV09, LY} and a couple of numerical schemes such as iterative method, Fourier method and test-function method
for recovering the source function \cite{BLLT, BLRX,  EV09, ZG}  and sampling-type methods for imaging the support  \cite{AHLS, GS, LMZ, JLZ2019}.
We shall establish a multi-frequency factorization method for imaging $D$ from a finite number of observation directions in the far field or a finite number of near-field positions. This approach was proposed by A. Kirsch in 1998 \cite{K98, KG08} with the multi-static data at a fixed energy.  Later it has been extensively studied in various inverse time-harmonic scattering problems using far-field patterns over all observation directions (or equivalently, the Dirichlet-to-Neumann map for elliptic equations).

For inverse source problems in the frequency domain, to the best of our knowledge, little is known if the source function of the Helmholtz equation depends on both frequency/wave-number and spatial variables. Essential difficulties arise from the expression \eqref{u-infty}, where the far-field pattern is no longer the Fourier transform of the source function. For source terms of the integral form \eqref{fxt} with complete knowledge on the radiating period, a factorization method was established in \cite{GGH2022} for imaging the so-called $\Theta$-convex polygon (see \cite{GS} for the original definition) of the support associated with a finite number of observation directions. This extends the multi-frequency factorization scheme of \cite{GS, GS17} from wave-number-independent sources to wave-number-dependent ones.
In particular, the multi-frequency data at a single direction can be used to rigorously characterize the smallest strip that contains the support $D$ and perpendicular to the observation direction. The goal of this paper is to establish the same factorization method but with only partial information on the radiation period (that is, either $t_{\min}$ or $t_{\max}$). The price we pay is to introduce a new parameter into the design of test functions and utilize the multi-frequency far-field data at two opposite directions. Moreover, we show that the unknown moment $t_{\min}$ or $t_{\max}$ can also be recovered by a indicator function on the new parameter in the test function and by using the data of two opposite directions. Another contribution of this paper is to address how to extend our approach to source functions with disconnected supports. A direct sampling method was proposed in \cite{GHZ} for treating (ISPs) with complete knowledge on the radiating period, which we think remains valid under the weaker assumptions of this paper.

It remains interesting to investigate the inverse source problems (ISPs) when both $t_{\min}$ and $t_{\max}$ are not given, in particular, to extract source information from a single direction/position.
Our approach does work, while it relies heavily on the information of the initial moment $t_{\min}$ or the terminal moment $t_{\max}$. Intuitively, we think that two different parameters should be incorporated into the test function if none of them is known, but the theoretical frame is still unclear to us. However, in the special case that $t_{\min}=t_{\max}=t_0$ and $S(x,t)=s(x)\delta(t-t_0)$
with an unknown moment $t_0\in \R_+$, it is possible to recover $D=\mbox{supp}(s)$ and $t_0$ using the same data considered in this paper.
Another open problem is to establish the multi-frequency factorization method for inverse medium and inverse obstacle scattering problems for which the right hand side $f(x,k)$ depends also on the solution $w(x,k)$. The lack of the positivity condition \eqref{F} brings difficulties in extending our idea to inverse medium problems, but the
framework of this paper can be used to handle high frequency or weak scattering models for inverse obstacle scattering problems. It is worth noting that the factorization method explored here differs from that for  acoustic wave scattering by obstacles in the time domain (see \cite{CHL2019, HL2020}). In these works the time-dependent causal scattered waves at infinitely many receivers are needed to defined a modified far-field operator. Our attention is to extract the support information of a time-dependent source from wave signals recorded at one or several receivers. The approach considered here should be seen as a frequency-domain analysis for inverse source problems in the time domain. In the recent work \cite{GHM}, our approach has been adapted to inverse moving point sources problems in the frequency domain via Fourier transform.

The remaining part is organized as follows. In Section \ref{sec2}, the factorizations of the multi-frequency far/near-field operator are reviewed. Section \ref{sec3} is devoted to the design of test functions and indicator functions using multi-frequency far-field data received at a single pair or a finite pair of symmetric directions.  Uniqueness and inversion algorithms will be also addressed.  In Section \ref{sec4},  the corresponding inverse problems using  multi-frequency near-field data are discussed.
We comment on possible extensions of our reconstruction method  to source functions with two disconnected supports in Section \ref{mul-compo}. Finally, a couple of numerical tests will be reported in Section \ref{num}.

Below we introduce some notations to be used throughout this paper. Unless otherwise stated, we always suppose that $D$ is bounded and connected. Given $\hat{x}\in \s^2$, we define
\ben
\hat{x}\cdot D:=\{t\in \R: t=\hat{x}\cdot y\;\mbox{for some}\; y\in D\}\subset \R.
\enn
Hence, $(\inf(\hat{x}\cdot D),  \sup(\hat{x}\cdot D))$ must be a finite and connected interval on the real axis.
A ball centered at $y\in \R^3$ with the radius $\epsilon>0$ will be denoted as $B_\epsilon(y)$. For brevity we write $B_\epsilon=B_\epsilon(0)$ when the ball is centered at the origin.
Obviously, $\hat{x}\cdot B_\epsilon(y)=(\hat{x}\cdot y-\epsilon,\hat{x}\cdot y+\epsilon)$.
In this paper the one-dimensional Fourier and inverse Fourier transforms are defined by
\begin{equation*}
(\mathcal{F}f)(k)=\frac{1}{\sqrt{2\pi}}\int_{\R}f(t)e^{-{\rm i} k t}\,dt,\quad
(\mathcal{F}^{-1}v)(t)=\frac{1}{\sqrt{2\pi}}\int_{\R}v(k)e^{{\rm i} k t}\,dk,
\end{equation*}
respectively.

\section{Review of the factorization of far-field and near-field operators}\label{sec2}
In this section, we will introduce the far/near-field operator with multi-frequency data at a fixed observation direction/position and review its connection with the data-to-pattern operator by a range identity. We refer to \cite{GS} for discussions on source functions independent of the wave-number and to \cite{GGH2022} for wave-number-dependent source functions with complete a priori information on the radiating period.
Introduce the central frequency  $k_c$ and  half of the bandwidth of the multi-frequency data as (see \cite{GS})
\ben
k_c:=\frac{k_{\min}+k_{\max}}{2},\quad K :=\frac{k_{\max}-k_{\min}}{2}.
\enn
These parameters were firstly introduced in \cite{GS},
allowing us to define far-field and near-field operators with the same domain (i.e., $L^2(0, K)$) of definition and range.

For every fixed direction $\hat{x}\in \s^2$, we define the far-field operator by (see \cite{GS})
\be \label{def:F}
(F\phi)(\tau)=(F_D^{(\hat{x})}\phi)(\tau):=\int_{0}^{K } w^\infty(\hat{x}, k_c+\tau-s)\,\phi(s)\,ds.
\en
Since $w^\infty(\hat{x},k)$ is analytic with respect to the wave number $k\in \R$, the operator
$F^{(\hat{x})}: L^2(0, K )\rightarrow L^2(0, K )$ is bounded.
For notational convenience we introduce the space
$$X_D:=L^2(D\times(t_{\min}, t_{\max})).$$
Denote by $\langle \cdot,  \cdot\rangle_{X_D}$ the inner product over $X_D$.
Below we show a factorization of the multi-frequency far-field operator at a fixed observation direction.
\begin{lemma}\label{Fac-F}(\cite{GGH2022})
We have $F=L\mathcal{T}L^*$, where $L=L_D^{(\hat{x})}: X_D\rightarrow L^2(0, K )$ is defined by
\be\label{def:L}
(Lu)(\tau)=\int_{t_{\min}}^{t_{\max}}\int_D e^{{\rm i} \tau (t-\hat{x}\cdot y)} u(y,t)\,dy \,dt,\qquad \tau\in [0, K]
\en
for all $u\in X_D$,
and $\mathcal{T}: X_D\rightarrow X_D$ is a multiplication operator defined by
\be
(\mathcal{T}u)(y,t):=\frac{1}{{\sqrt{2\pi}}} e^{{\rm i} k_c (t-\hat{x}\cdot y)}\,S(y, t)\,u(y, t) .
\en
\end{lemma}
The operator $L=L_D^{(\hat{x})}$ will be referred to as the data-to-pattern operator, because it maps the time-dependent source function to the multi-frequency far-field data at a fixed observation direction, that is,
\ben
w^\infty(\hat{x}, k)=\frac{1}{\sqrt{2\pi}}(L_D^{(\hat{x})}S)(k).
\enn
Lemma \ref{Fac-F} implies that the far-field operator $F$ is
self-adjoint. It is also
 positive if the positivity constraint \eqref{F} holds.
 Define $F_{\#}:=|\real F|+|\mbox{Im}\; F|$.
Using the range identity proved in \cite{GGH2022}, we obtain the relation
\be\label{RI}
\mbox{Range}\, (F_{\#}^{1/2})=\mbox{Range}\,(L).\en

Let $\chi_y(k)\in L^2(0, K )$
be a $y$-dependent test function with $y\in \R^3$ and fix the observation direction $\hat{x}\in \s^2$.
Denote by $(\lambda_n^{(\hat{x})}, \psi_n^{(\hat{x})})$ an eigensystem of the positive and self-adjoint operator $ F_{\#}$, which is uniquely determined by the multi-frequency far-field patterns $\{w^\infty(\hat{x}, k): k\in [k_{\min}, k_{\max}]\}$. Applying Picard's theorem and range identity, we obtain
\be\label{In1}
\chi_y \in \range (L)\quad\mbox{if and only if}\quad \sum_{n=1}^\infty\frac{|\langle \chi_y, \psi_n^{(\hat{x})} \rangle|^2}{ |\lambda_n^{(\hat{x})}|}<\infty.
\en
In the subsequent Section 3 we shall choose proper test functions $\chi_y$ such that $\chi_y \in \range (L)$ if and only if the sampling variable $y$ belongs to a domain associated with the support of the source function.
The support of the Fourier transform of the function $Lu$ can be evaluated as follows (see \cite{GGH2022}):
\begin{equation}\label{rl-f}
	\mbox{supp} (\mathcal{F}(Lu)) \subset \left[t_{\min} - \sup(\hat{x}\cdot D),\; t_{\max} - \inf(\hat{x}\cdot D) \right]\quad\mbox{for all}\quad u\in X_D.
\end{equation}

The above results can be carried over the near-field case naturally.
Given a fixed observation point ${x}\in \partial B_R:=\{x\in \R^3, |x|=R\}$, we define the near-field operator by
\be \label{def:N}
(\mathcal{N}\phi)(\tau)=(\mathcal{N}_D^{({x})}\phi)(\tau):=\int_{0}^{K } w({x}, k_c+\tau-s)\,\phi(s)\,ds.
\en
Since $w({x},k)$ is also analytic with respect to the wave number $k\in \R$, the operator
$\mathcal{N}: L^2(0, K )\rightarrow L^2(0, K )$ is bounded.
Below we show a factorization of the near-field operator.
\begin{lemma}\label{Fac-N}(\cite{GGH2022})
We have $\mathcal{N}=\tilde{L}\tilde{T}\tilde{L}^*$, where $\tilde{L}=\tilde{L}_D^{({x})}: X_D\rightarrow L^2(0, K )$ is defined by
\ben
(\tilde{L}u)(\tau)=\int_{t_{\min}}^{t_{\max}}\int_D e^{i\tau (t+|x-y|)} u(y,t)\,dy dt,\qquad \tau\in [0, K]
\enn
for all $u\in X_D$,
and $\tilde{T}: X_D\rightarrow X_D$ is a multiplication operator defined by
\ben
(\tilde{T}u)(y,t):=\frac{e^{{\rm i} k_c (t+|x-y|)}}{\sqrt{32\pi^3} |x-y|}\,u(y, t)\, S(y, t),\qquad |x|=R.
\enn
\end{lemma}
Analogously, the operator $\tilde{L}$ can be interpreted as the data-to-near-field operator, since it maps the time-dependent source function to the multi-frequency data at the receiver $x$, that is,
\ben
w(x,k)=\left(\tilde{L}_D^{(x)} \frac{S(y,t)}{\sqrt{32\pi^3}|x-y|}\right)(k).
\enn
Again applying the range identity of \cite{GGH2022} yields the relation
\be\label{RN}
\mbox{Range}\, (\mathcal{N}_{\#}^{1/2})=\mbox{Range}\,(\tilde{L}).
\en

Let $\chi_y(k)\in L^2(0, K )$ be a test function relying on the sampling variable $y\in \R^3$, and
denote by $(\lambda_n^{({x})}, \psi_n^{({x})})$ an eigensystem of the positive and self-adjoint operator $ \mathcal{N}_{\#}$.
Similar to the far-field case (see \eqref{In1}), one obtains
\ben
\chi_y \in \range (\tilde{L})\quad\mbox{if and only if}\quad \sum_{n=1}^\infty\frac{|\langle \chi_y, \psi_n^{({x})} \rangle|^2}{ |\lambda_n^{({x})}|}<\infty.
\enn
In our analysis, we also need to estimate the support of the Fourier transform of $\tilde{L}^{(x)}_D u$ (see \cite{GGH2022}):
\begin{equation}\label{rl-n}
	\mbox{supp }\mathcal{F}(\tilde{L}^{(x)}_D u) \subset \left[\inf\limits_{z\in D} |x-z| +t_{\min},\; \sup\limits_{z\in D} |x-z| +t_{\max}\right],\quad\mbox{for all}\quad u\in X_D.
\end{equation}

\section{Test and indicator functions with multi-frequency far-field data}\label{sec3}

The aim of this section is to define test and indicator functions  with multi-frequency far-field data when the terminal moment $t_{\max}$  or the initial moment $t_{\min}$  is unknown. We shall consider the two cases separately.
\subsection{The terminal moment $t_{\max}$ is unknown}
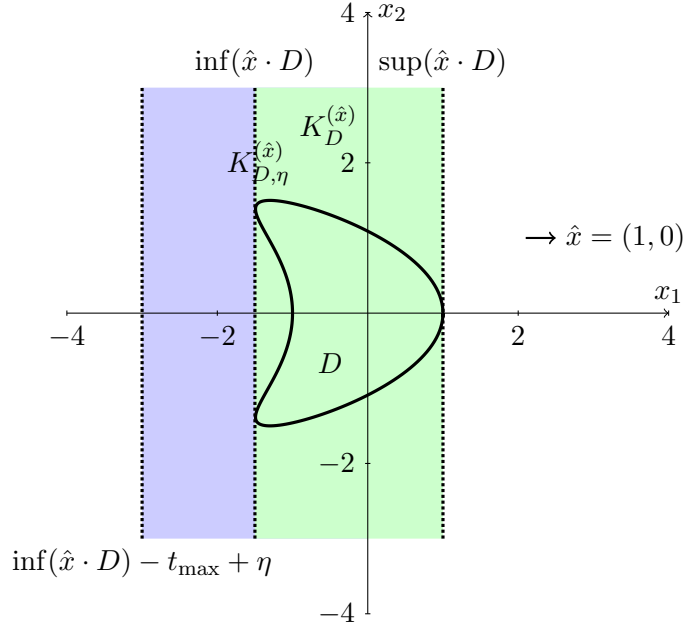
\begin{figure}[H]
		\centering
		\scalebox{1}{
		\begin{tikzpicture}
		\draw[line width=3cm,color=blue!20] (-1.5,-3) -- (-1.5,3);
		\draw[line width=2.5cm,color=green!20] (-0.25,-3) -- (-0.25,3);
		\draw (0,2.5) node [left] {$K_{D}^{(\hat{x})}$};
		\draw (-0.9,2) node [left] {$K _{D,\eta}^{(\hat{x})}$};
		\draw[->] (-4,0) -- (4,0) node[above] {$x_1$} coordinate(x axis);
		\draw[->] (0,-4) -- (0,4) node[right] {$x_2$} coordinate(y axis);
		\foreach \x/\xtext in {-4,-2, 2, 4}
		\draw[xshift=\x cm] (0pt,1pt) -- (0pt,-1pt) node[below] {$\xtext$};
		\foreach \y/\ytext in {-4,-2, 2, 4}
		\draw[yshift=\y cm] (1pt,0pt) -- (-1pt,0pt) node[left] {$\ytext$};
	
		 \draw[domain = -2:360,very thick][samples = 200] plot({cos(\x)+0.65*cos(2* \x)-0.65},{1.5*sin(\x)});
		\draw [very thick, densely dotted] (-1.5,-3) -- (-1.5,3);
		\draw [very thick, densely dotted] (1,-3) -- (1,3);
		\draw [very thick, densely dotted] (-3,-3) -- (-3,3);
		\draw (-1.5,3) node [above] {$\inf (\hat{x}\cdot D)$};
		\draw (1,3) node [above] {$\sup (\hat{x}\cdot D)$};
		\draw (-3,-3) node [below] {$\inf (\hat{x}\cdot D)-t_{\max}+\eta$};
		\draw (-0.5,-0.4) node [below] {$D$};

		\draw [thick,->] (2.1,1) -- (2.5,1);
		\draw (2.5,1) node[right] {$\hat{x} = (1,0)$};
		
		\end{tikzpicture}
		}
		\caption{Illustration of the strips $K_D^{(\hat{x})}$ (green area) and $K _{D,\eta}^{(\hat{x})}$ (union of green and blue area) with $\hat{x} = (1, 0)$ when $\eta<t_{\max}$.  }
		\label{strip}
	\end{figure}

Let $\hat{x}\in \s^2$ be a fixed observation direction and assume that the initial moment $t_{\min}$ is known.
For $y\in \R^3$ and $\epsilon>0$, define  the test function $\phi^{(\hat{x})}_{y, \eta,\epsilon}\in L^2(0, K )$  by
\be\label{Test}
\phi^{(\hat{x})}_{y,\eta, \epsilon}(k)=\frac{1}{(\eta -t_{\min}) \,|B_\epsilon(y)|}\int_{t_{\min}}^{\eta}\int_{B_\epsilon(y)}e^{{\rm i} k  (t-\hat{x}\cdot z)} dzdt,\qquad k\in[0, K],
\en
where $\eta>t_{\min} $ is a parameter and $|B_\epsilon(y)|=4/3\pi \epsilon^3$ denotes the volume of the ball $B_\epsilon(y)\subset \R^3$.
The parameter $\eta$ will be taken in a small neighborhood on the right hand side of  $t_{\min}$. As $\epsilon\rightarrow 0$,  there holds the uniform convergence in $L^2(0, K)$:
\begin{equation}\label{test}
	\begin{aligned}
		\phi^{(\hat{x})}_{y,\eta, \epsilon}(k)\rightarrow
		\phi^{(\hat{x})}_{y,\eta}(k)&:=\frac{1}{\eta -t_{\min}}\left(\int_{t_{\min}}^{\eta }e^{{\rm i} k  t} dt\right) \;e^{-{\rm i} k  \hat{x}\cdot y} \\
		&=\frac{1 }{ik (\eta-t_{\min }) } \left( e^{{\rm i} k  \eta }- e^{{\rm i} k  t_{\min}}\right)\,e^{-{\rm i} k  \hat{x}\cdot y}.
	\end{aligned}
\end{equation}
\begin{remark}
The choice of $\eta=t_{\max}$ was taken in \cite{GGH2022} when both $t_{\min}$ and $t_{\max}$ are known. If $t_{\max}$ is unknown, we design the $\eta$-dependent test function \eqref{Test} for recovering both $k_{\max}$ and the support of the source function, which makes this paper quite different from our previous work \cite{GGH2022}.
\end{remark}
Below we describe the supporting interval of the Fourier transform of the test functions defined by \eqref{Test}.
\begin{lemma} For $\epsilon>0$,
we have
\be\label{a1}
&&[\mathcal{F} \phi^{(\hat{x})}_{y, \eta,\epsilon}](\xi)>0\quad \mbox{if}\quad \xi\in (t_{\min} - \hat{x}\cdot y-\epsilon,\; \eta - \hat{x}\cdot y+\epsilon ),\\ \label{a2}
&&[\mathcal{F} \phi^{(\hat{x})}_{y,\eta, \epsilon}](\xi)=0\quad \mbox{if}\quad \xi\notin  (t_{\min} - \hat{x}\cdot y-\epsilon,\; \eta - \hat{x}\cdot y+\epsilon ).
\en If $\epsilon=0$, it holds that
\be \label{a3}
&&[\mathcal{F} \phi^{(\hat{x})}_{y,\eta}](\xi)=
\left\{\begin{array}{lll}
\sqrt{2\pi}/(\eta - t_{\min}) \quad &&\mbox{if}\quad \xi\in (t_{\min} - \hat{x}\cdot y,\; \eta - \hat{x}\cdot y ), \\
0 \quad&&\mbox{if  otherwise}.
\end{array}\right.
\en
\end{lemma}
\begin{proof}
Setting $\xi = t-\hat{x}\cdot z$, we can rewrite the function $\phi^{(\hat{x})}_{y,\eta, \epsilon}$ as
\ben
\phi^{(\hat{x})}_{y,\eta, \epsilon}(\tau)=\int_{\R} e^{{\rm i} \tau\xi}g_{\eta,\epsilon}(\xi, \hat{x})\,d\xi,
\quad
g_{\eta,\epsilon}(\xi,\hat{x})=\frac{1}{(\eta - t_{\min})|B_{\epsilon}(y)| }\int_{t_{\min} - \xi}^{\eta - \xi} \int_{\Gamma(t, \hat{x})} ds(z)dt,
\enn
with $\G(t, \hat{x})=\{z\in B_\epsilon(y): \hat{x}\cdot z=t\}$. Hence, $\mathcal{F} \phi^{(\hat{x})}_{y,\eta, \epsilon}=\sqrt{2\pi}\;g_{\eta,\epsilon}(\cdot,\hat{x})$.
Observing that
$$\sup\big(\hat{x}\cdot  B_\epsilon(y)\big)=\hat{x}\cdot y+\epsilon,\quad  \inf\big(\hat{x}\cdot  B_\epsilon(y)\big)=\hat{x}\cdot y-\epsilon,$$ we obtain $\eqref{a1}$ and $\eqref{a2}$ from the expression of $g_{\eta,\epsilon}(\cdot, \hat{x})$. If $\epsilon=0$, there holds
\ben
\phi^{(\hat{x})}_{y,\eta}(k)=\frac{1}{\eta - t_{\min} }\int_{t_{\min}}^{\eta } e^{{\rm i} k  (t-\hat{x}\cdot y)}\,dt
=\frac{1}{\eta - t_{\min} }\int_{\R} e^{{\rm i} k \xi} g_{\eta}(\xi)\,d\xi,
\enn
where
\ben
g_{\eta}(\xi):=\left\{\begin{array}{lll}
1,\quad &&\mbox{if} \quad \xi\in (t_{\min} - \hat{x}\cdot y,\; \eta - \hat{x}\cdot y), \\
0, \quad&&\mbox{if\; otherwise}.
\end{array}\right.
\enn
Therefore, $[\mathcal{F} \phi^{(\hat{x})}_{y,\eta}](\xi)=\sqrt{2\pi}\,g_{\eta}(\xi)/(\eta - t_{\min}) $.
\end{proof}
Define the unbounded and parallel strips (see Figure \ref{strip})
\be\label{K}
K_D^{(\hat{x})}&:=&\{y\in \R^3: \inf(\hat{x}\cdot D) < \hat{x}\cdot y < \sup (\hat{x}\cdot D) \}\subset \R^3,
\\ \label{tildeK}
K _{D,\eta}^{(\hat{x})}&:=&\{y\in \R^3: \inf(\hat{x}\cdot D)- t_{\max} + \eta < \hat{x}\cdot y < \sup (\hat{x}\cdot D) \}\subset \R^3,
\en
whose directions are perpendicular to the observation direction $\hat{x}$.
The region $K_D^{(\hat{x})}\subset \R^3$ represents the smallest strip containing $D$ and perpendicular to the vector $\hat{x}\in\s^2$. By the definition \eqref{tildeK}, it is obvious that
\ben
K _{D,\eta}^{(-\hat{x})}&:=&\{y\in \R^3: \inf(\hat{x}\cdot D)< \hat{x}\cdot y < \sup (\hat{x}\cdot D)+t_{\max}-\eta \}\subset \R^3,
\enn
and that $K _{D}^{(\hat{x})}=
K_{D,\eta}^{(\hat{x})}\cap K_{D,\eta}^{(-\hat{x})}  $ if $\eta\in(t_{\min}, t_{\max}]$.%
\begin{lemma}\label{lem3.4} \begin{itemize}
\item[(i)] For $y\in K_D^{(\hat{x})}$,  there exists $\epsilon_0=\epsilon_0(y)>0$ such that  $\phi^{(\hat{x})}_{y,\eta, \epsilon}\in \range (L_D^{(\hat{x})})$ for all $\epsilon\in(0, \epsilon_0)$ and $\eta \in (t_{\min},t_{\max}]$.
\item[(ii)] If $y\notin \overline{K _{D,\eta}^{(\hat{x})}} $, we have
 $\phi^{(\hat{x})}_{y,\eta, \epsilon}\notin \range (L_D^{(\hat{x})})$ for all $\epsilon>0$ and $\eta \in (t_{\min},t_{\max}]$.
\end{itemize}
\end{lemma}
\begin{proof}
($i$) If $y\in K_D^{(\hat{x})}$, there must exist some $z\in D$ and $\epsilon_0>0$ such that  $\hat{x}\cdot y=\hat{x}\cdot z$ and  $B_\epsilon(z)\subset D$ for all $\epsilon\in(0, \epsilon_0)$.  Moreover, we have $\phi^{(\hat{x})}_{y,\eta, \epsilon}=\phi^{(\hat{x})}_{z,\eta, \epsilon}$. Set
\ben
u(x, t):=\left\{\begin{array}{lll}
\frac{1}{|B_\epsilon(z)|\;(\eta-t_{\min})} ,  \quad &&\mbox{if}\; x\in B_\epsilon(z),  t\in [t_{\min}, \eta],\\
0\quad&&\mbox{if otherwise}.
\end{array}\right.
\enn
It is obvious  that $u(x,t)\in L^2(D\times (t_{\min}, t_{\max}))$ for $t_{\min}< \eta \leq t_{\max}$. By the definition of $L_D^{(\hat{x})}$ (see \eqref{def:L}), it is easy to see
 $\phi^{(\hat{x})}_{z,\eta, \epsilon}=L_D^{(\hat{x})}u$.

($ii$) Given $y\notin \overline{K _{D,\eta}^{(\hat{x})} }$, we suppose on the contrary that
$\phi^{(\hat{x})}_{y,\eta, \epsilon}= L_D^{(\hat{x})}g$
with some $g\in L^2(D\times(t_{\min}, t_{\max}))$, i.e.,
\begin{equation}\label{eq}
\phi^{(\hat{x})}_{y,\eta, \epsilon}(\tau)=\int_{t_{\min}}^{t_{\max}}\int_{D} e^{{\rm i} \tau (t-\hat{x}\cdot z)} g(z,t)\,dz dt,\qquad \tau\in [0, K].
\end{equation}
By the analyticity in $\tau$, the above relation remains valid for all $\tau\in \R$. Hence, the supporting intervals of the Fourier transform of both sides of \eqref{eq} must coincide.
Using \eqref{a1} and \eqref{rl-f}, we obtain
\begin{equation*}
\left[t_{\min} - \hat{x}\cdot y-\epsilon,\; \eta - \hat{x}\cdot y+\epsilon \right]\, \subset\,\left[ t_{\min} - \sup(\hat{x}\cdot D),\; t_{\max} - \inf(\hat{x}\cdot D)  \right],
\end{equation*}
leading to
\begin{equation}\label{eq:ep}
\inf(\hat{x}\cdot D)+\epsilon-t_{\max} +\eta \leq \hat{x}\cdot y \leq \sup(\hat{x}\cdot D) -\epsilon,\quad\mbox{for all}\; \epsilon>0.
\end{equation}
This implies that $y\in \overline{K _{D,\eta}^{(\hat{x})} }$,
a contradiction to the assumption $y\notin \overline{K _{D,\eta}^{(\hat{x})} }$. This proves  $\phi^{(\hat{x})}_{y,\eta, \epsilon}\notin \range (L_D^{(\hat{x})})$ for all $\epsilon>0$.
\end{proof}

In contrast with our previous work \cite{GGH2022},  Lemma \ref{lem3.4} cannot be used to characterize the strip $K_D^{(\hat{x})}$ when $t_{\max}$ is not given, due to the lack of a inclusion relation between $\phi^{(\hat{x})}_{y,\eta, \epsilon}$ and $\range (L_D^{(\hat{x})})$ for $y\in K _{D,\eta}^{(\hat{x})}\backslash \overline{K _D^{(\hat{x})}}$.
The price we pay is to make use of the multi-frequency data at two opposite directions, which will be discussed in details in Section \ref{sec:3.3}.
We proceed with the definition of the region $\Omega_\eta$ with $\eta\in(t_{\min}, t_{\max}]$, given by
\be\label{Om}
\Omega_\eta=\{y\in \R^3: \mbox{there exists an $\epsilon_0(y)>0$ such that}\;\phi_{y,\eta,\epsilon}^{(\hat{x})}\in \mbox{Range}(L_D^{(\hat{x})})\;\mbox{for all}\; \epsilon\in(0, \epsilon_0)\}.
\en
It obviously follows from Lemma \ref{lem3.4} that $K_D^{(\hat{x})} \subset \Omega_\eta \subset K _{D,\eta}^{(\hat{x})}$ if $\eta\in(t_{\min}, t_{\max}]$.
Later we shall define an indicator function for imaging the region $\Omega_\eta$ from the multi-frequency far-field data at a single observation direction. This implies that, when $\hat{x}=(1,0)$, the right boundary of $\Omega_\eta$ given by $\{y\in \R^3: \hat{x}\cdot y= \sup(\hat{x}\cdot D)\}$ (see Fig. \ref{strip}) illustrates partial information  on the boundary of the source support.
 Moreover, we have $K_D^{(\hat{x})} = \Omega_\eta = K _{D,\eta}^{(\hat{x})}$ if and only if $\eta = t_{\max}.$

\subsection{The initial moment $t_{\min}$ is unknown}
	When the initial time point $t_{\min}$ is unknown, the test function $\phi^{(\hat{x})}_{y,\eta, \epsilon}$ can be correspondingly defined as
	\begin{equation}
		\phi^{(\hat{x})}_{y,\eta, \epsilon}(k):=\frac{1}{(t_{\max} - \eta) \,|B_\epsilon(y)|}\int_{\eta}^{t_{\max}}\int_{B_\epsilon(y)}e^{{\rm i} k  (t-\hat{x}\cdot z)} dzdt,\qquad k\in [0, K],
	\end{equation}
	where $\eta<t_{\max}$ is a parameter lying in a small neighborhood on the left of the terminal time point $t_{\max}$. As $\epsilon\rightarrow 0$,  there holds the convergence
	\begin{equation}
		\phi^{(\hat{x})}_{y,\eta, \epsilon}(k)\rightarrow \phi^{(\hat{x})}_{y,\eta}(k):=\frac{1}{t_{\max}-\eta}\int_{\eta}^{t_{\max}}e^{ik(t-\hat{x}\cdot y)} dt.
	\end{equation}
In this case the definition of $K_{D,\eta}^{(\hat{x})}$ will be changed into
	\begin{equation}
		K_{D,\eta}^{(\hat{x})}:=\{y\in \R^3: \inf(\hat{x}\cdot D) < \hat{x}\cdot y < \sup (\hat{x}\cdot D) - t_{\min} + \eta\}\subset \R^3.
	\end{equation}
The strip $K _D^{(\hat{x})}$ is still a  subset of
$K_{D,\eta}^{(\hat{x})}$ if $\eta\in [t_{\min}, t_{\max})$.	
Analogously to Lemma \ref{lem3.4} we have
\begin{lemma}\label{cor3.3}
\begin{itemize}
\item[(i)] For $y\in K_D^{(\hat{x})}$,  there exists an $\epsilon_0=\epsilon_0(y)>0$ such that  $\phi^{(\hat{x})}_{y,\eta, \epsilon}\in \range (L_D^{(\hat{x})})$ for all $\epsilon\in(0, \epsilon_0)$ and $\eta \in [t_{\min},t_{\max})$.
\item[(ii)] If $y\notin \overline{K _{D,\eta}^{(\hat{x})} }$, we have
 $\phi^{(\hat{x})}_{y,\eta, \epsilon}\notin \range (L_D^{(\hat{x})})$ for all $\epsilon>0$ and $\eta \in [t_{\min},t_{\max})$.
\end{itemize}
\end{lemma}
The region $\Omega_\eta$ defined by \eqref{Om} still satisfies the inclusion relations $K_D^{(\hat{x})} \subset \Omega_\eta \subset K _{D,\eta}^{(\hat{x})}$ if $\eta\in [t_{\min}, t_{\max})$. In the case that $\hat{x}=(1,0)$, the left boundary of $\Omega_\eta$, which coincides with $\{y\in \R^3: \hat{x}\cdot y= \inf(\hat{x}\cdot D)\}$, yields information on the support of the source function.

\subsection{Indicator functions and inversion algorithms}\label{sec:3.3}
\subsubsection{Uniqueness and algorithm at two opposite observation directions}
We first recall from the previous two subsections that the test function $\phi^{(\hat{x})}_{y,\eta, \epsilon}$ with a small $\epsilon>0$ can be utilized to characterize the region $\Omega_\eta$. Define the auxiliary indicator function

\be\label{indicator4}
I_{\eta,\epsilon}^{(\hat{x})}(y):=\sum_{n=1}^\infty\frac{|\langle \phi^{(\hat{x})}_{y,\eta, \epsilon}, \psi_n^{(\hat{x})} \rangle|_{L^2(0, K )}^2}{ |\lambda_n^{(\hat{x})}|}, \qquad y\in \R^3.
\en
By Picard's theorem, for every $\eta \in (t_{\min},t_{\max})$ we have
$I_{\eta,\epsilon}^{(\hat{x})}(y) < +\infty$ if $y\in \Omega_\eta$ and $\epsilon\in(0, \epsilon_0(y))$. However, the region $\Omega_\eta$ provides only partial information on the strip $K_D^{(\hat{x})}$. For the purpose of completely characterizing $K_D^{(\hat{x})}$, we utilize the multi-frequency data at two opposite directions $\hat{x}$ and $-\hat{x}$
to define the indicator function
\begin{equation}\label{w-eps}
	W_{\epsilon}^{(\hat{x})} (y) = \left[I_{\eta,\epsilon}^{(\hat{x})}(y) +I_{\eta,\epsilon}^{(-\hat{x})}(y) \right]^{-1}.
\end{equation}
Combining the range identity \eqref{RI} and Lemma \ref{lem3.4} we can characterize the strip $K_D^{(\hat{x})}$ through the indicator function \eqref{w-eps}.
\begin{theorem}[Determination of the strip $K_D^{(\hat{x})}$]\label{Th:factorization}

\begin{itemize}
\item[(i)] If $y\in K_D^{(\hat{x})}$, there exists an $\epsilon_0=\epsilon_0(y)>0$ such that
$W_{\epsilon}^{(\hat{x})}(y)$ is strictly positive for all $\epsilon\in(0, \epsilon_0)$.
\item[(ii)] If $y\notin \overline{K_D^{(\hat{x})} }$,  there holds
$W_\epsilon^{(\hat{x})}(y)=0$ for all $\epsilon>0$.
\end{itemize}
\end{theorem}

\begin{proof}
(i) Given $y\in K_D^{(\hat{x})}$, by Lemmas \ref{lem3.4} and \ref{cor3.3} there exists an $\epsilon_0(y)>0$ such that
$$\phi^{(\hat{x})}_{y,\eta, \epsilon}\in \range (L_D^{(\hat{x})}), \; \phi^{(-\hat{x})}_{y,\eta, \epsilon}\in \range (L_D^{(-\hat{x})})\quad\mbox{for all}\;\epsilon\in(0, \epsilon_0),\quad \eta \in (t_{\min},t_{\max}).$$
From the range identity \eqref{RI}, we know $0<I_{\eta,\epsilon}^{(\hat{x})}(y)<+\infty$ and $0<I_{\eta,\epsilon}^{(-\hat{x})}(y)<+\infty$. Thus, $W_{\epsilon}^{(\hat{x})}(y)$ must be a finite positive number.

(ii) For $y\notin \overline{K_D^{(\hat{x})} }$, it holds that either
$y\notin \overline{K _{D,\eta}^{(\hat{x})} }$ or
$y\notin \overline{K _{D,\eta}^{(-\hat{x})} }$.
Without loss of generality we assume the former case $y\notin \overline{K _{D,\eta}^{(\hat{x})} }$, which implies $y\notin \Omega_\eta$.  In this case, $\phi^{(\hat{x})}_{y,\eta, \epsilon}\notin \range (L_D^{(\hat{x})})$ and thus $I_{\eta,\epsilon}^{(\hat{x})}(y)=+\infty$ for all $\epsilon>0$. Consequently, $W_\epsilon^{(\hat{x})}(y)=0$ for all $\epsilon>0$.
\end{proof}

Since $\phi^{(\hat{x})}_{y,\eta, \epsilon}$ convergences uniformly to $\phi^{(\hat{x})}_{y,\eta}$ over the  finite wavenumber interval $[k_{\min}, k_{\max}]$, we shall use the limiting function $\phi^{(\hat{x})}_{y,\eta}$ in place of $\phi^{(\hat{x})}_{y, \eta,\epsilon}$ in the aforementioned indicator function. Consequently, we introduce a new indicator function
\begin{equation}\label{indicator1}
W^{(\hat{x})}(y):=\left[I_{\eta,0}^{(\hat{x})}(y) +I_{\eta,0}^{(-\hat{x})}(y) \right]^{-1}\sim
\left[\sum_{n=1}^N\frac{|\langle \phi^{(\hat{x})}_{y,\eta}, \psi_n^{(\hat{x})} \rangle|_{L^2(0, K )}^2}{ |\lambda_n^{(\hat{x})}|} + \frac{|\langle \phi^{(-\hat{x})}_{y,\eta}, \psi_n^{(-\hat{x})} \rangle|_{L^2(0, K )}^2}{ |\lambda_n^{(-\hat{x})}|} \right]^{-1},
\end{equation}
where $y\in \R^3$ and the integer $N\in \N$ is a truncation number.  Taking the limit $\epsilon\rightarrow 0$ in Theorem \ref{Th:factorization}, it follows that
\be\label{W}
W^{(\hat{x})}(y)=\left\{\begin{array}{lll}
\geq 0\quad&&\mbox{if}\quad y\in K_D^{(\hat{x})},\\
0\quad&&\mbox{if}\quad y\notin \overline{K_D^{(\hat{x})}}.
\end{array}\right.
\en
Hence, the values of $W^{(\hat{x})}$ in the strip \rot{$K_D^{(\hat{x})}$} should be relatively bigger than those elsewhere.

Uniqueness results in identifying the strip $K_D^{(\hat{x})}$ and the terminal time point $t_{\max}$ when $t_{\min}$ is given are summarized as follows.
\begin{theorem}\label{Th:3.1}[Uniqueness at two opposite directions]
Assume the terminal moment $t_{\max}$ is unknown and let $\hat{x} \in \s^2$ be an arbitrarily fixed observation direction. Then the strip $K_D^{(\hat{x})}$ and the terminal time point $t_{\max}$ can be uniquely determined by the multi-frequency far-field data $\{w^\infty(\pm\hat{x}, k): k\in(k_{\min}, k_{\max})\}$.
\end{theorem}
\begin{proof}
The unique determination of the strip $K_D^{(\hat{x})}$ follows from Theorem \ref{Th:factorization}. Below we shall prove the uniqueness in identifying $t_{\max}$ by contradiction. Suppose that there are two different terminal time points $t_{\max}^{(1)}<t_{\max}^{(2)}$ but corresponding to identical multi-frequency far-field data at the observation directions $\pm \hat{x}$. From the expressions
\ben
w^\infty(\pm\hat{x},k)=\frac{1}{\sqrt{2\pi}}\int_{t_{\min}}^{t^{(j)}_{\max}}\int_D e^{{\rm i} k (t\mp\hat{x}\cdot y)} S(y,t)\,dy \,dt,\quad j=1,2,
\enn it is not difficult to deduce that
\ben
0=\int_{t^{(1)}_{\max}}^{t^{(2)}_{\max}}e^{{\rm i} k t}\left(\int_D  S(y,t)\,dy\right) \,dt=0\quad\mbox{for all}\quad k\in(k_{\min}, k_{\max}).
\enn
By the analyticity on $k$, the above identity holds true for all $k\in \R$. Hence
\ben
\int_D S(y,t)\,dy=0 \quad\mbox{for all}\quad t\in (t_{\max}^{(1)},\;t_{\max}^{(2)}),
\enn
which contradicts the positivity condition \eqref{F} for $S$. This proves $t_{\max}^{(1)}=t_{\max}^{(2)}$.
\end{proof}

Having determined the strip $K_D^{(\hat{x})}$ from Theorem \ref{Th:factorization}, one can recover
the unknown terminal moment $t_{\max}$ by plotting the one-dimensional function $\eta\mapsto
I_{\eta,\epsilon}^{(\hat{x})}(z)$ with some fixed $z\in K_D^{(\hat{x})}$.
\begin{theorem}[Determination of $t_{\max}$]\label{Th:max} Suppose that the initial time point $t_{\min}$ is known.
Choose $z\in K_D^{(\hat{x})}$ such that $\hat{x}\cdot z-\inf(\hat{x}\cdot D)=\epsilon_0$ for some small number $\epsilon_0>0$. Then
\ben
[I_{\eta,\epsilon}^{(\hat{x})}(z)]^{-1}=\left\{\begin{array}{lll}
\mbox{a finite positive number}, &&\mbox{if}\quad\eta\in(t_{\min},t_{\max}],\;\epsilon\in(0, \epsilon_0/2],\\
0, && \mbox{if}\quad\eta>t_{\max}+\epsilon_0,\; \epsilon\in(0, \epsilon_0/2].
\end{array}\right.
\enn
\end{theorem}
\begin{proof}
By the proof of Lemma \ref{lem3.4} (i), it holds that $\phi^{(\hat{x})}_{y,\eta, \epsilon}\in \range (L_D^{(\hat{x})})$ for all $\epsilon\in(0, \epsilon_0/2]$ and $\eta \in (t_{\min},t_{\max}]$. Hence, $0<I_{\eta,\epsilon}^{(\hat{x})}(z)<\infty$
and $[I_{\eta,\epsilon}^{(\hat{x})}(z)]^{-1}>0$
 for such $\eta$ and $\epsilon$. On the other hand,
one can repeat the arguments in the proof of Lemma \ref{lem3.4} (ii) to prove that $\phi^{(\hat{x})}_{z,\eta, \epsilon}\notin \range (L_D^{(\hat{x})})$ for all $\epsilon\in(0,\epsilon_0/2]$ if $\eta>t_{\max}+\epsilon_0$, which together with the Picard's theorem implies
 $I_{\eta,\epsilon}^{(\hat{x})}(z)=+\infty$ and
 $[I_{\eta,\epsilon}^{(\hat{x})}(z)]^{-1}=0$.
 \end{proof}
\begin{remark} (i)
Choosing $\epsilon_0>0$ to be sufficiently small, one deduces  from Theorem \ref{Th:max} that the function  $\eta\mapsto
[I_{\eta,\epsilon}^{(\hat{x})}(z)]^{-1}$ must decay fast in a small neighborhood on the right hand side of $\eta=t_{\max}$, which indicates an approximation of $t_{\max}$. We refer to the Section \ref{sub-tmax} for numerical examples.
(ii) The results of Theorems \ref{Th:3.1} and \ref{Th:max} can be
established analogously when $t_{\min}$ is unknown and $t_{\max}$ is given.
\end{remark}

\subsubsection{Algorithm and uniqueness at a finite pair of opposite  directions}
If the wave signals are detected at a finite number of observation directions $\{\pm \hat{x}_j: j=1,2,\cdots, M\}$, we shall make use of the following indicator function:
\begin{equation}\label{W2}
W(y)\rot{= \left[\sum_{j=1}^M \frac{1}{W^{(\hat{x}_j)}(y) }\right]^{-1}}= \left[\sum_{j=1}^M\sum_{n=1}^N\frac{|\langle \phi^{(\hat{x}_j)}_{y,\eta}, \psi_n^{(\hat{x}_j)} \rangle|_{L^2(0, K )}^2}{ |\lambda_n^{(\hat{x}_j)}|} + \frac{|\langle \phi^{(-\hat{x}_j)}_{y,\eta}, \psi_n^{(-\hat{x}_j)} \rangle|_{L^2(0, K )}^2}{ |\lambda_n^{(-\hat{x}_j)}|} \right]^{-1}, y\in \R^3.
\end{equation}
Define the $\Theta$-convex hull of $D$ associated with the directions $\{\pm\hat{x}_j: j=1,2,\cdots, M\}$ as
\ben
\Theta_D:=\bigcap_{j=1,2,\cdots, M} K_D^{(\hat{x}_j)}.
\enn
By Lemma \ref{lem3.4}, the above convex hull can be uniquely determined by the multi-frequency far-field data
$\{w^\infty(\pm\hat{x}_j, k): k\in [k_{\min}, k_{\max}], j=1,2,\cdots, M\}$.

\begin{theorem}[Algorithm with a finite pair of opposite directions]\label{TH-hull}
We have $W(y) \geq 0$ if $y\in \Theta_D$ and $W(y)=0$ if
$y\notin \overline{\Theta_D}$.
\end{theorem}
\begin{proof}
If $y\in \Theta_D$, then $y\in K_D^{(\hat x_j)}$ for all $j=1,2,...,M$, yielding that $\hat x_j \cdot y \in \hat x_j \cdot D$. Hence, one deduces from  Theorem \ref{Th:factorization} that  $0\leq W^{(\hat x_j)}(y)<\infty$ for all $j=1,2,...,M$, implying that $W(y) \geq 0$. On the other hand, if $y\notin \overline{\Theta_D}$, there must exist some unit vector $\hat{x}_l$ such that $y\notin \overline{K_D^{(\hat x_l)}}$. Again using Theorem \ref{Th:factorization}, we get $$
[W^{(\hat{x}_l)}(y)]^{-1}=I_{\eta}^{(\hat{x}_l)}(y) +I_{\eta}^{(-\hat{x}_l)}(y)=\infty,$$ which proves $W(y)=0$ for $y\notin \Theta_D$.
\end{proof}
The values of $W(y)$ are expected to be large for $y\in \tcr{\Theta_D}$ and small for those $y\notin \Theta_D$. 
As a by-product of the above factorization method, we obtain a uniqueness result with multi-frequency far-field data for all observation directions. Denote by $\mbox{ch}(D)$ the convex hull of $D$, that is,  the intersections of all half spaces containing $D$.

\begin{theorem}[Uniqueness with all observation directions]\label{TH-u}
Let the assumption \eqref{F} hold. Then $\mbox{ch}(D)$ can be uniquely determined by the multi-frequency far-field patterns $\{w^\infty(\hat{x}, k):  k\in[k_{\min}, k_{\max}],\, \hat{x}\in \s^2\}.$
\end{theorem}
\begin{proof}
Given a fixed direction $\hat{x}\in \s^2$, it follows from Lemma \ref{lem3.4} that the strip $K_D^{(\hat{x})}$ can be uniquely determined. Since the $\Theta$-convex hull of $D$ associated with all directions $\hat{x}\in \s^2$ coincides with $\mbox{ch}(D)$, we obtain uniqueness in recovering the convex hull of $D$.
\end{proof}

\section{Inversion algorithm with multi-frequency near-field data}\label{sec4}



\begin{figure}[H]
		\centering
		\scalebox{0.7}{
		\begin{tikzpicture}

		\node (x) at (2,0) {};
		\filldraw[green!20,even odd rule](x)circle(3.8)(x)circle(1);
		\filldraw[blue!20,even odd rule](x)circle(4.8)(x)circle(3.8);

		\draw (0.5,-2) node [right] {$A_{D}^{({x})}$};
		\draw (-0.8,-2) node [left] {$A_{D,\eta}^{({x})}$};
	
		\draw[domain = -2:360,very thick][samples = 200] plot({cos(\x)+0.65*cos(2* \x)-0.65},{1.5*sin(\x)});
		\draw (x) node [right]{$x$};
		\fill (x)	 circle (2pt);	
		\draw [dotted] let \p1 = ($ (x) - (1,0) $),  \n2 = {veclen(\x1,\y1)}
		in (x) circle (\n2);
		\draw [dotted] let \p1 = ($ (x) - (-1.8,0) $),  \n2 = {veclen(\x1,\y1)}
		in (x) circle (\n2);
		\draw [dotted] let \p1 = ($ (x) - (-2.8,0) $),  \n2 = {veclen(\x1,\y1)}
		in (x) circle (\n2);

		\draw (2,0.8) node [above] {$\inf\limits_{z\in D} |x-z|$};
		\draw (1.5,2.8) node [above] {$\sup\limits_{z\in D} |x-z|$};
		\draw (1.5,5.5) node [below] {$\sup\limits_{z\in D} |x-z|+t_{\max}-\eta$};
		\draw (-0.5,-0.4) node [below] {$D$};

		
		\end{tikzpicture}
		}
		\caption{Illustration of the annuluses $A_D^{({x})}$ (green area) and $A_{D,\eta}^{({x})}$ (union of green and blue area) defined in \eqref{annulus-a} and \eqref{annulus-b}. }
		\label{annuluses}
	\end{figure}
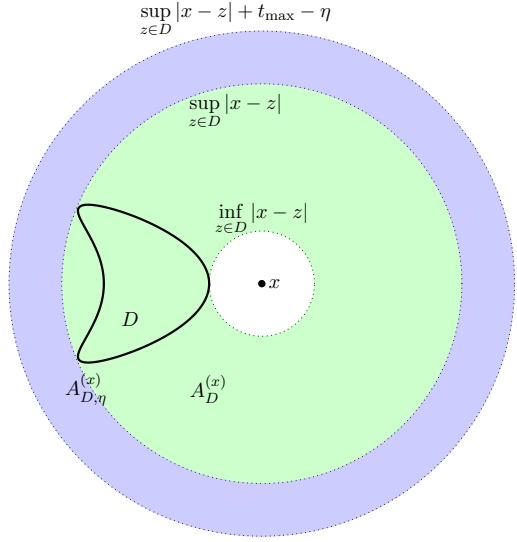

Suppose that the initial time point $t_{\min}$ is given, but the terminal time point $t_{\max}$ is unknown. Choose the test functions
\begin{equation*}
\tilde{\phi}^{(x)}_{y,\eta,\epsilon}(k) =\frac{1}{(\eta -t_{\min}) \,|B_\epsilon(y)|}\int_{t_{\min}}^{\eta}\int_{B_\epsilon(y)} e^{ik(t+|x-z|)} dzdt,\quad k\in [0,K],\eta\in (t_{\min},t_{\max}].
\end{equation*}
As $\epsilon\rightarrow 0$,  there holds the convergence
\begin{equation}
	\tilde{\phi}^{(x)}_{y,\eta,\epsilon}(k) \to \tilde{\phi}^{(x)}_{y,\eta}(k) := \frac{1}{\eta -t_{\min}} \int_{t_{\min}}^{\eta} e^{ik(t+|x-y|)} dt,\quad k\in [0,K],\eta\in (t_{\min},t_{\max}].
\end{equation}
Introduce two annuluses centered at the receiver $x\in \partial B_R$ (see Figure \ref{annuluses}):
\begin{equation} \label{annulus-a}
	A^{(x)}_D := \{y\in \R^3: \inf\limits_{z\in D}|x-z| < |x-y| < \sup\limits_{z\in D}|x-z| \}
\end{equation}
and
\begin{equation}\label{annulus-b}
	A^{(x)}_{D,\eta} := \{y\in \R^3: \inf\limits_{z\in D}|x-z| < |x-y| < \sup\limits_{z\in D}|x-z| +t_{\max}-\eta \}.
\end{equation}
Let the operator $\tilde{L}=\tilde{L}^{(x)}_D$ be defined as in Lemma \ref{Fac-N}.
By arguing similarly to Lemma \ref{lem3.4}, 
we deduce from the range identity \eqref{RN} that
\begin{theorem}\label{lem4.1}
 Suppose that the initial time point $t_{\min}$ is known, $D\subset B_R$ and $|x|=R$.
 \begin{itemize}
\item[(i)] For $y\in A_D^{(x)}$,  there exists $\epsilon_0>0$ such that  $\tilde{\phi}^{({x})}_{y,\eta, \epsilon}\in \range (\tilde{L}_D^{({x})})$ for all $\epsilon\in(0, \epsilon_0)$ and $\eta \in (t_{\min},t_{\max}]$.
\item[(ii)] If $y\notin \overline{A_{D,\eta}^{({x})} }$, we have
 $\tilde{\phi}^{({x})}_{y,\eta, \epsilon}\notin \range (\tilde{L}_D^{({x})})$ for all $\epsilon>0$ and $\eta \in (t_{\min},t_{\max}]$.


\end{itemize}
\end{theorem}

For $y\in \R^3$, introduce the auxiliary function
\begin{equation}  \label{indic-near-1}
\widetilde{W}^{(x)}(y):=\left[\sum_{n=1}^N\frac{|\langle \tilde{\phi}^{(x)}_{y,\eta}, \tilde{\psi}_n^{(x)} \rangle|_{L^2(0, K )}^2}{ |\tilde{\lambda}_n^{(x)}|}
+
\frac{|\langle \tilde{\phi}^{(-x)}_{y,\eta}, \tilde{\psi}_n^{(-x)} \rangle|_{L^2(0, K )}^2}{ |\tilde{\lambda}_n^{(-x)}|}
 \right]^{-1},
\end{equation}
where $(\tilde{\lambda}_n^{(x)}, \tilde{\psi}_n^{(x)})$ is an eigensystem of the near-field operator $\mathcal{N}^{(x)}_D$ and the integer $N\in \N$ is a truncation number.
Note that the multi-frequency data at two opposite receivers $\pm x\in \partial B_R$ are required in computing the function
$\widetilde{W}^{(x)}$.
As the counterpart to  the relation \eqref{W} and Theorem \ref{Th:max}, one can
 show in the near-field case that
\begin{corollary}\label{coro4.3}
(i) It holds that
\begin{equation*}
\widetilde{W}^{(x)}(y)=\left\{\begin{array}{lll}
\geq 0\quad&&\mbox{if}\quad y\in A_D^{(x)}\cap A_D^{(-x)},\\
0\quad&&\mbox{if}\quad y\notin \overline{A_D^{(x)}\cap A_D^{(-x)}}.
\end{array}\right.
\end{equation*}
(ii) Choose $z_0\in A_D^{(x)}\cap A_D^{(-x)}$ such that $\sup_{z\in D} |z-x|-|x-z_0|=\epsilon_0$ for some small number $\epsilon_0>0$. Then
\ben
[\tilde{I}_{\eta,\epsilon}^{(x)}(z_0)]^{-1}=\left\{\begin{array}{lll}
\mbox{a finite positive number}, &&\mbox{if}\quad\eta\in(t_{\min},t_{\max}],\;\epsilon\in(0, \epsilon_0/2],\\
0, && \mbox{if}\quad\eta>t_{\max}+\epsilon_0,\; \epsilon\in(0, \epsilon_0/2],
\end{array}
\right.
\enn
where
\be   \label{indic-near}
\tilde{I}_{\eta,\epsilon}^{(x)}(y):=
\sum_{n=1}^\infty\frac{|\langle \tilde{\phi}^{(x)}_{y,\eta,\epsilon}, \tilde{\psi}_n^{(x)} \rangle|_{L^2(0, K )}^2}{ |\tilde{\lambda}_n^{(x)}|},\quad y\in \R^3.
\en
\end{corollary}
By the first assertion, the values of $\widetilde{W}^{(\hat{x})}$ in the domain $A_D^{(x)}\cap A_D^{(-x)}$ should be relatively bigger than those elsewhere, which yield information of the source support. The second assertion can be used to calculate the terminal time point $t_{\max}$ by properly choosing $z_0$ with a sufficiently small $\epsilon_0>0$.  If the wave signals are detected at a finite couple of symmetric observation points $\{\pm {x}_j\in S_R: j=1,2,\cdots, M\}$,  the indicator function
\begin{equation} \label{indic-near-m}
\widetilde{W}(y)\rot{= \left[\sum_{j=1}^M \frac{1}{\widetilde{W}^{({x}_j)}(y) }\right]^{-1}}= \left[\sum_{j=1}^M\sum_{n=1}^N\frac{|\langle \tilde{\phi}^{({x}_j)}_{y,\eta}, \tilde{\psi}_n^{({x}_j)} \rangle|_{L^2(0, K )}^2}{ |\tilde{\lambda}_n^{({x}_j)}|}  \right]^{-1},\qquad y\in \R^3,
\end{equation}
can be used to approximate the region $\bigcap_{j=1,2,\cdots,M} \big\{A_D^{(x_j)}\cap A_D^{(-x_j)}\big\}$. As the radius of the receivers $|x_j|=R\rightarrow\infty$,
this region will converge to the $\Theta$-convex set associated with the directions $\pm \hat{x}_j$ for $j=1,2,\cdots,M$.

\begin{remark}
\begin{itemize}
\item[(i)] In the case that
the initial moment $t_{\min}$ is unknown,
we define the test function
\begin{equation*}
\tilde{\phi}^{(x)}_{y,\eta,\epsilon}(k) =\frac{1}{(t_{\max}-\eta) \,|B_\epsilon(y)|}\int^{t_{\max}}_{\eta}\int_{B_\epsilon(y)} e^{ik(t+|x-z|)} dzdt,\quad k\in [0,K],\eta\in [t_{\min},t_{\max}),
\end{equation*}
which converges to
\begin{equation}
	\tilde{\phi}^{(x)}_{y,\eta,\epsilon}(k) \to \tilde{\phi}^{(x)}_{y,\eta}(k) := \frac{1}{t_{\max} - \eta } \int^{t_{\max}}_{\eta} e^{ik(t+|x-z|)} dt
	\end{equation}
as $\epsilon\rightarrow 0$.
We retain the notation $A^{(x)}_D$ defined as in \eqref{annulus-a}, but change the definition of $A^{(x)}_{D,\eta}$ into
\begin{equation*}
	A^{(x)}_{D,\eta} := \{y\in \R^3: \inf\limits_{z\in D}|x-z| +t_{\min}-\eta < |x-y| < \sup\limits_{z\in D}|x-z|  \}.
\end{equation*}
Then one can prove the results in Theorem \ref{lem4.1} and Corollary \ref{coro4.3} in the same manner.
\item[(ii)] In contrast with Theorem \ref{TH-u}, we cannot obtain the uniqueness result in recovering $\mbox{ch}(D)$ even if the near-field data are measured at all receivers lying on $|x|=R$.
 \end{itemize}
\end{remark}

\section{Discussions on source functions with two disconnected supports}\label{mul-compo}
In the previous sections the support $D$ of the source function $S$ is always supposed to be connected. In this section we assume that $D=D_1\cup D_2\subset \mathbb R^3$ contains two disjoint sub-domains $D_j$ ($j=1,2$) which can be separated by some plane.
For simplicity we only consider the far-field measurement data at multi-frequencies.
The indicator function \eqref{W} can be used to image the strip $K_D^{(\hat{x})}$ if the multi-frequency data are observed at two opposite directions $\pm\hat{x}$. With all observation directions the convex hull of $D$  can be recovered from the indicator \eqref{W2}. Physically, it would be more interesting to determine $\mbox{ch}(D_j)$ for each $j=1,2$, whose union is usually only a subset of $\mbox{ch}(D)$.
Analogously to Lemma \ref{lem3.4} and Theorem \ref{Th:factorization},  we can prove the following results.
\begin{corollary} Let $\hat{x}\in \s^2$ be fixed and suppose one of $t_{\min}$ and $t_{\max}$ is unknown. Set $T:=t_{\max}-t_{\min}>0$.
\begin{itemize}
\item[(i)] For $y\in{K_{D_1}^{(\hat{x})}}\cup {K_{D_2}^{(\hat{x})}} $,  we have $\phi^{(\hat{x})}_{y,\eta,\epsilon}\in \range (L_D^{(\hat{x})})$ for all $\epsilon\in(0, \epsilon_0)$ with some $\epsilon_0>0$.
\item[(ii)] If $y\notin \overline{{K_{D_1,\eta}^{(\hat{x})}}}\cup \overline{{K_{D_2,\eta}^{(\hat{x})}}} $, we have
$\phi^{(\hat{x})}_{y,\eta,\epsilon}\notin \range (L_D^{(\hat{x})})$ for all $\epsilon>0$,
 provided one of the following conditions holds
\be\label{ab}
\quad\quad \quad\;(a)\; \inf (\hat{x}\cdot D_2)-\sup (\hat{x}\cdot D_1)>T;\; (b) \; \inf (\hat{x}\cdot D_1)-\sup (\hat{x}\cdot D_2)>T .
\en
\item[(iii)] Let the indicator function $W^{(\hat{x})}$ be defined by \eqref{indicator1}. Under one of the conditions in \eqref{ab} it holds that
\ben
W^{(\hat{x})}(y)=\left\{\begin{array}{lll}
\geq 0\quad&&\mbox{if}\quad y\in  {K_{D_1}^{(\hat{x})}}\cup {K_{D_2}^{(\hat{x})}} ,\\
0\quad&&\mbox{if}\quad y\notin \overline{{K_{D_1}^{(\hat{x})}}\cup {K_{D_2}^{(\hat{x})}} }.
\end{array}\right.
\enn
\end{itemize}
\end{corollary}
Note that the conditions in \eqref{ab} imply that $(\hat{x}\cdot D_1)\cap (\hat{x}\cdot D_2)=\emptyset$. Furthermore, the Fourier transform of $L_D^{(\hat{x})} f$ with $f\in X_D$ is supported in the following two disjoint intervals (see Lemma \ref{lem1} below for a detailed proof)
\ben
\Big[t_{\min} - \sup(\hat{x}\cdot D_1),\, t_{\max}-\inf(\hat{x}\cdot D_1)\Big]\; \bigcup \;\Big[t_{\min} - \sup(\hat{x}\cdot D_2),\, t_{\max}-\inf(\hat{x}\cdot D_2)\Big].
\enn
For small $T>0$, there exists at least one observation directions $\hat{x}\in \s^2$ such that the relations in \eqref{ab} hold, because $D_1$ and $D_2$ can be separated by some plane by our assumption. If the conditions \eqref{ab} hold for all observation directions $\hat{x}\in \s^2$, one can make use of the indicator function \eqref{W2} to get an image of the set
 $\bigcap_{j=1,2,\cdots M} \{ K_{D_1}^{(\hat{x}_j)} \cup K_{D_2}^{(\hat{x}_j)}\}$, which is usually larger than $\mbox{ch}(D_1)\cup \mbox{ch}(D_2)$.
This means that our approach can only be used to recover partial information of  $\mbox{ch}(D_j)$,
provided the source radiating period $T$ is sufficiently small in comparison with the distance between $D_1$ and $D_2$.
The numerical experiments performed in Section \ref{num} confirm the above theory; see Figures \ref{fig:5}-\ref{fig:8}.

Physically, the conditions in \eqref{ab} ensure that the time-dependent signals recorded at $\hat{x}$ has two disconnected supports which correspond to the wave fields emitting from $D_1$ and $D_2$, respectively. If one can
split the multi-frequency far-field patterns at a single observation direction, it is still possible to recover $\{ K_{D_1}^{(\hat{x})} \cup K_{D_2}^{(\hat{x})}\}$ even if the conditions in \eqref{ab} cannot be fulfilled. Below we prove that the multi-frequency far-field patterns excited by two disconnected source terms can be split under additional assumptions.

To rigorously formulate the splitting problem, we go back to the Helmholtz equation \eqref{Helmholtz},
where $D=D_1\cup D_2$ contains two disjoint bounded and connected sub-domains $D_1$ and $D_2$, that is,
\ben
\Delta u+k^2 u=-f_1(\cdot,k)-f_2(\cdot,k)\quad \mbox{in}\quad \R^3,
\enn
where $\mbox{supp} f_j(\cdot ,k)=D_j$ for any $k>0$  and
\ben
f_j(x,k):=\frac{1}{\sqrt{2\pi}}\int_{t_{\min}}^{t_{\max}} S_j(x, t) e^{ik t} dt,\quad x\in D_j.
\enn
Let $u_j$ be the unique radiating solution to
\ben
\Delta u_j+k^2 u_j=-f_j(x,k)\qquad\mbox{in}\quad \R^3,\quad j=1,2.
\enn
Denote by $u^\infty(\hat{x}, k), u_j^\infty(\hat{x}, k)$ the far-field patterns of $u$ and $u_j$ at some fixed observation direction $\hat{x}=x/|x|$, respectively.
It is obvious that $u^\infty=u_1^\infty+u_2^\infty$, where
\be\label{FF}
u_j^\infty(\hat{x}, k)=\frac{1}{\sqrt{2\pi}}\int_{t_{\min}}^{t_{\max}}\int_{D_j} S_j(y, t) e^{ik (t-\hat{x}\cdot y)}\,dy dt,\quad k\in \R_+.
\en
The splitting problem in the frequency domain can be formulation as follows: Given a fixed observation direction $\hat{x}\in\s^2$,
split $\{u_j^\infty(\hat{x}, k): k\in \R\}$ from the data $\{u^\infty(\hat{x}, k):k\in \R\}$ for $j=1,2$.

Set  $l_j:=\sup(\hat{x}\cdot D_j)-\inf(\hat{x}\cdot D_j)$ and $\Lambda_j=T+\ell_j$ with $T=t_{\max}-t_{\min}>0$.
We make the following assumptions on the time-dependent source function $S(x,t)$.
\begin{itemize}
\item[(i)] $S(x,t)\geq c_0>0$ for $(x,t)\in \overline{D}\times [t_{\min}, t_{\max}]$.
\item[(ii)] $S$ is analytic on $\overline{D}\times [t_{\min}, t_{\max}]$ and the boundary $\partial D$ is analytic.
\item[(iii)] Either $\inf(\hat{x}\cdot D_1)< \inf(\hat{x}\cdot D_2)$, or $\sup(\hat{x}\cdot D_1)> \sup(\hat{x}\cdot D_2)$.
\end{itemize}
\begin{lemma}\label{lem1} Under the assumption (i), the supporting interval of  $\mathcal{F}(u_j^\infty)$ is $I_j:=[t_{\min} - \sup(\hat{x}\cdot D_j), t_{\max} - \inf(\hat{x}\cdot D_j) ]$ and  the function $t\mapsto (\mathcal{F}u_j^\infty)(t)$ is positive
on $I_j$. Moreover,  $(\mathcal{F}u_j^\infty)(t)$ is  analytic  in $t\in I_j$ under the additional assumption (ii).
\end{lemma}
\begin{proof}
The far-field expression \eqref{FF} can be rewritten as
\ben
u_j^\infty(\hat{x}, k)=\int_{\R}e^{ik\xi}g_j(\xi)\,d\xi,\qquad j=1,2,
\enn
where
\be\label{g}
  g_j(\xi):=\frac{1}{\sqrt{2\pi}}\int_{t_{\min}}^{t_{\max}} \int_{\Gamma_j(t-\xi)} S_j(y, t)ds(y) dt = \frac{1}{\sqrt{2\pi}}\int^{t_{\max}-\xi}_{t_{\min}-\xi} \int_{\Gamma_j(t)} S_j(y, t+\xi)ds(y) dt.
\en
Here $\Gamma_j(t)\subset D_j$ is defined as
$
\Gamma_j(t):=\{y\in D_j: \hat{x}\cdot y=t\}.
$
By the assumption (i) we deduce from \eqref{g} with $\xi=t_{\max} - \inf(\hat{x}\cdot D_j)-\epsilon$,  $\epsilon\in (0, \Lambda_j)$ that
\ben
g_j(\xi)&=&
 \int^{\inf(\hat{x}\cdot D_j)+\epsilon}_{\inf(\hat{x}\cdot D_j)+\epsilon-T} \int_{\Gamma_j(t)} S_j(y, t+\xi)ds(y) dt\\
 &=&\int^{\inf(\hat{x}\cdot D_j)+\epsilon}_{\inf(\hat{x}\cdot D_j)} \int_{\Gamma_j(t)} S_j(y, t+\xi)ds(y) dt\\
 &>& 0,
 \enn
because
\ben
&& t+\xi\in(t_{\max}-\epsilon, t_{\max})\quad\mbox{if}\quad t\in \Big(\inf(\hat{x}\cdot D_j), \inf(\hat{x}\cdot D_j)+\epsilon\Big),
\enn
and
\ben
\int_{\Gamma_j(t)} S_j(y, t+\xi)ds(y)=0\quad \mbox{if}\quad t<\inf(\hat{x}\cdot D_j)
\enn
due to the fact that $\Gamma_j(t)=\emptyset$.
Since $g_j$ coincides with the Fourier transform of $u_j^\infty$ by a factor, this proves the first part of the lemma. The analyticity of  $(\mathcal{F}u_j^\infty)(t)$ in $t\in I_j$ follows from \eqref{g} under the assumption (ii).
\end{proof}
Next we show that the multi-frequency far-field measurement data at a fixed observation direction can be uniquely split. Note that the splitting is obvious under the conditions in \eqref{ab}, because by Lemma \ref{lem1} the Fourier transform of $u^\infty(\hat{x}, k)$ has two disconnected components.
\begin{theorem}
Suppose that there are two time-dependent sources $S$ and $\tilde{S}$ with \mbox{supp}\;$\tilde{S}=\overline{\tilde{D}}\times [t_{\min}, t_{\max}]$ and $\tilde{D}=\tilde{D}_1\cup\tilde{D}_2$. Here the source function $\tilde{S}$ and its support $\overline{\tilde{D}}$ are also required to satisfy the assumptions $(i)$-$(iii)$. Let $\tilde{u}_j^\infty$ be defined by \eqref{FF} with $D_j, S_j$ replaced by $\tilde{D}_j$, $\tilde{S}_j:=\tilde{S}|_{\tilde{D}_j\times (t_{\min}, t_{\max})}$, respectively. Then the relation
\be\label{uinfinity}
\quad \quad \quad u^{\infty}(\hat{x}, k)=u_1^\infty(\hat{x}, k)+u_2^\infty(\hat{x}, k)=\tilde{u}_1^\infty(\hat{x}, k)+\tilde{u}_2^\infty(\hat{x}, k),\; k\in(k_{\min}, k_{\max})
\en implies that $u_j^\infty(\hat{x}, k)=\tilde{u}_j^\infty(\hat{x}, k)$ for $k\in[k_{\min}, k_{\max}]$ and $j=1,2$.
\end{theorem}
\begin{proof} By the analyticity of $u_j^\infty$, $\tilde{u}_j^\infty$ in $k\in \R$, the function $u^\infty$ can be analytically extended to the whole real axis.
Denote by $[T_{\min}, T_{\max}]$ the supporting interval of
the Fourier transform of $u^\infty$ with respect to $k$. In view of assumption (iii) and  Lemma \ref{lem1},  we may suppose without loss of generality that $T_{\max}=t_{\max} - \Lambda_{\min}$ with
\ben
\Lambda_{\min}=\inf(\hat{x}\cdot D_1)=\inf(\hat{x}\cdot \tilde{D}_1)<
\inf(\hat{x}\cdot D_2),\quad \Lambda_{\min}<
\inf(\hat{x}\cdot \tilde{D}_2).
\enn
If otherwise, there must hold
$T_{\min}=t_{\min} - \Lambda_{\max}$ with
\ben
\Lambda_{\max}=\sup(\hat{x}\cdot D_1)=\sup(\hat{x}\cdot \tilde{D}_1)>
\sup(\hat{x}\cdot D_2),\quad \Lambda_{\max}>
\sup(\hat{x}\cdot \tilde{D}_2)
\enn
and the proof can be carried out similarly.

Define $w_j=u_j^\infty-\tilde{u}_j^\infty$ for $j=1,2$. Using \eqref{uinfinity}, we get $w_1(\hat{x}, k)=-w_2(\hat{x}, k)$ for all $k\in \R$. Hence,  their Fourier transforms must also coincide, i.e., $[\mathcal{F}w_1](t)=- [\mathcal{F}w_2](t)$ for all $t\in \R$.
Taking $\delta<\min\{\inf(\hat{x}\cdot D_2)-\Lambda_{\min}, \; \inf(\hat{x}\cdot \tilde{D}_2)-\Lambda_{\min} \}$. Again using Lemma \ref{lem1}, we obtain
\ben
0=[\mathcal{F}w_2](t)=[\mathcal{F}w_1](t)\quad\mbox{for all}\quad t\in \left[T_{\max}-\delta, T_{\max}\right]
\enn
because the interval $[T_{\max}-\delta, T_{\max}]$ lies in the exterior of the supporting intervals of both $\mathcal{F} u_2^\infty$ and $\mathcal{F} \tilde{u}_2^\infty$.
Combining this with the analyticity of  $[\mathcal{F}w_1](t)$ in $t$, we get
\be\label{eq:2}
\quad\quad\quad[\mathcal{F}w_1](t)=0\;\;\mbox{for all}\;\;t\in\Big[T_{\min}, t_{\max}- \max\left\{\sup(\hat{x}\cdot D_1), \sup(\hat{x}\cdot \tilde{D}_1)\right\}\Big].
\en
If  $\sup(\hat{x}\cdot D_1)<\sup(\hat{x}\cdot \tilde{D}_1)$, it is seen from Lemma \ref{lem1} that
\be\label{eq:1}
[\mathcal{F}^{-1} u_1^\infty](t^*)=0,\qquad
[\mathcal{F}^{-1} \tilde{u}_1^\infty](t^*)>0
\en
where
\ben
 t^*=t_{\min} -\sup(\hat{x}\cdot D_1)-\epsilon \in I_j=\Big[ t_{\min}-\sup(\hat{x}\cdot \tilde{D}_1),  t_{\max}-\inf(\hat{x}\cdot \tilde{D}_1)\Big].
\enn

Obviously, the relations in \eqref{eq:1} contradicts the fact that  $[\mathcal{F}w_1](t^*)=0$ by \eqref{eq:2}. This proves  $\sup(\hat{x}\cdot D_1)\geq\sup(\hat{x}\cdot \tilde{D}_1)$. The relation $\sup(\hat{x}\cdot D_1)\leq\sup(\hat{x}\cdot \tilde{D}_1)$ can be proved analogously. Hence, $\sup(\hat{x}\cdot D_1)=\sup(\hat{x}\cdot \tilde{D}_1):=\Lambda_{\max}$ and
\ben
[\mathcal{F}w_1](t)=0\quad\mbox{for all}\quad t\in\Big[t_{\min}-\Lambda_{\max},t_{\max}-\Lambda_{\min} \Big].
\enn
Using again Lemma \ref{lem1} we find that $\mathcal{F}u_1^\infty$ and $\mathcal{F}\tilde{u}_1^\infty$ also vanish for $t\notin [t_{\min}-\Lambda_{\max},t_{\max}-\Lambda_{\min} ]$. Therefore, $w_1\equiv 0$ and $u_1^\infty\equiv \tilde{u}_1^\infty$, which implies $u_2^\infty\equiv \tilde{u}_2^\infty$.
\end{proof}

\begin{remark} Once $u_j^\infty(\pm\hat{x}, k)$ ($j=1,2$) can be computed from $u^\infty(\pm\hat{x}, k)$, one can apply the factorization scheme proposed in Sections 2-4 to get an image of the strip $K_{D_j}^{(\hat{x})}$ for $j=1,2$. The numerical implementation of the multi-frequency far-field splitting at a single observation direction is beyond the scope of this paper. We refer to 
\cite{GHS} for a numerical scheme splitting far-field patterns over all directions at a fixed frequency.
\end{remark}

\section{Numerical examples}\label{num}
In this section, we will  conduct a series of numerical experiments to validate our algorithm in $\R^2$ and $\R^3$. In practical scenarios, the time-domain data should be transformed into the multi-frequency data via the
inverse Fourier transform. For the sake of streamlining the numerical procedures in simulation, we will  only perform computational tests for the Helmholtz equation. Our primary objective is to extract information about the supports of sources. This goal is achieved through the utilization of multi-frequency far/near-field data recorded at either a single pair of opposite observation direction/point or a finite pair of observation directions/points.

\subsection{Numerical examples for the far-field case in $\R^2$ and $\R^3$}
 Assuming a wave-number-dependent source term $f(x, k)$, as defined in (\ref{fxt}), we can synthesize the far-field pattern using equation (\ref{u-infty}) by
\begin{equation}
\begin{aligned}
w^{\infty}(\hat x, k)=\frac{1}{\sqrt{2\pi}}\int_{t_{\min}}^{t_{\max}}  \int_{D} e^{ik(t-\hat x\cdot y)}S(x,t)\, dy dt,\; k\in[k_{\min}, k_{\max}].
\end{aligned}
\end{equation}
Below we will describe the process of inversion algorithm.
The frequency interval $[k_{\min},k_{\max}]$ can be discretized by defining $$k_n=(n-0.5)\Delta k, \quad \Delta k:=\frac{K}{N}, \quad n=1,2,\cdots,N.$$
We approximate the far-field operator in (\ref{def:F}) using the rectangle method:
\begin{equation}
    ({F}^{(\hat x)}\phi)(\tau_n) \approx \sum_{m=1}^{N} w(\hat x, k_c+\tau_n-s_m)\phi(s_m)\Delta k,
\end{equation}
 where $\tau_n:=n\Delta k$ and $s_m:=(m-0.5)\Delta k$, $n,m=1,2,\cdots,N$.
A discrete approximation of the far-field operator ${F}^{(\hat x)}$ is given  by the Toeplitz matrix
\be \label{matF}
 F^{(\hat x)}:= \Delta k \begin{pmatrix}
w^{\infty}({\hat x},k_c+k_1) & w^{\infty}({\hat x},k_c-k_1) & \cdots   & w^{\infty}({\hat x},k_c-k_{N-1})  \\
w^{\infty}({\hat x},k_c+k_2) & w^{\infty}({\hat x},k_c+k_1) & \cdots  &w^{\infty}({\hat x},k_c-k_{N-2})   \\
\vdots & \vdots  &\vdots  &\vdots \\
w^{\infty}({\hat x},k_c+k_{N-1}) & w^{\infty}({\hat x},k_c+k_{N-2}) &  \cdots  & w^{\infty}({\hat x},k_c-k_1)\\
w^{\infty}({\hat x},k_c+k_N) & w^{\infty}({\hat x},k_c+k_{N-1}) &  \cdots & w^{\infty}({\hat x},k_c+k_1)\\
\end{pmatrix},
\en
where $F^{(\hat x)}$ is a $N\times N$ complex matrix. Here, we adopt {$2N-1$ samples $w^{\infty}(\hat x, k_c+k_n), n=1,2,\cdots,N$ and $w^{\infty}(\hat x, k_c-k_n), n=1,2,\cdots,N-1$}, of the far-field pattern.
\noindent Denoting by  $\left\{ ( {\tilde \lambda^{( x)}_n}, \psi^{(x)}_n): n=1,2,\cdots,N \right\}$ an eigen-system of the matrix $F^{(\hat x)}$ (\ref{matF}), then one deduces that  an eigen-system of the matrix $(F^{(\hat x)})_\#:= |\real F^{(\hat x)})|+|\ima(F^{(\hat x)})|$ is $\left\{ ( \lambda^{(\hat x)}_n, \psi^{(\hat x)}_n): n=1,2,\cdots,N \right\}$, where $ \lambda^{(\hat x)}_n:=|\real (\tilde \lambda^{(\hat x)}_n)| +|\ima (\tilde \lambda^{( \hat x)}_n)|$.
For any test point $y\in \R^3$, test function $\phi_{y,\eta}^{(\hat x)}$  (\ref{test}) can be discretized as
\be \label{testn}
\phi_{y,\eta}^{(\hat x)}:= \left(\frac{1}{\eta-t_{\min}}\int_{t_{\min}}^{\eta} e^{i\tau_1(t-\hat x\cdot y)}dt, \,\cdots, \, \frac{1}{\eta-t_{\min}}\int_{t_{\min}}^{\eta} e^{i\tau_N(t-\hat x\cdot y)}dt\right).
\en
We truncate the auxiliary indicator function $I_{\eta,\epsilon}^{(\hat{x})}(y)$   (\ref{indicator4}) by
\be\label{indicator5}
I_{\eta}^{(\hat{x})}(y)\sim \sum_{n=1}^N\frac{| \phi^{(\hat{x})}_{y,\eta} \cdot \psi_n^{(\hat{x})}|^2}{ |\lambda_n^{(\hat{x})}|}, \qquad y\in \R^3,
\en
where $\cdot$ denotes the inner product in $\R^3$ and $N$ is consistent with the dimension of the Toeplitz matrix \eqref{matF}.
Furthermore, an approximation of the indicator function $W(y)$ in (\ref{W2}) is

\begin{equation} \label{indicator7}
W(y)\rot{= \left[\sum_{j=1}^M \frac{1}{W^{(\hat{x}_j)}(y) }\right]^{-1}}\sim \left[\sum_{j=1}^M\sum_{n=1}^N\frac{| \phi^{(\hat{x}_j)}_{y,\eta}, \psi_n^{(\hat{x}_j)} |^2}{ |\lambda_n^{(\hat{x}_j)}|} + \frac{| \phi^{(-\hat{x}_j)}_{y,\eta}, \psi_n^{(-\hat{x}_j)} |^2}{ |\lambda_n^{(-\hat{x}_j)}|} \right]^{-1}, \; y\in \R^3.
\end{equation}
It can also be reduced to the approximation of $W^{(\hat x)}(y)$ in (\ref{w-eps}) as $M=1$.
The visualizations of the strips $K_{D,\eta}^{(\hat x)}$, $ K_D^{(\hat x)}$ (see Figure \ref{strip}) and $\Theta$-convex hull of source support  $D$ are attainable through the plots of $1/I^{(\hat x)}_{\eta}$ defined as (\ref{indicator5}) and  $W^{(\hat x)}(y)=W(y)$ defined as (\ref{indicator7}) for $M=1$.
We suppose  $k_{\min}=0$ for the sake of simplicity, then
the bandwidth of frequency can be extended from $(0,k_{\max})$ to $(-k_{\max}, k_{\max})$ by $w^{\infty}(\hat x, -k)=\overline{w(\hat x, k)}$. Thus, one deduces from these new measurement data with $k_{\min}=-k_{\max}$ that $k_c=0$ and $K=k_{\max}$.  The frequency band is represented by the interval $(0, 8\pi/3)$ and discretized by  $N=16$ and $\Delta k=\pi/6$. The source function is chosen as  $S(x,t)=3|x|^2(t+1)$ in $\R^3$. The radiating time is supposed to be $[t_{\min},t_{\max}]=[0,2]$ in $\R^2$ and $[t_{\min},t_{\max}]=[0,0.5]$ in $\R^3$.  We always take $\eta=0.1$ unless other specified. In the figures below, the shape of source support $D$ will be highlighted by the red/pink solid line and the indicator function values are all  normalized to their respective maximum values. The shapes of the source support related in our two dimensional numerical examples are listed below:
\begin{itemize}
\item Peanut: \quad \quad \quad $x(t)=(x_1,x_2)+ar\sqrt{b\cos^2 t+1}(\cos t, \sin t), \quad t\in[0,2\pi]$,
\item Round square: $x(t)=(x_1,x_2)+r( \cos^3 t+\cos t,  \sin^3 t+\sin t), \quad t\in[0,2\pi]$,
\item Kite:\quad \quad \quad \quad \;\,$x(t)=(x_1,x_2)+(r\cos t+ar(\cos 2t-1), br\sin t), \quad t\in[0,2\pi]$,
\item Ellipse: \quad \quad \quad $x(t)=(x_1,x_2)+(a\cos t, b\sin t), \quad t\in[0,2\pi]$.
\end{itemize}

\subsubsection{Recovery of $K_D^{(\hat{x})}$ and $\Theta_D$ with far-field measurements in $\R^2$.}
\begin{figure}[H]
\centering
\subfigure[$\eta=0.3$]{
\includegraphics[scale=0.22]{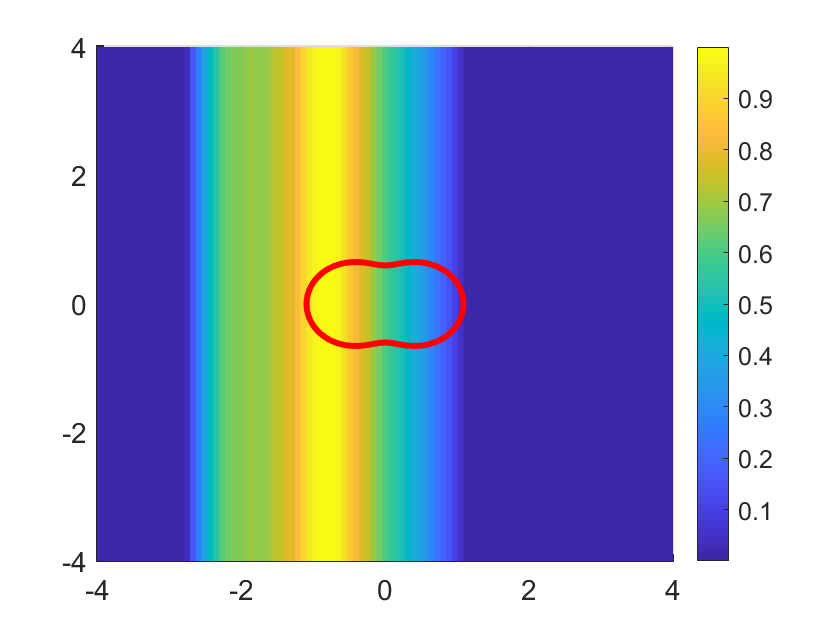}
}
\subfigure[$\eta=1$ ]{
\includegraphics[scale=0.22]{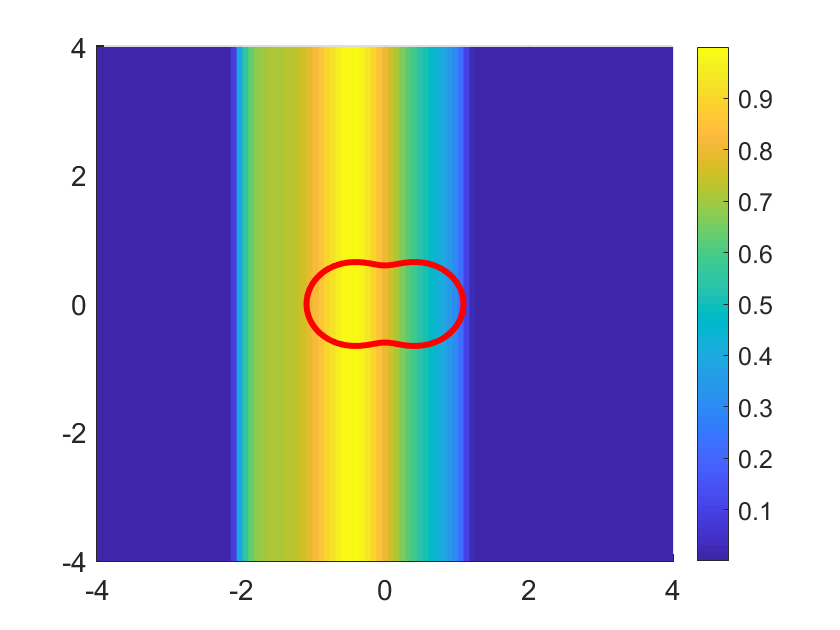}

}
\subfigure[$\eta=1.8$]{
\includegraphics[scale=0.22]{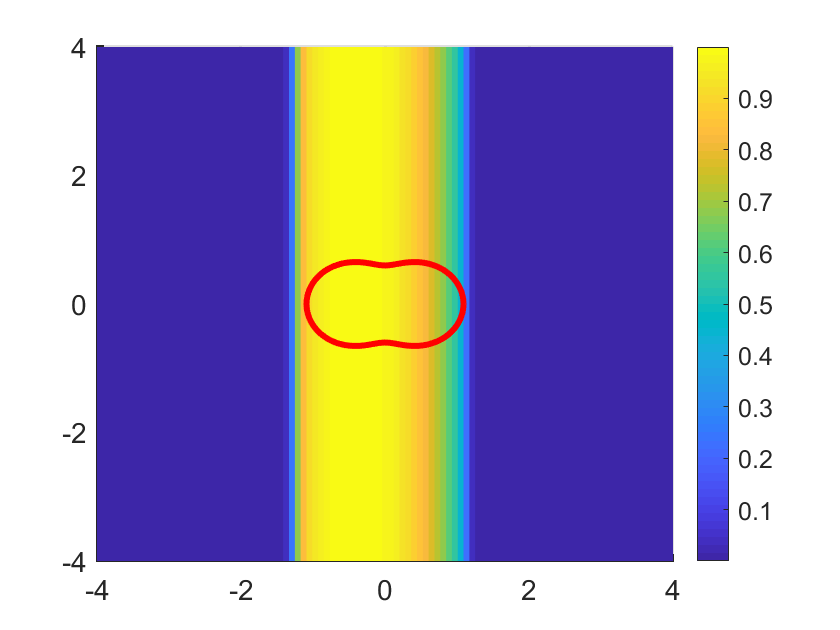}

}
\caption{Reconstructions using multi-frequency far-field data from a single observation direction for a peanut-shaped source support with the  auxiliary indicator function $1/I^{(\hat x)}_{\eta}$ defined in (\ref{indicator5}). 
} \label{fig:1}
\end{figure}

Firstly, we focus on reconstructing strips  $K_{D,\eta}^{(\hat x)}$ in (\ref{K}) and $K_D^{(\hat x)}$ in (\ref{tildeK}) perpendicular to the observation direction and  containing a peanut-shaped source support using  a single observation direction. The search domain is selected as $[-4,4]\times[-4,4]$. In Figure \ref{fig:1}, we present reconstructions $K_{D,\eta}^{(\hat x)}$ using multi-frequency far-field data from a single observation direction $\hat x=(1,0)$ for a peanut-shaped support with the  auxiliary indicator function $1/I^{(\hat x)}_{\eta}(y)$ defined in (\ref{indicator5}). We test various values for $\eta$, namely $\eta=0.3, 1, 1.8$. It is evident that both the right boundary of the strips $\{y\in \R^2: \hat x\cdot y=\sup (\hat x\cdot D)\}$  and the left boundary  $\{y\in \R^2: \hat x\cdot y=\inf (\hat x\cdot D)-t_{\max}+\eta\}$ are perfectly reconstructed. The width of strips  $K_{D,\eta}^{(\hat x)}$ in Figure \ref{fig:1} is $t_{\max}-\eta$ units wider than the width of the smallest strip $K_D^{(\hat x)}$ containing the peanut-shaped source support.
As $\eta$ ($\eta\leq t_{\max}$) increases, strips $K_{D,\eta}^{(\hat x)}$ get progressively closer to strip $K_D^{(\hat x)}$.  In the limiting case where $\eta= t_{\max}$, it follows that  $ K_{D,\eta}^{(\hat x)}=K_D^{(\hat x)}$.

\begin{figure}[H]
\centering
\subfigure[$1/I_{\eta}^{(\hat{x})}(y)$]{
\includegraphics[scale=0.22]{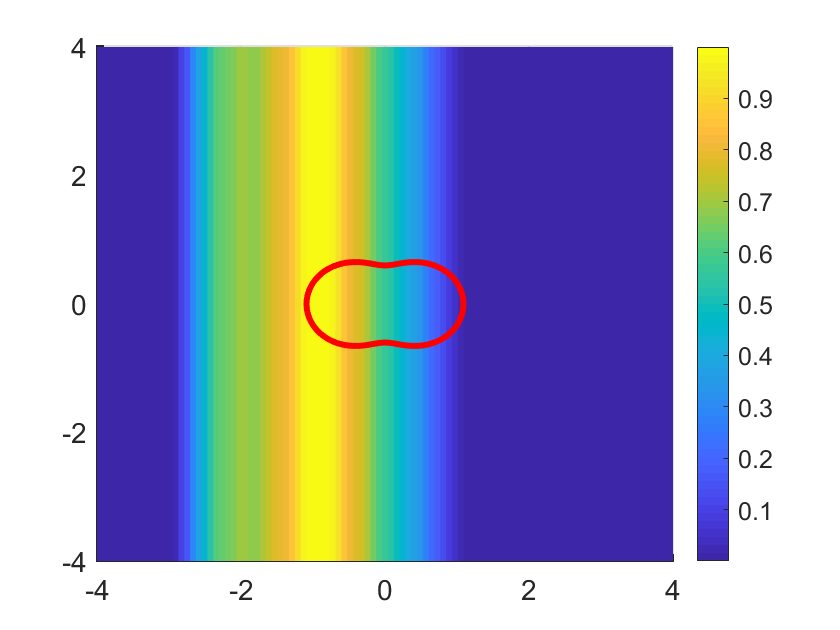}
}
\subfigure[$1/I_{\eta}^{(-\hat{x})}(y)$ ]{
\includegraphics[scale=0.22]{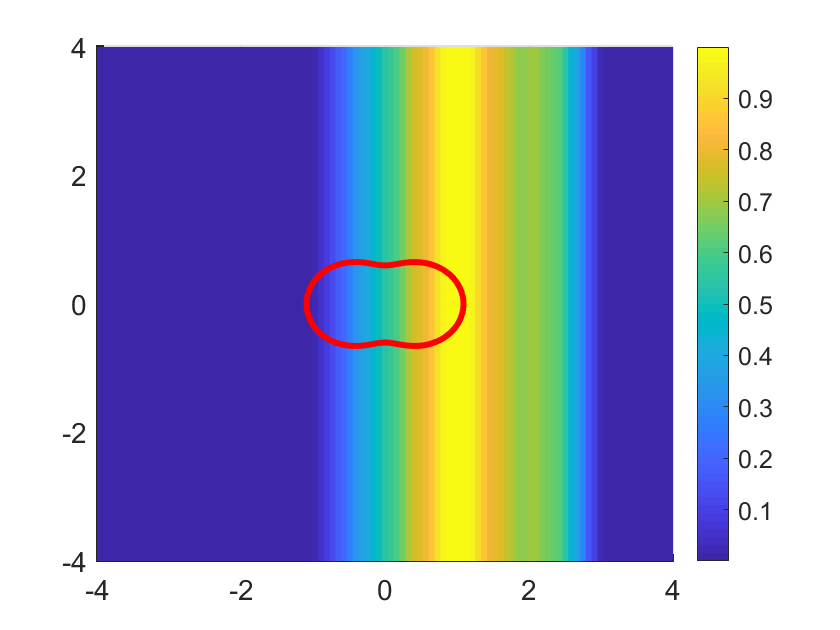}

}
\subfigure[$W^{(\hat{x})}(y)$]{
\includegraphics[scale=0.22]{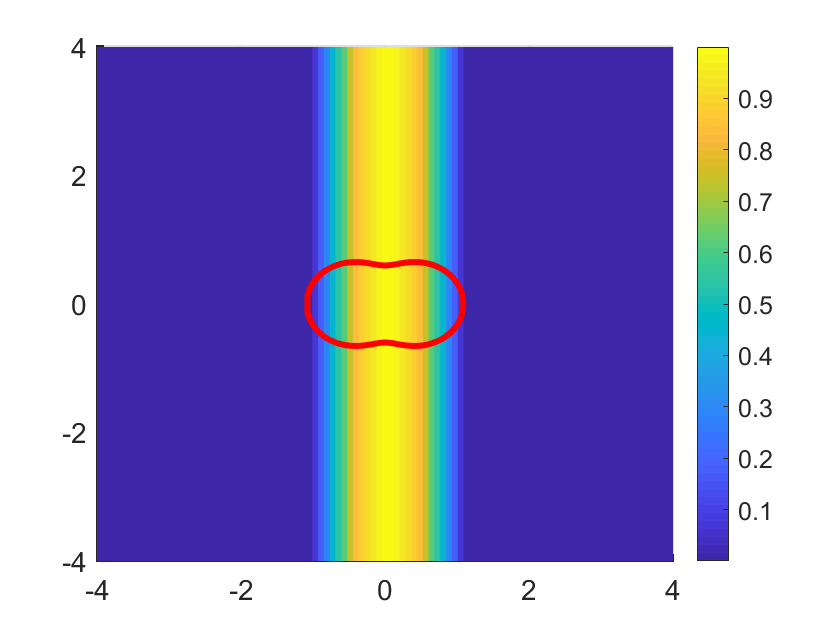}

}
\subfigure[$\delta=3\times10^{-3}$]{
\includegraphics[scale=0.22]{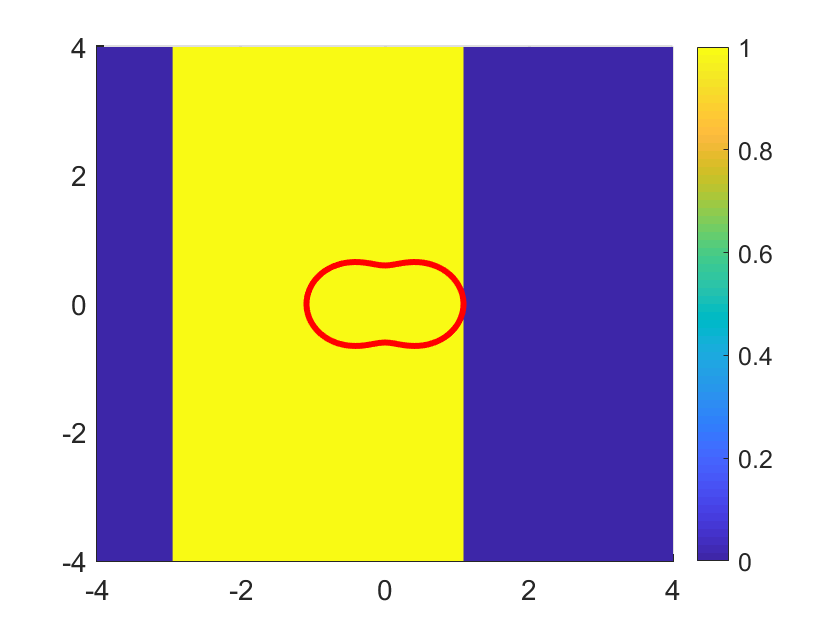}

}
\subfigure[$\delta=3\times10^{-3}$]{
\includegraphics[scale=0.22]{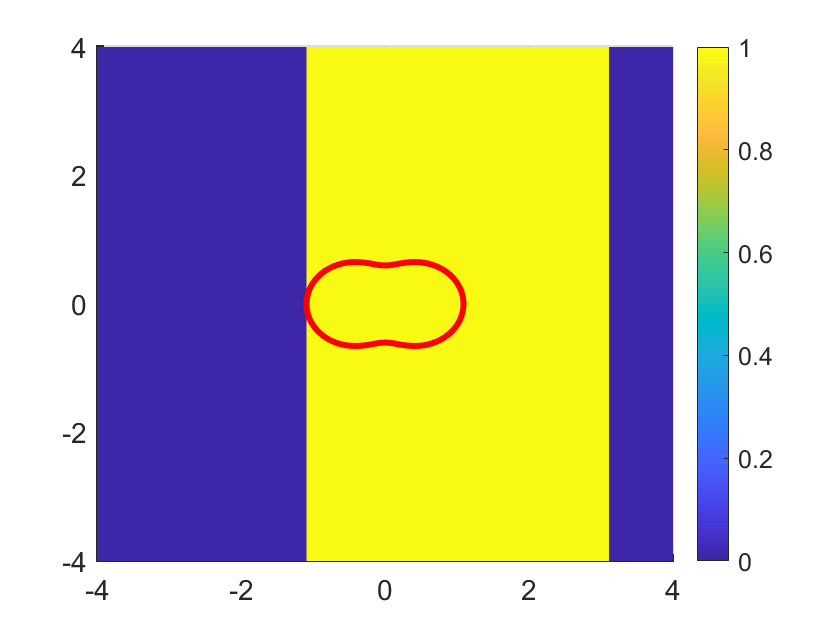}

}
\subfigure[$\delta=3\times10^{-3}$]{
\includegraphics[scale=0.22]{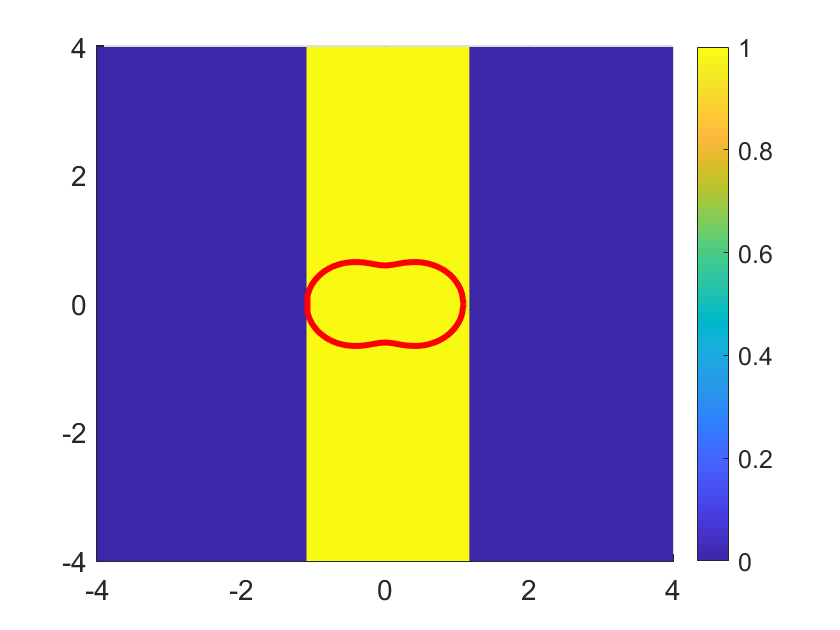}

}
\caption{Reconstructions using multi-frequency far-field data from a single observation direction for a peanut-shaped support.  We employ an auxiliary indicator function  as defined in (\ref{indicator5}) and indicator function as defined in (\ref{indicator7}) with $M=1$.
} \label{fig:2}
\end{figure}

\noindent Figure \ref{fig:2} illustrates the reconstructions  using multi-frequency far-field data from a single observation direction $\hat x=(1,0)$ for a peanut-shaped support. We employ the auxiliary indicator functions $1/I_{\eta}^{(\hat{x})}(y)$ in Figure \ref{fig:2}(a) and  $1/I_{\eta}^{(-\hat{x})}(y)$ in Figure \ref{fig:2} (b). Additionally,  an indicator function $W^{(\hat{x})}(y)=W(y)$  (\ref{indicator7}) defined with $M=1$, is utilized in Figure \ref{fig:2}(c). It can be observed that the strip $K_D^{(\hat x)}$ in Figure \ref{fig:2}(c) is the intersection of the strip $K_{D,\eta}^{(\hat x)}$ in Figure \ref{fig:2}(a) and the strip $ K_{D,\eta}^{(-\hat x)}$ in Figure \ref{fig:2}(b), as proven by our theory in Theorem \ref{Th:factorization}.  To enhance the intuitiveness of inversion results, we introduce a threshold $\delta > 0$ in Figures \ref{fig:2}(d)–(f), respectively. It is evident that $1/I_{\eta}^{(\pm\hat{x})}(y)$ can effectively capture the strip $\hat x \cdot D$ with a shift $t_{\max}-\eta=1.9$ of the lower boundary along the observation direction $\pm\hat{x}$. The strip $K_D^{(\hat x)}$ has been accurately recovered by the normalized indicator function.

\begin{figure}[H]
\centering
\subfigure[$M=2$]{
\includegraphics[scale=0.22]{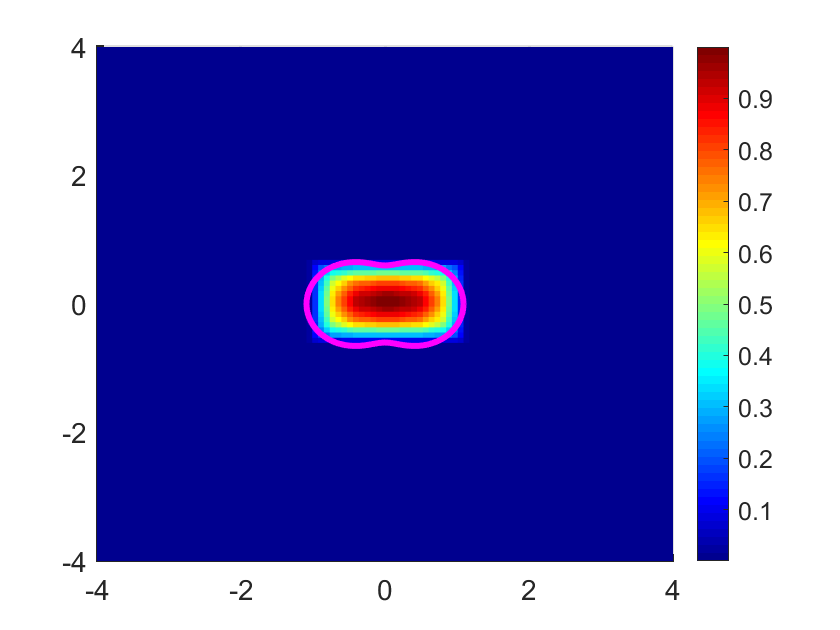}
}
\subfigure[$M=4$ ]{
\includegraphics[scale=0.22]{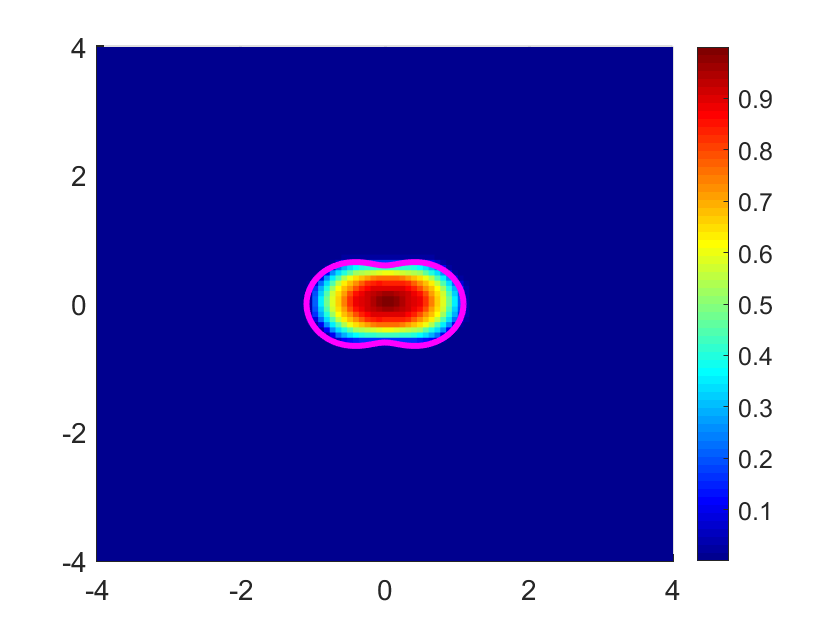}

}
\subfigure[$M=8$]{
\includegraphics[scale=0.22]{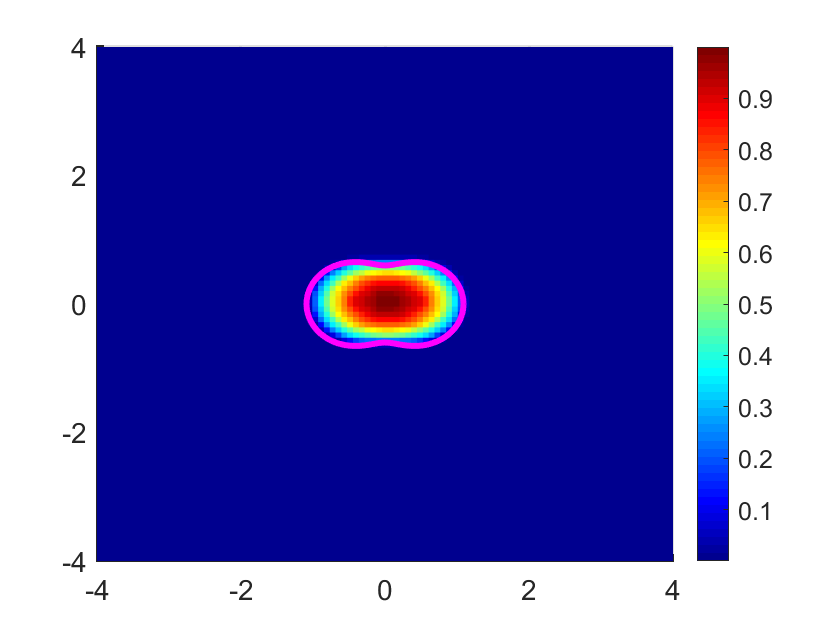}

}
\caption{Reconstructions using multi-frequency far-field data from $M$ pairs of opposite observation directions for a peanut-shaped source support.} \label{fig:3}
\end{figure}

\begin{figure}[H]
\centering
\subfigure[$M=2$]{
\includegraphics[scale=0.22]{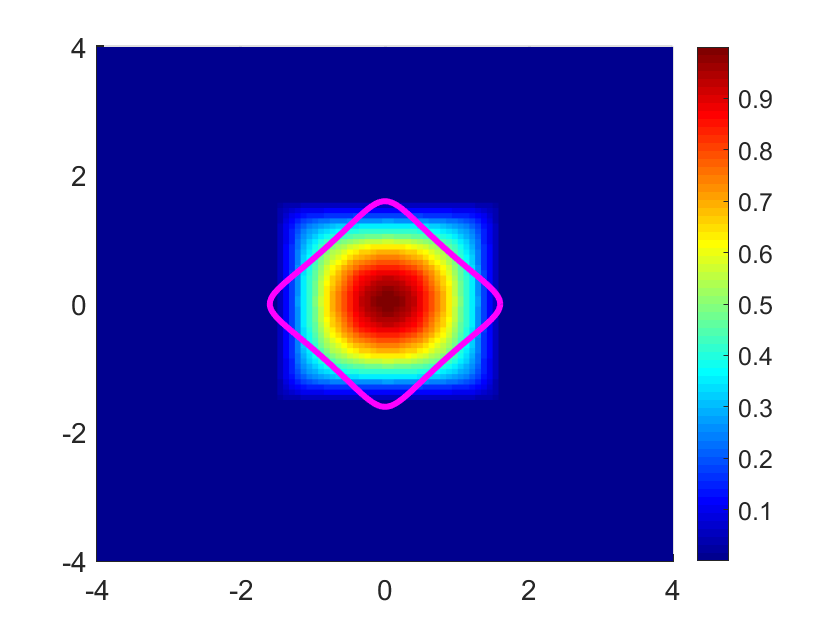}
}
\subfigure[$M=4$ ]{
\includegraphics[scale=0.22]{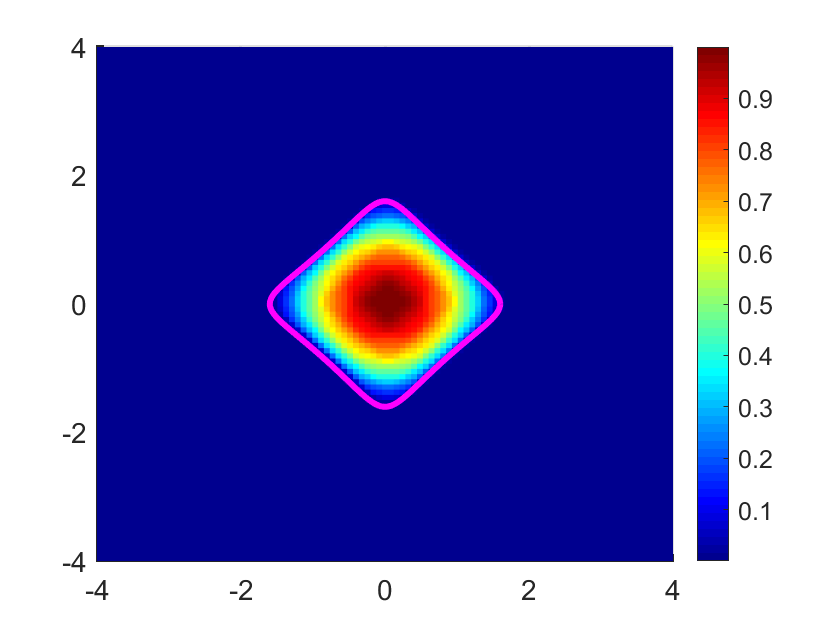}

}
\subfigure[$M=8$]{
\includegraphics[scale=0.22]{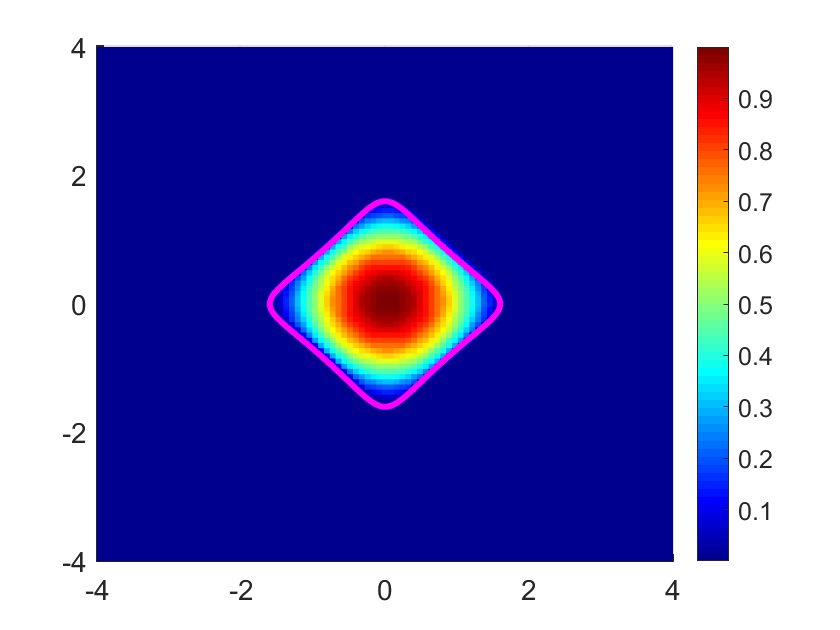}

}
\caption{Reconstructions using multi-frequency far-field data from $M$ pairs of opposite  observation directions for a rounded-square-shaped source support.
} \label{fig:4}
\end{figure}

%
%
Next, we conduct numerical tests  using $M$ pairs of opposite observation directions $\hat x_m= (\cos \theta_m, \sin \theta_m)$ with $\theta_m=\frac{(m-1)\pi}{M}, m=1, \cdots, M.$  The reconstructions of a peanut-shaped source support are presented in Figure \ref{fig:3} and those for a rounded-square-shaped source support in Figure \ref{fig:4}, respectively.
It is evident from Figures \ref{fig:3}(a) and \ref{fig:4}(a) that the locations of the supports are precisely captured by the intersection of two strips $K^{(\pm \hat{x})}_{D,\eta}$ from two different pairs of opposite observation directions. The two observations are perpendicular to each other with $M=2$, hence, yielding the smallest rectangle/square containing the source support. In fact, the position of the source can be reconstructed from multi-frequency data taken from any two observation directions. In Figures \ref{fig:3}(b) and \ref{fig:4}(b), we  utilize multi-frequency far-field date from $M=4$ pairs of opposite observation directions.  It is evident that both the location and shape are accurately captured. Similarly, by setting $M=8$, the shape of the supports can be more precisely inverted in Figures \ref{fig:3}(c) and \ref{fig:4}(c).

\begin{figure}[H]
\centering
\subfigure[$\eta=0.1$]{
\includegraphics[scale=0.22]{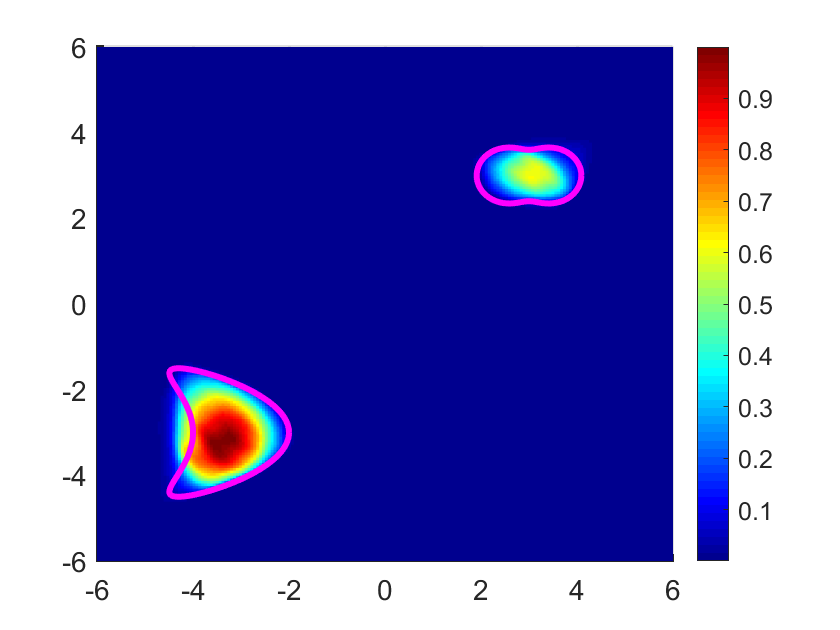}
}
\subfigure[$\eta=0.5$ ]{
\includegraphics[scale=0.22]{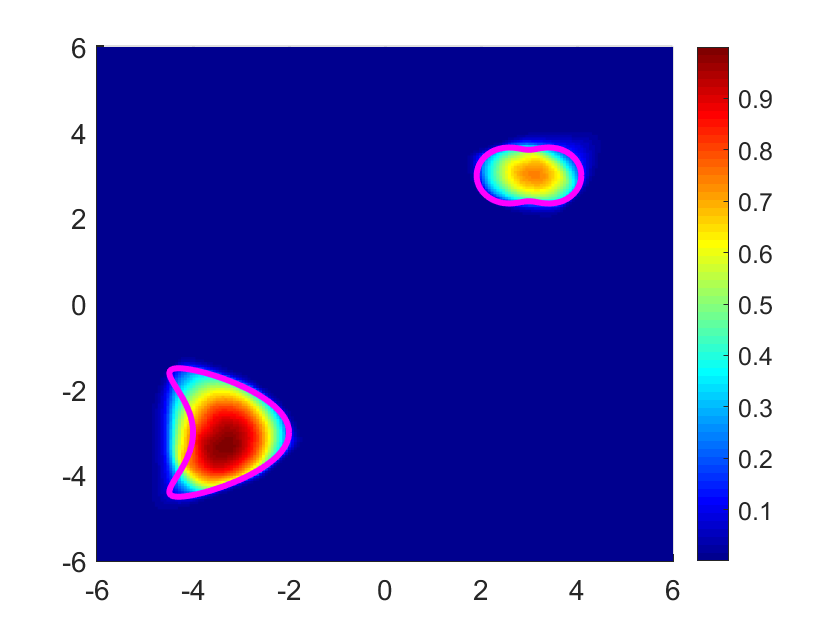}

}
\subfigure[$\eta=0.9$]{
\includegraphics[scale=0.22]{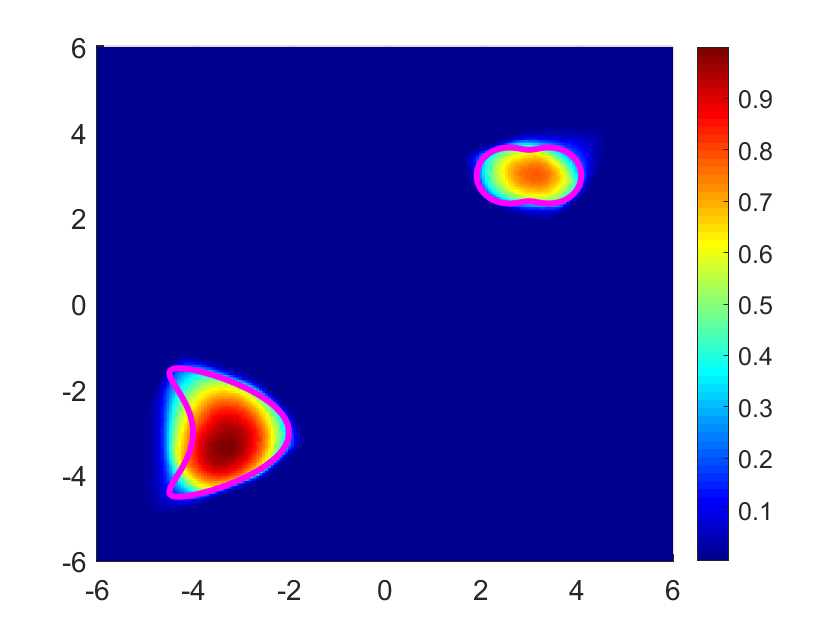}

}
\caption{Reconstructions using multi-frequency far-field data from $8$ pairs of opposite observation directions  for a peanut-shaped and a kite-shaped support with various $\eta$. We set the radiating period $[t_{\min}, t_{\max}]=[0,1]$.
} \label{fig:5}
\end{figure}

\begin{figure}[H]
\centering
\subfigure[$T=2$]{
\includegraphics[scale=0.22]{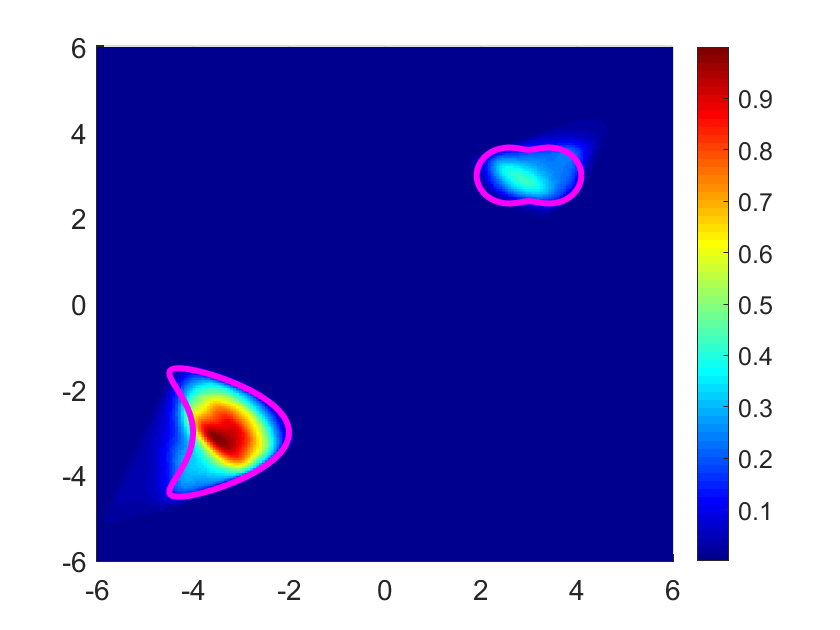}
}
\subfigure[$T=3$ ]{
\includegraphics[scale=0.22]{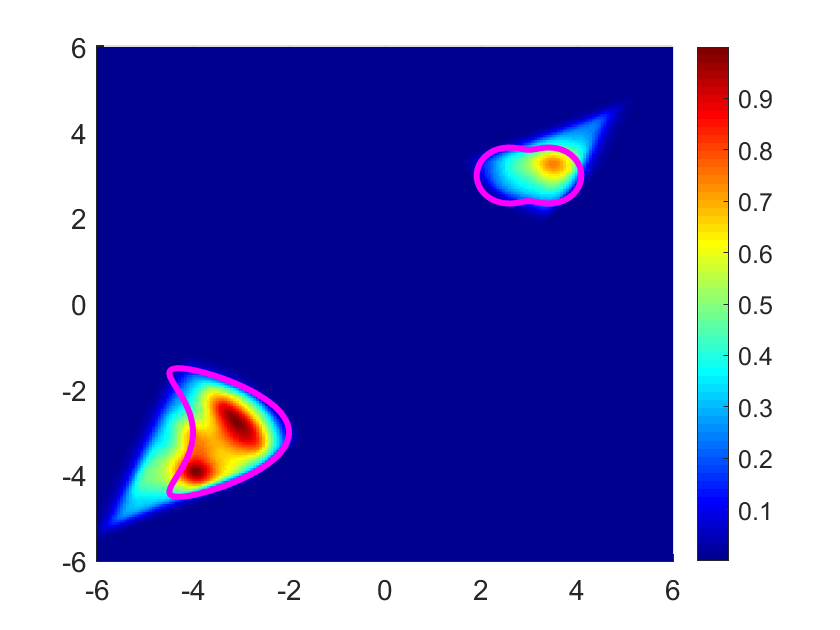}

}
%
\subfigure[$T=5$]{
\includegraphics[scale=0.22]{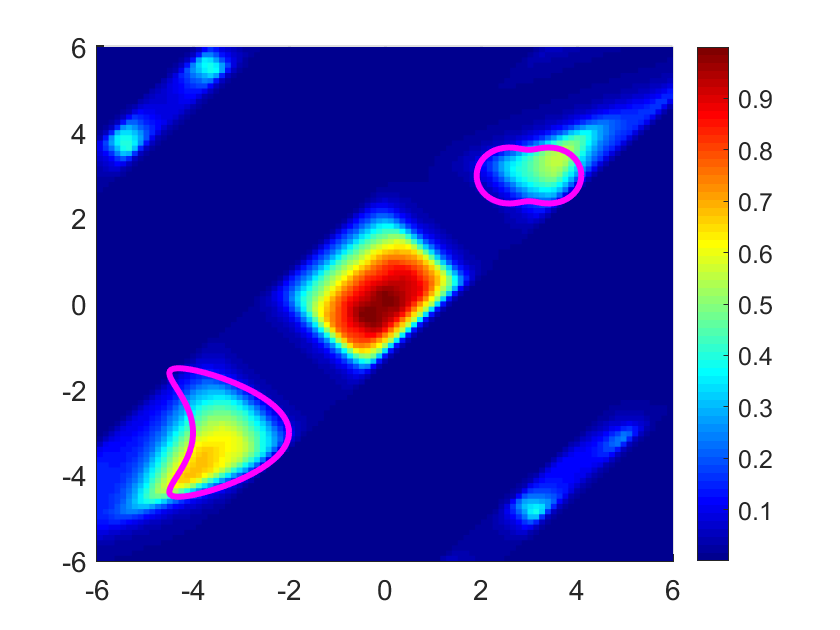}

}
\caption{Reconstructions using multi-frequency far-field data from $8$ pairs of opposite observation directions  for a peanut-shaped and a kite-shaped support with various radiating periods $[t_{\min}, t_{\max}]=[0,T]$.
} \label{fig:6}
\end{figure}

\begin{figure}[H]
\centering
\subfigure[$\eta=0.1$]{
\includegraphics[scale=0.22]{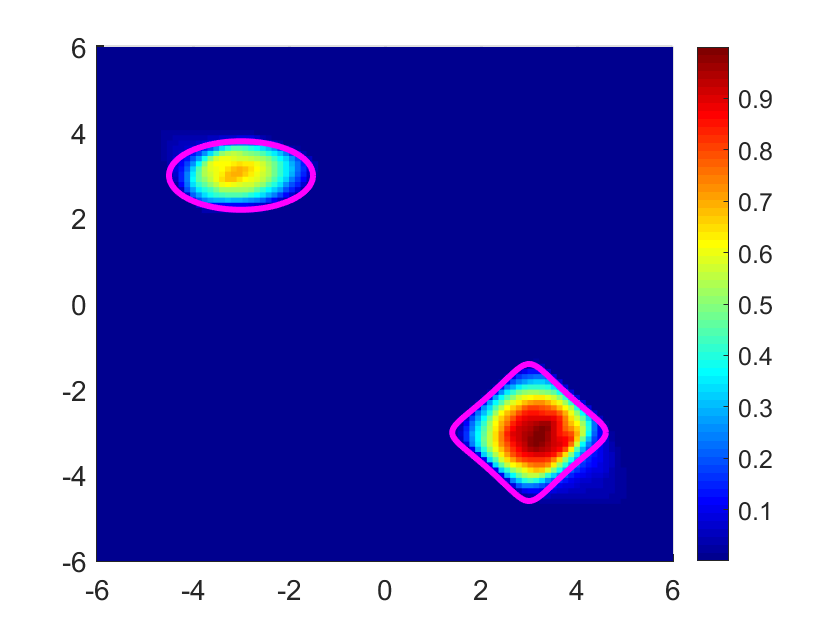}
}
\subfigure[$\eta=0.5$ ]{
\includegraphics[scale=0.22]{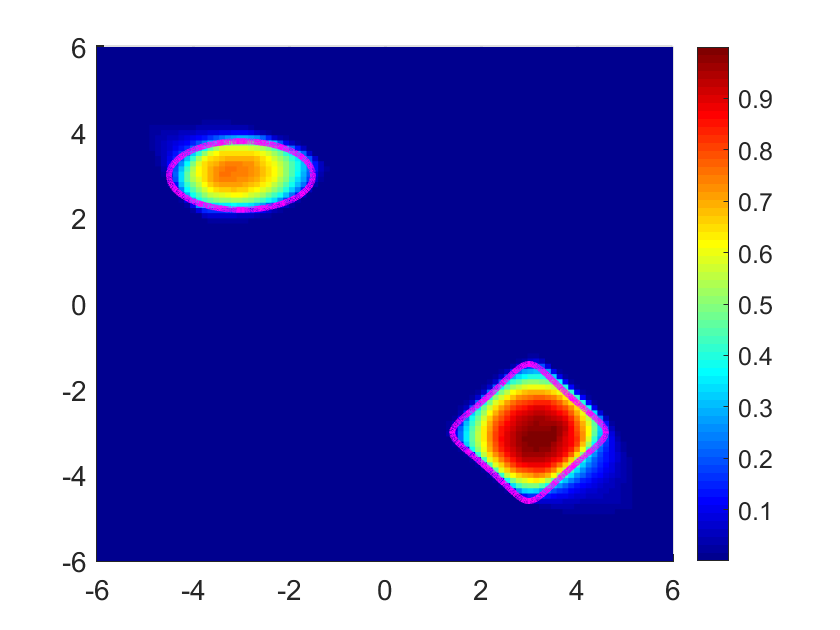}

}
\subfigure[$\eta=0.9$]{
\includegraphics[scale=0.22]{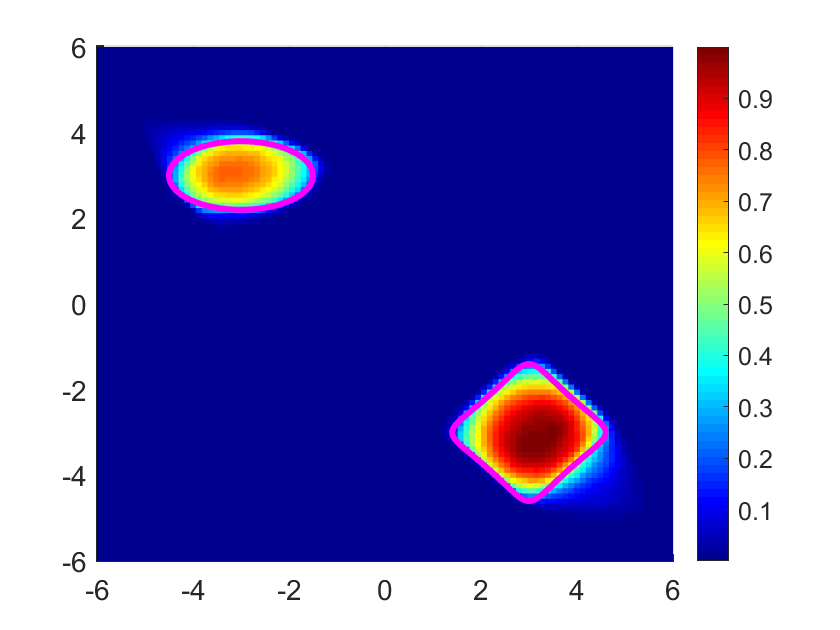}

}
\caption{Reconstructions using multi-frequency far-field data from $8$ pairs of opposite observation directions  for an elliptic-shaped and a rounded-square-shaped support with various $\eta$. We set the radiating period $[t_{\min}, t_{\max}]=[0,1]$.
} \label{fig:7}
\end{figure}

\begin{figure}[H]
\centering
\subfigure[$T=2$]{
\includegraphics[scale=0.22]{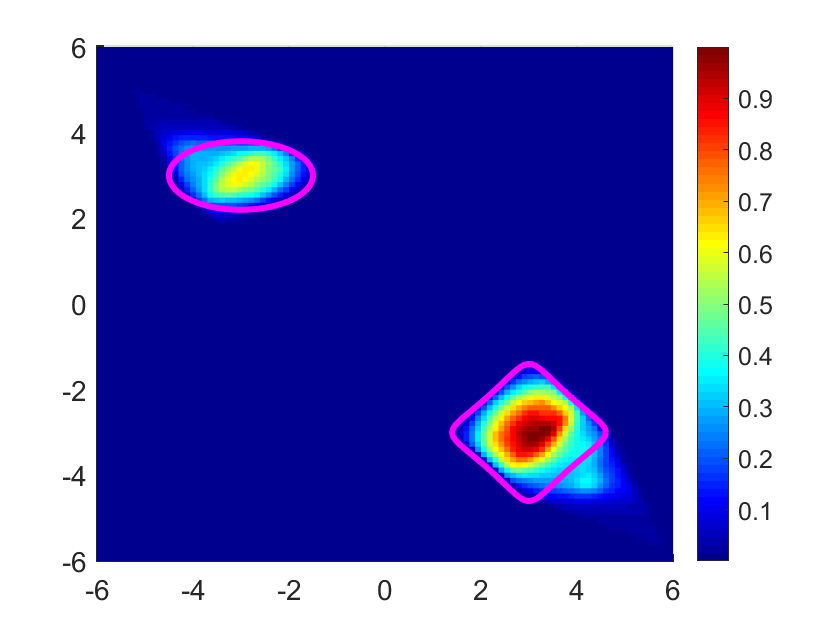}
}
\subfigure[$T=3$ ]{
\includegraphics[scale=0.22]{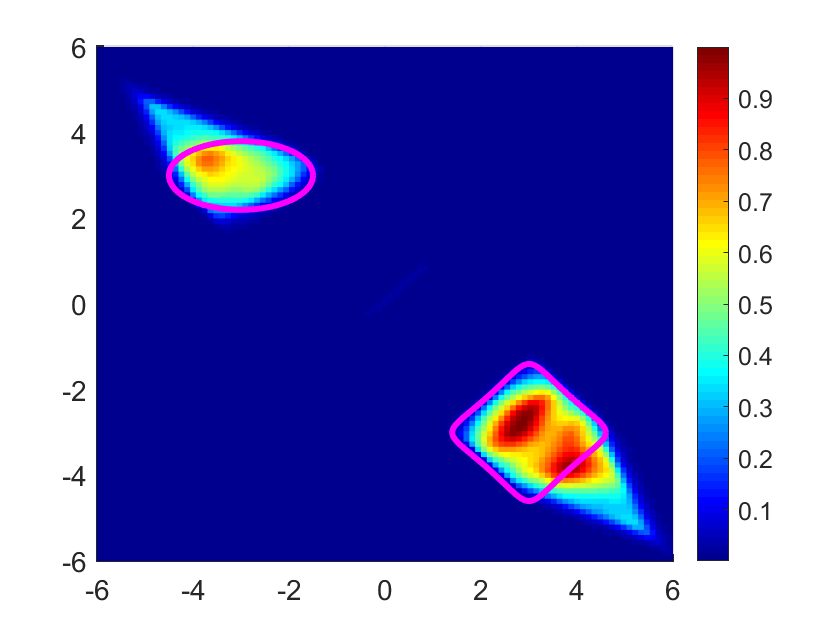}

}
%
\subfigure[$T=5$]{
\includegraphics[scale=0.22]{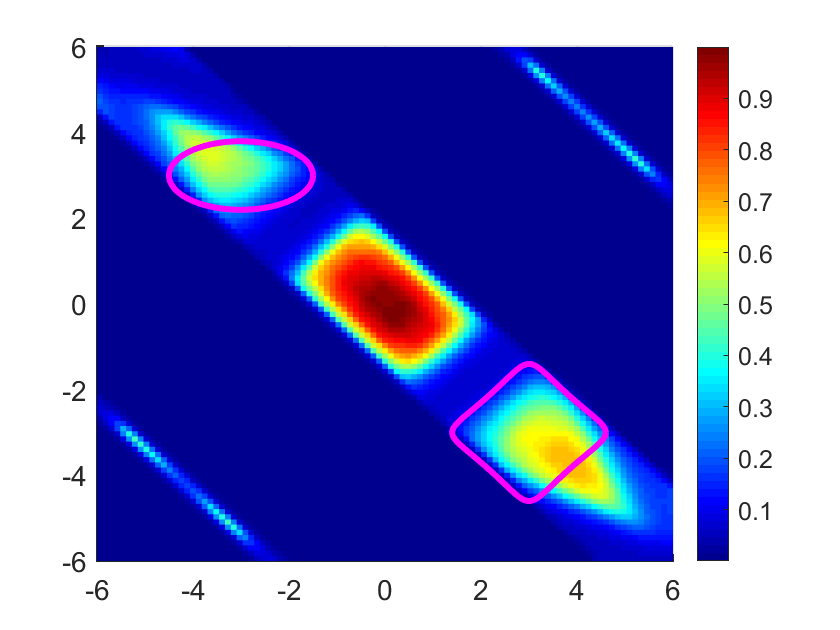}

}
\caption{Reconstructions using multi-frequency far-field data from $8$ pairs of opposite  observation direction  for an elliptic-shaped and a rounded-square-shaped support with various radiating periods $[t_{\min}, t_{\max}]=[0,T]$.
} \label{fig:8}
\end{figure}

Subsequently, we consider a source support with two disconnected components,  one kite-shaped and the other peanut-shaped depicted in Figures \ref{fig:5}-\ref{fig:6}, one elliptic-shaped and the other rounded-square-shaped in Figures \ref{fig:7}-\ref{fig:8}.  We examine various values of $\eta$ in the test function  and different radiating periods $[t_{\min}, t_{\max}]$ of the source. The search domain is chosen as $[-6, 6]\times[-6,6]$.  In Figures \ref{fig:5}-\ref{fig:6},  the peanut and kite are centered at the points $(3,3)$ and (-3,-3), respectively.  In Figures \ref{fig:7}-\ref{fig:8}, the ellipse and round square are centered at the points $(-3,3)$ and (3,-3), respectively.
We present visualizations of indicator function (\ref{indicator7}) with $M=8$ for a peanut-shaped and a kite-shaped support with $\eta=0.1, 0.5, 0.9$ in Figure \ref{fig:5}. The radiating period of the source is set to  $[t_{\min}, t_{\max}]=[0,1]$. It is evident that both the shapes and locations are accurately reconstructed. The inverted results become more  accurate with the increasing $\eta$ ($\eta\leq t_{\max}=1$).
If we enlarge the radiating period, can we still obtain results as accurate as those in Figure \ref{fig:5}? In Figure \ref{fig:6}, we set $[t_{\min}, t_{\max}]=[0,T]$ with $T=2,3,5$.  It turns out that increasing radiating period leads to distorted numerical reconstructions. This is due to the reason that the conditions (\ref{ab}) are no longer satisfied and the wave signals from two components are interfered with each other for large $T$; see
Figure \ref{fig:6}(c).
The accuracy of Figure \ref{fig:6}(b) can be improved by increasing the number of observation directions.
In Figures \ref{fig:7} and \ref{fig:8}, we perform extensive tests by reconstructing an elliptic-shaped component and a rounded-shaped square component, which further confirm the correctness of our algorithm.

\subsubsection{Determination of $t_{\max}$ from far-field measurements in $\R^2$.}\label{sub-tmax}

\begin{figure}[H]
\centering
\subfigure[$\hat x=(1,0), \theta=\pi$]{
\includegraphics[scale=0.22]{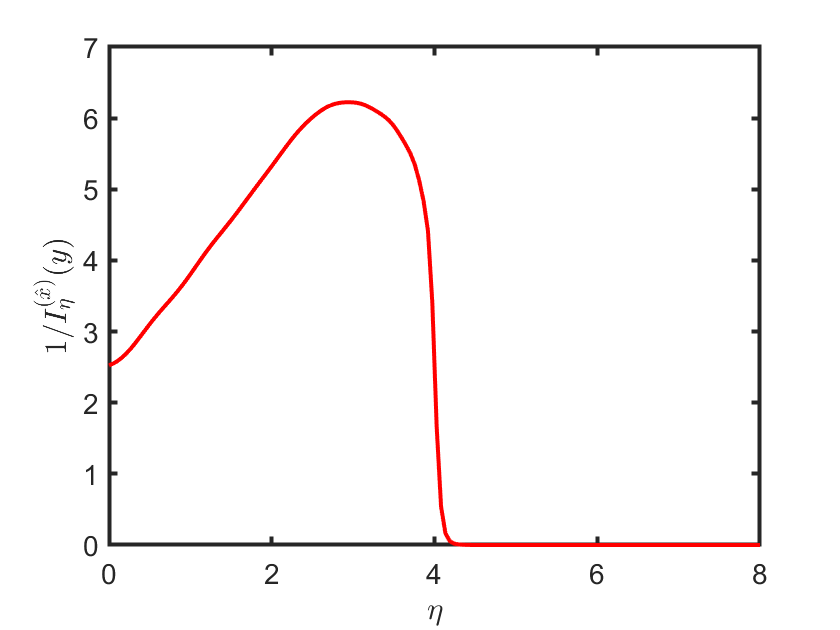}
}
\subfigure[$\hat x=(\sqrt{2}/2, \sqrt{2}/2),\theta=5\pi/4$ ]{
\includegraphics[scale=0.22]{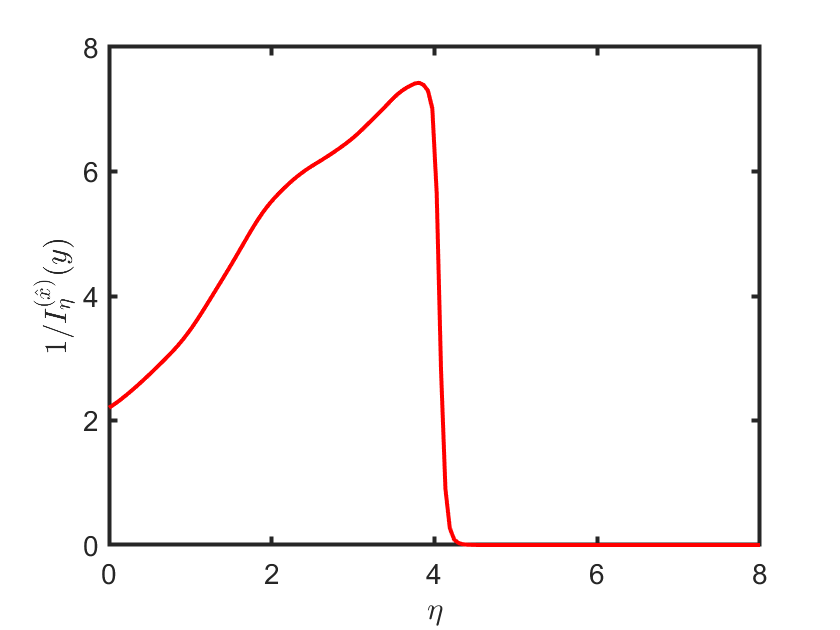}

}
\subfigure[$\hat x=(0,1), \theta=3\pi/2$]{
\includegraphics[scale=0.22]{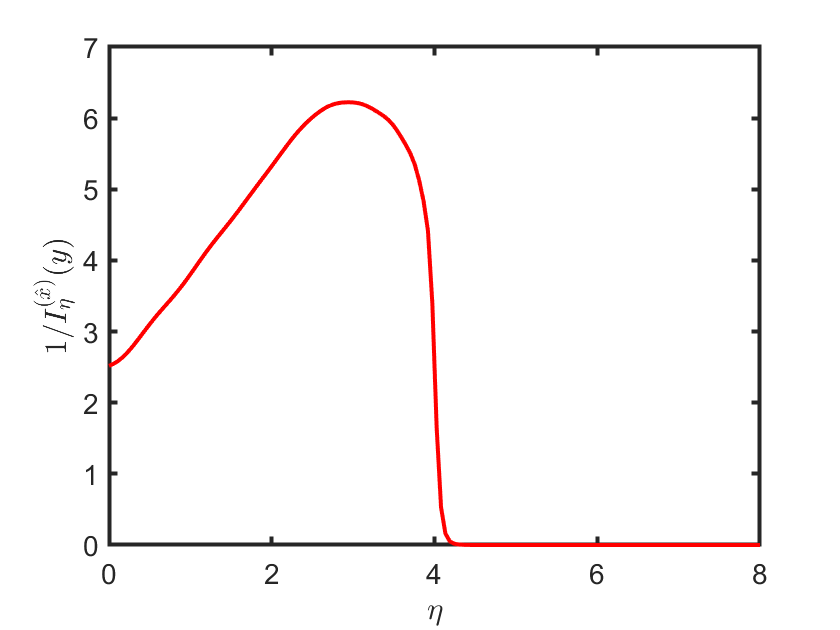}

}
\caption{Determination of $t_{\max}$ with different observation directions $\hat x$ and a fixed test point $y(\theta)=0.8( \cos^3 \theta+\cos \theta,  \sin^3 \theta+\sin \theta) \in K^{(\hat x)}_D$ such that $\hat x \cdot y-\inf (\hat x \cdot D)=\epsilon_0=0.01$ by  plotting the one-dimensional function  $\eta\rightarrow 1/ I_{\eta}^{(\hat x)}(y)$ where $I_{\eta}^{(\hat x)}(y)$ is defined in (\ref{indicator5}). The radiating period of the source is $[t_{\min}, t_{\max}]=[0,4]$ for a rounded-square-shaped source support.
} \label{fig:tmax-1}
\end{figure}

Assume that the initial moment $t_{\min}$ for source radiating is known, but the terminal moment $t_{\max}$ is unknown. Suppose that $\partial D$ is the square parameterized by $x(\theta)=0.8( \cos^3 \theta+\cos \theta,  \sin^3 \theta+\sin \theta)$ and the radiating period is $[t_{\min}, t_{\max}]=[0,4]$. Our objective is to determine the terminal time $t_{\max}$ using the far-field measurement $w^{\infty}(\hat x, k)$ from a specific observation direction $\hat x$. The frequency band is represented by the interval $(0, 8\pi/3)$ and discretized by  $N=32$ and $\Delta k=\pi/12$. 
In our tests, we use different directions and fix the test point $y(\theta)=0.8( \cos^3 \theta+\cos \theta,  \sin^3 \theta+\sin \theta)\in K^{(\hat x)}_D$ to plot the one-dimensional function  $\eta\rightarrow 1/ I_{\eta}^{(\hat x)}(y)$ in Figure \ref{fig:tmax-1}.
Notably, the one-dimensional function $\eta\rightarrow 1/ I_{\eta}^{(\hat x)}(y)$ exhibits a rapid decay near  $\eta=4$ and the value of $1/ I_{\eta}^{(\hat x)}(y)$ tends to $0$ for all $\eta>t_{\max}+\epsilon_0$. This aligns with our theoretical prediction in Theorem \ref{Th:max}, stating that the function $\eta\rightarrow 1/I_{\eta,\epsilon}^{(\hat x)}(y)$ must have a rapid decay  near $\eta=t_{\max}$, indicates  the value of $t_{\max}$. Consequently, the terminal time can be clearly seen as $t_{\max}=\eta=4$.

\begin{figure}[H]
\centering
\subfigure[$\epsilon_0=0, y=y_0+\epsilon_0\hat x$]{
\includegraphics[scale=0.22]{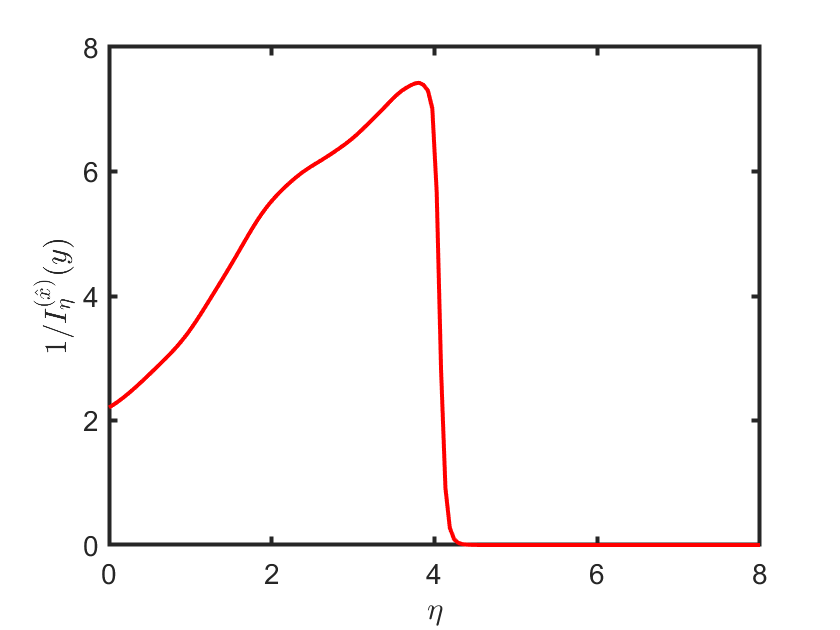}
}
%
\subfigure[$\epsilon_0=0.1, y=y_0+\epsilon_0\hat x$]{
\includegraphics[scale=0.22]{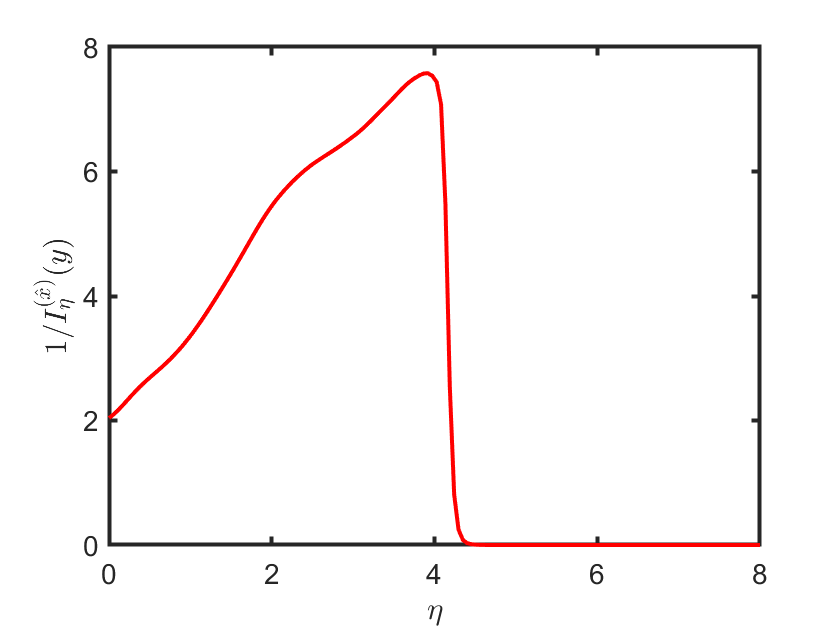}
}
\subfigure[$\epsilon_0=1, y=y_0+\epsilon_0\hat x$]{
\includegraphics[scale=0.22]{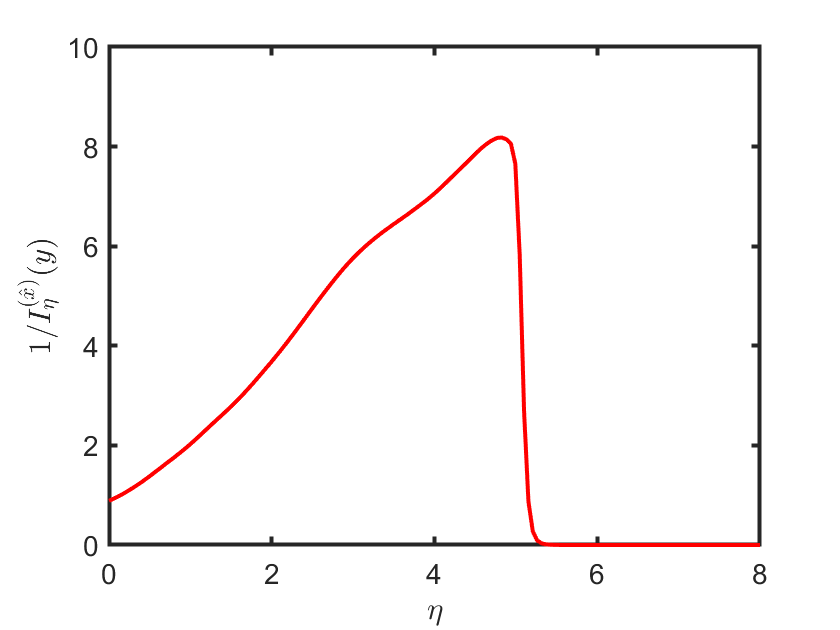}
}
\caption{ Determination of $t_{\max}$ with a fixed observation direction $\hat x=(-\sqrt{2}/2, \sqrt{2}/2)$ and  different test point $y=y_0+\epsilon_0 \hat x $ such that $\hat x\cdot y-\inf (\hat x \cdot  D)=\epsilon_0$ where we set $y_0=0.8( \cos^3 7\pi/4+\cos 7\pi/4,  \sin^3 7\pi/4+\sin 7\pi/4)$  by  plotting the one-dimensional function  $\eta\rightarrow 1/ I_{\eta}^{(\hat x)}(y)$. The radiating period of the source is $[t_{\min}, t_{\max}]=[0,4]$ for a rounded-square-shaped source support.
} \label{fig:tmax-2}
\end{figure}

Next we test various test points $y=y_0+\epsilon_0 \,\hat x$, where $\hat x\cdot y-\inf (\hat x \cdot D)=\epsilon_0=0, 0.1,1$. Here, $y_0$ is set as $0.8( \cos^3 \frac{7\pi}{4}+\cos \frac{7\pi}{4}, \sin^3 \frac{7\pi}{4}+\sin \frac{7\pi}{4})$ and $\hat x=(-\sqrt{2}/2, \sqrt{2}/2)$. We visualize the behavior of the one-dimensional function $\eta\rightarrow 1/ I_{\eta}^{(\hat x)}(y)$ in Figure \ref{fig:tmax-2}. The plot illustrates that the function $\eta\rightarrow 1/ I_{\eta}^{(\hat x)}(y)$ exhibits rapid decay on the right hand of $\eta=t_{\max}+\epsilon$, and the value of $1/ I_{\eta}^{(\hat x)}(y)$ tends to $0$ for all $\eta>t_{\max}+\epsilon_0$. Consequently, it becomes evident that the algorithm for determining $t_{\max}$ is reliable only if $\epsilon_0$ is sufficiently small (see Figures \ref{fig:tmax-2} (a) and (b)).
Figure \ref{fig:tmax-2} (c) shows that the reconstruction of $t_{\max}$ is distorted if $\epsilon_0=1$ is not a small number.

\subsubsection{Numerical examples with far-field measurements in $\R^3$} \label{sec:far}

\begin{figure}[H]
\centering
\subfigure[A spherical support]{
\includegraphics[scale=0.22]{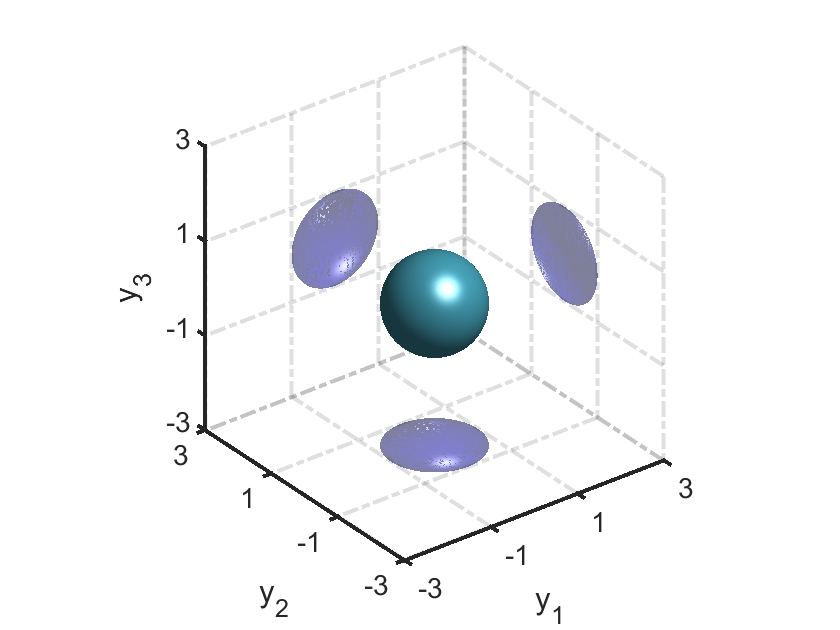}
}
\subfigure[A slice at $y_2=0$ ]{
\includegraphics[scale=0.22]{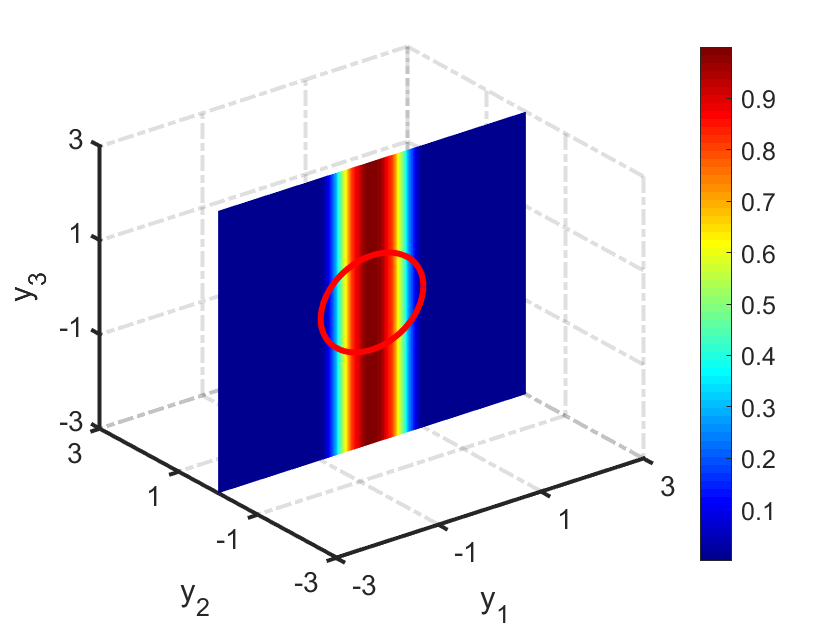}

}
\subfigure[Iso-surface and a sphere]{
\includegraphics[scale=0.22]{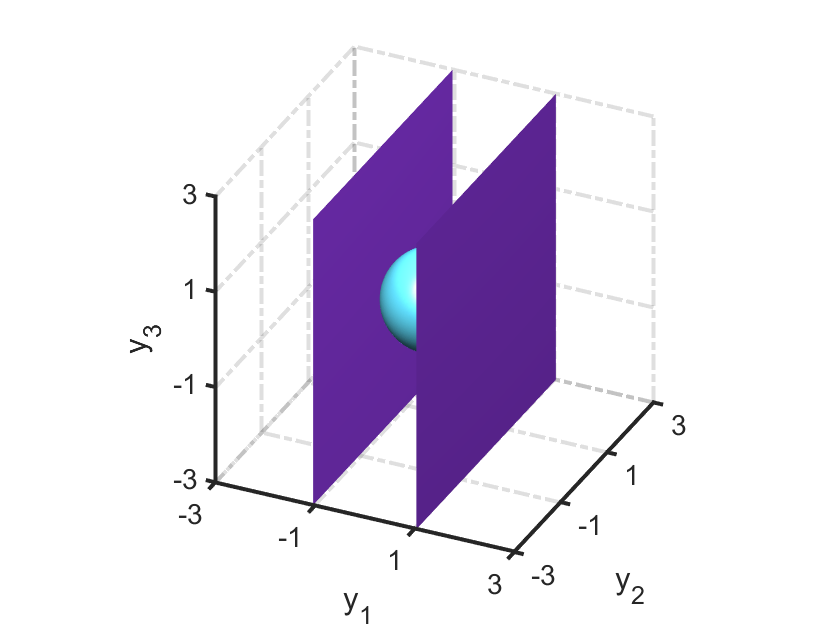}

}
\caption{Reconstructions for a spherical support using multi-frequency far-field data from one pair of opposite observation direction $(\pm1,0,0).$
} \label{fig:2-1}
\end{figure}

\begin{figure}[H]
\centering
\subfigure[A cubic support]{
\includegraphics[scale=0.22]{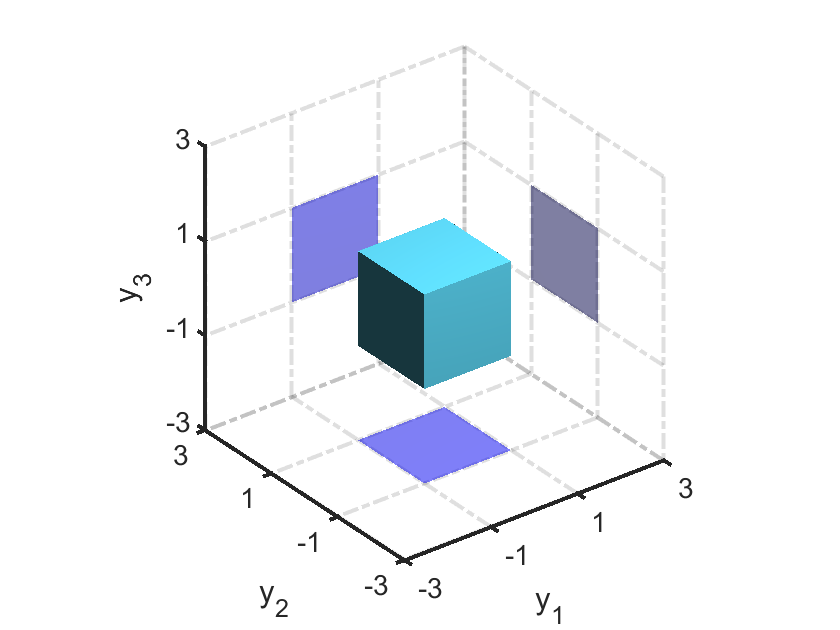}
}
\subfigure[A slice at $y_2=0$ ]{
\includegraphics[scale=0.22]{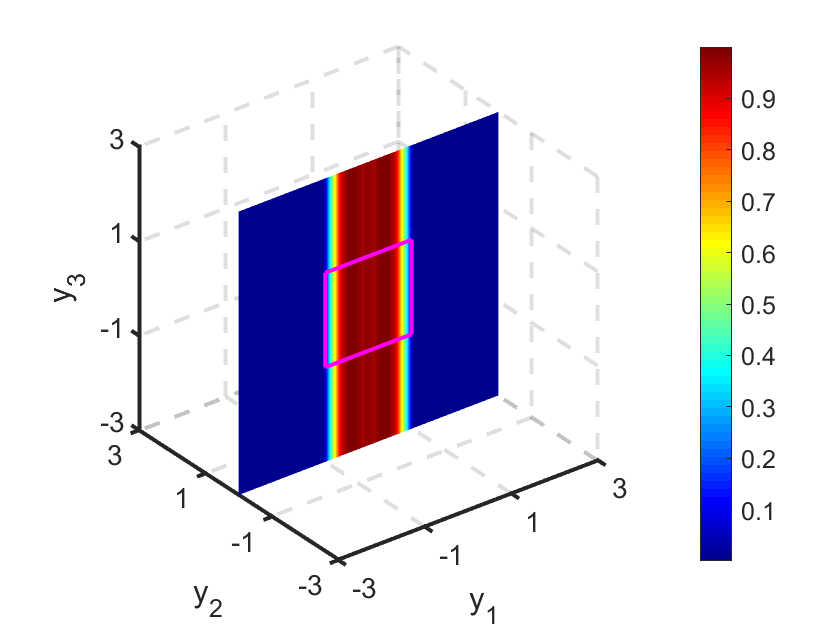}

}
\subfigure[iso-surface and a cube]{
\includegraphics[scale=0.22]{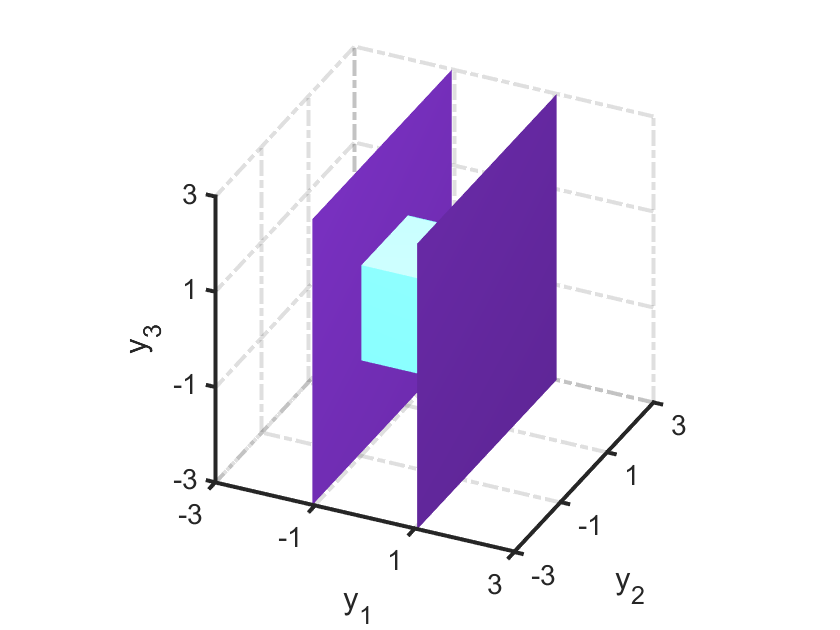}

}
\caption{Reconstructions for a cubic support using multi-frequency far-field data from one pair of opposite observation direction $(\pm1,0,0).$
} \label{fig:2-5}
\end{figure}

In this subsection, we explore the reconstructions of a spherical support and a cubic support using multi-frequency far-field measurements from single/sparse observation directions in $\R^3$, respectively.
The spherical support of the source $D$ is characterized by the set $D=\{x\in \R^3: x_1^2+x_2^2+x_3^2\leq1\}$ (see Figure \ref{fig:2-1}(a)). Meanwhile, the cubic support of the source $D$ is described by the set $D=\{x\in \R^3: |x_1|\leq1, |x_2|\leq1, |x_3|\leq1 \}$ (see Figure \ref{fig:2-5}(a)). The search domain is chosen as $[-3,-3]^3$.

Firstly,  we present the reconstructions of a spherical support from one observation $\hat x=(1,0,0)$ by plotting the indicator function $W^{(\hat x)}(y)$ in (\ref{indicator7}) with $M=1$ in Figure \ref{fig:2-1}. Figure \ref{fig:2-1}(a) shows the geometry of a spherical source support and its projections  onto  three coordinate planes. Figure \ref{fig:2-5}(b) illustrates a slice of reconstruction at $y_2=0$, from which we conclude that the smallest strip $K_D^{(\hat x)}$ containing the source and perpendicular  to the observation direction on the cross-section at $y_2=0$ is perfectly reconstructed. The sphere is precisely sandwiched between two hyperplanes/iso-surfaces  in Figure \ref{fig:2-1}(c). In fact, we obtain a smallest slab containing the source support and perpendicular to the observation direction from a single pair of opposite observation directions multi-frequency data in $\R^3$.


\begin{figure}[H]
\centering
\subfigure[$M=3$]{
\includegraphics[scale=0.3]{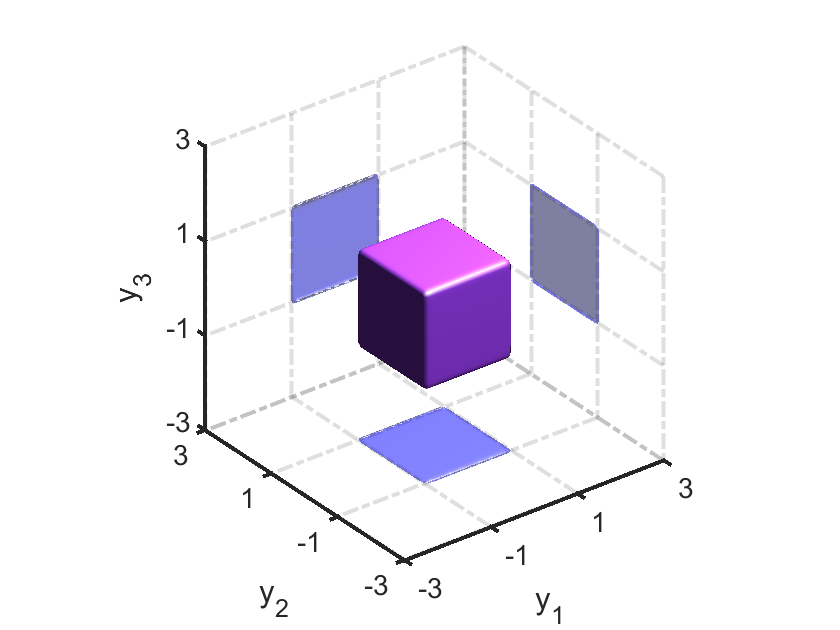}
}
\subfigure[$M=7$ ]{
\includegraphics[scale=0.3]{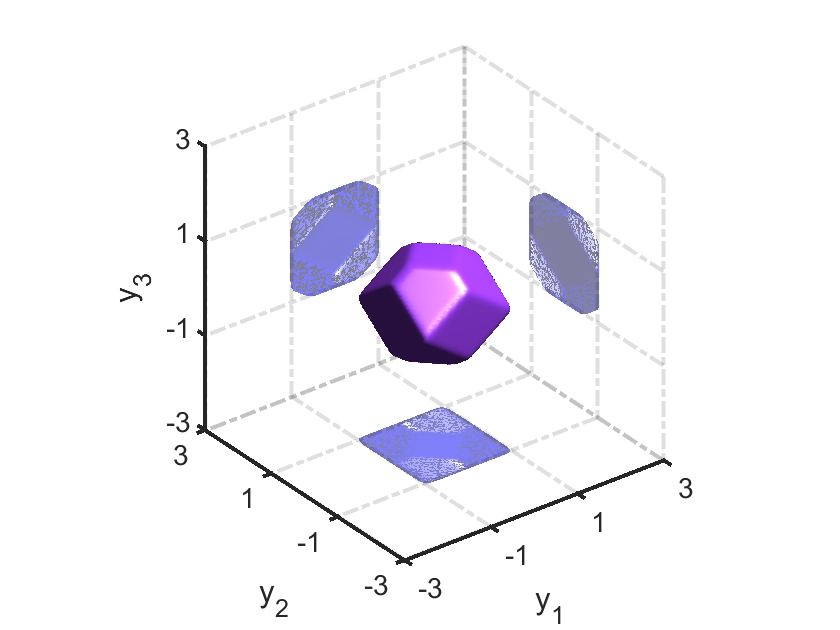}

}
%

\subfigure[$M=10$ ]{
\includegraphics[scale=0.3]{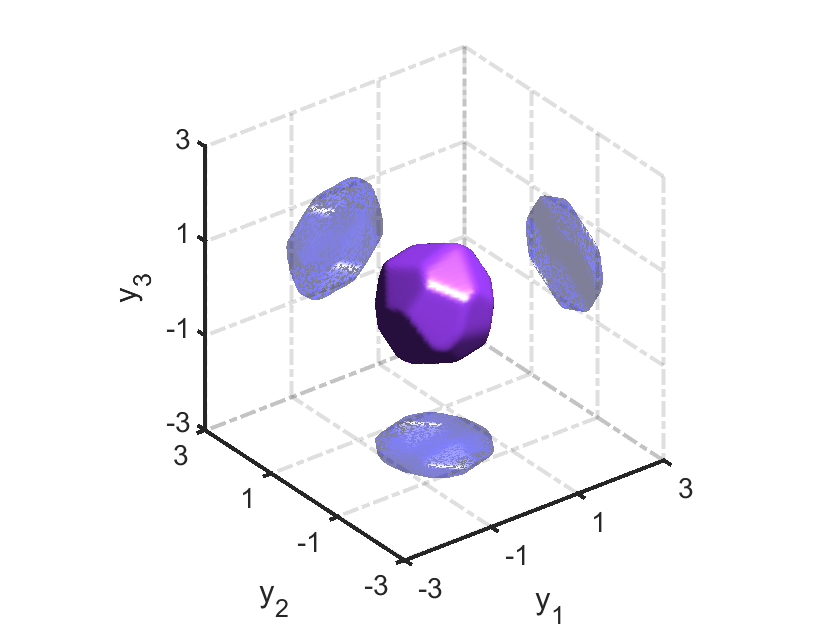}

}
\subfigure[$M=15$]{
\includegraphics[scale=0.3]{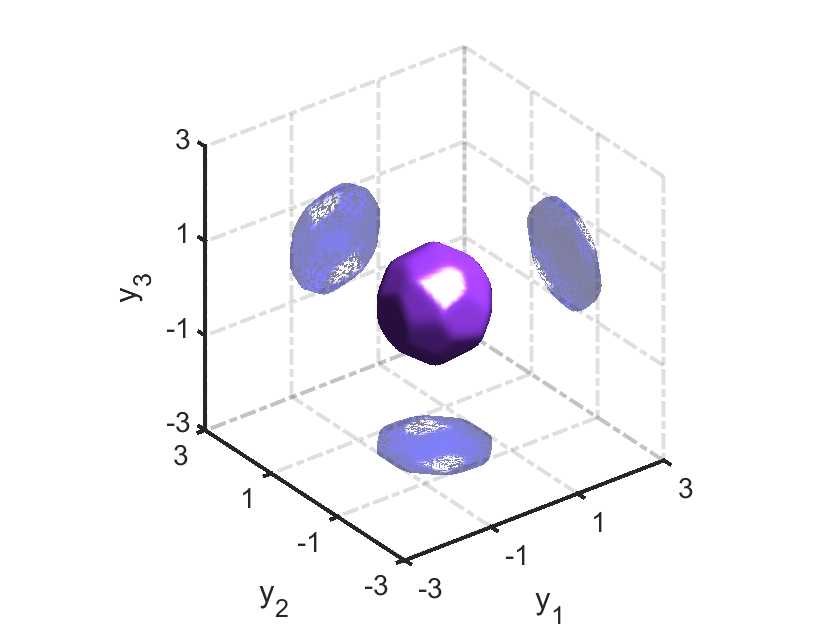}

}
%
\caption{Reconstructions for a spherical support using multi-frequency far-field data from multiple observation directions. Here we take the number of pairs of opposite observation directions $M=3,7, 10,15$.
} \label{fig:2-3}
\end{figure}

\begin{figure}[H]
\centering
\subfigure[A slice at $y_1=0$]{
\includegraphics[scale=0.22]{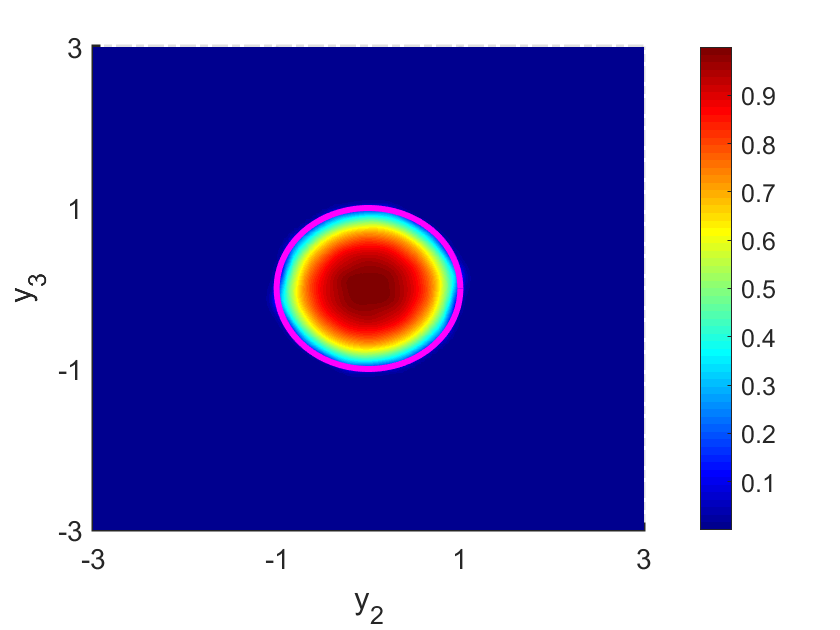}
}
\subfigure[A slice at $y_2=0$ ]{
\includegraphics[scale=0.22]{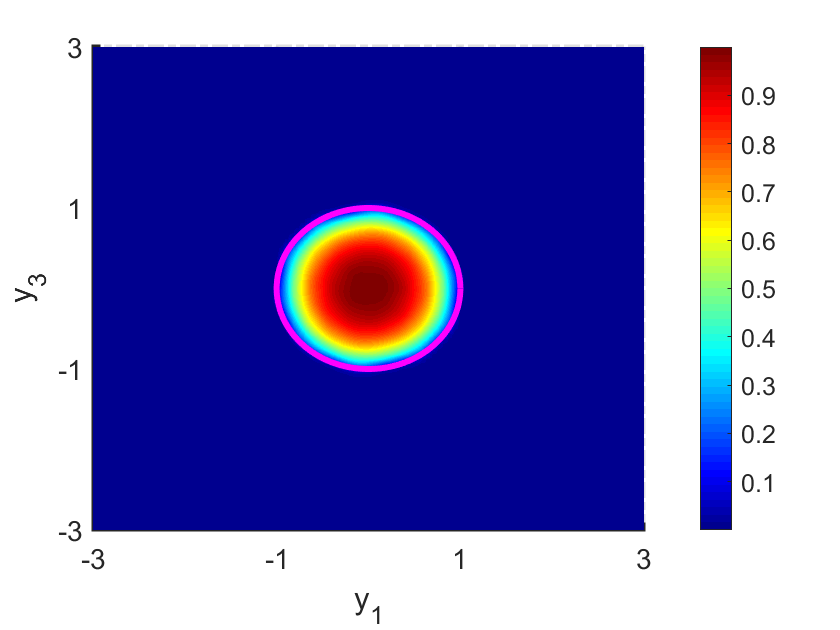}

}
\subfigure[A slice at $y_3=0$]{
\includegraphics[scale=0.22]{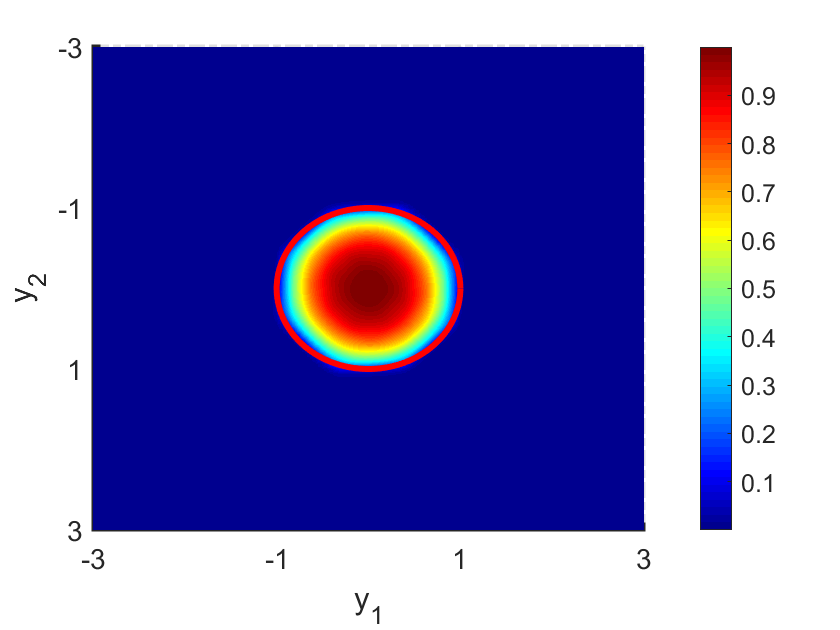}

}
\caption{Slices of reconstruction for a spherical support using multi-frequency far-field data from $15$ pairs of opposite observation directions.
} \label{fig:2-4}
\end{figure}

Secondly, we present the numerical results for reconstructing the location and shape of the spherical source support from multi-frequency data obtained through sparse observation directions, as illustrated in Figure \ref{fig:2-3}. We consider different scenarios with varying number of pairs of opposite observation directions, specifically $M=3, 7, 10, 15$.
To provide a comprehensive insight into the reconstructions presented in Figure \ref{fig:2-3}, projections onto the $oy_1y_2$, $oy_1y_3$, and $oy_2 y_3$ planes are included.  In Figure \ref{fig:2-3}(a), the reconstruction is derived from $3$ properly-selected pairs of opposite directions $\hat x=\{(\pm1,0,0), (0,\pm1,0), (0,0,\pm1)\}$. It is evident that the spherical support is accurately embedded within the cube, consistent with our theoretical predictions. Indeed, we theoretically establish the location of the source support using far-field data from any three different observation directions in $\R^3$.
Figure \ref{fig:2-3}(b) illustrates the reconstruction from a set of  7 properly-selected pairs of opposite observation direction  $\hat x=\{(1,0,0), (0,1,0), (0,0,1), (\sqrt {2}/2, \sqrt {2}/2, 1/2), (\sqrt {2}/2, -\sqrt {2}/2, 1/2), (-\sqrt {2}/2, \sqrt {2}/2, 1/2),$ $(-\sqrt {2}/2, -\sqrt {2}/2, 1/2) \}$ and the corresponding $-\hat x$.  For Figures \ref{fig:2-3}(d)-(e), we choose 10 and 15 observation directions uniformly distributed on the upper hemisphere with a radius of 1. Clearly, as the number of pairs of opposite observation directions $M$ increases, the inversion results become progressively closer to the spherical support.  Upon examining these 2D visualizations,  a noteworthy transformation is observed in the projections onto the three coordinate planes, evolving from a square $[-1,1]^2$  in Figure \ref{fig:2-4}(a) to gradually assume an approximate circular shape $y_i^2+y_j^2\leq1, i\neq j, i,j=1,2,3$ in Figure \ref{fig:2-4}(d).
These findings not only confirm the accuracy of our three-dimensional reconstructions but also emphasize the significance of utilizing multiple observation directions for enhanced precision.

\begin{figure}[H]
\centering
\subfigure[$M=3$]{
\includegraphics[scale=0.3]{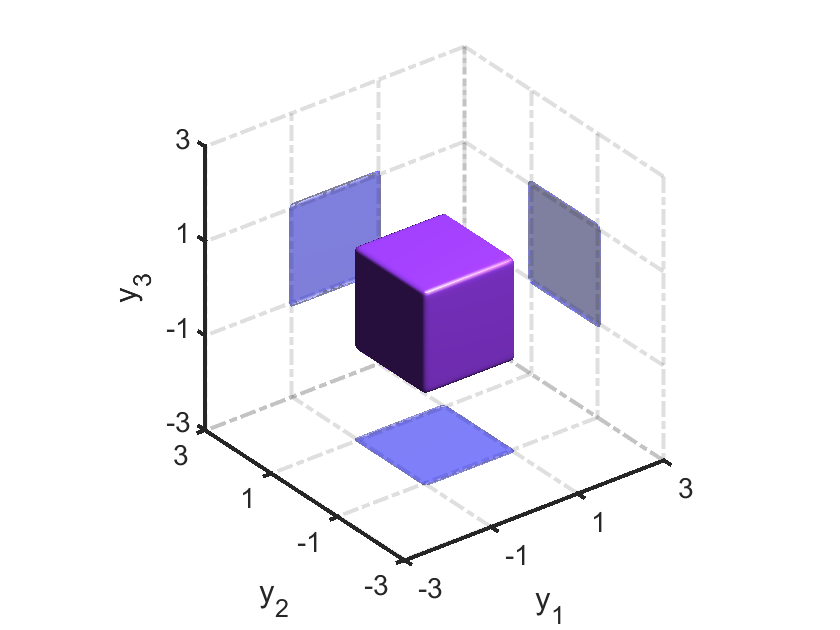}
}
\subfigure[$M=7$ ]{
\includegraphics[scale=0.3]{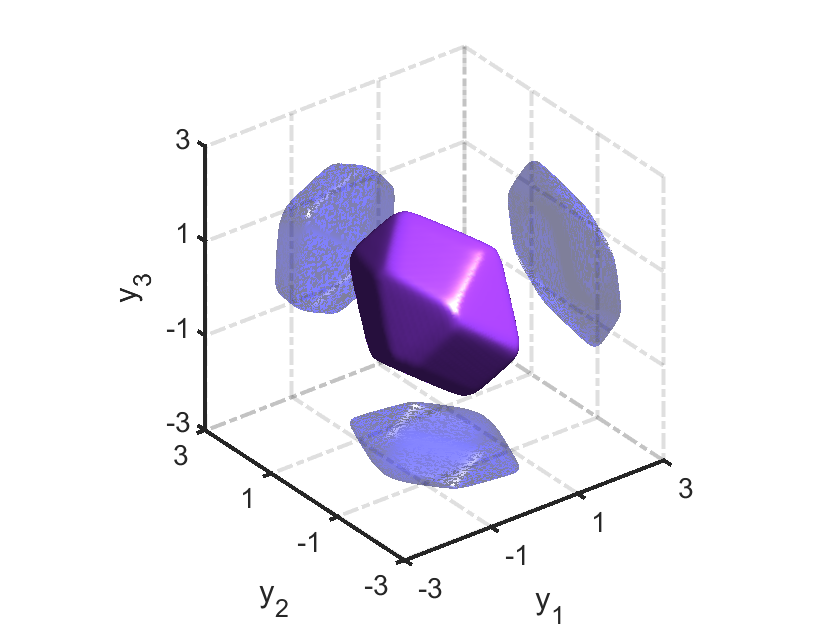}

}
\subfigure[$M=10$]{
\includegraphics[scale=0.3]{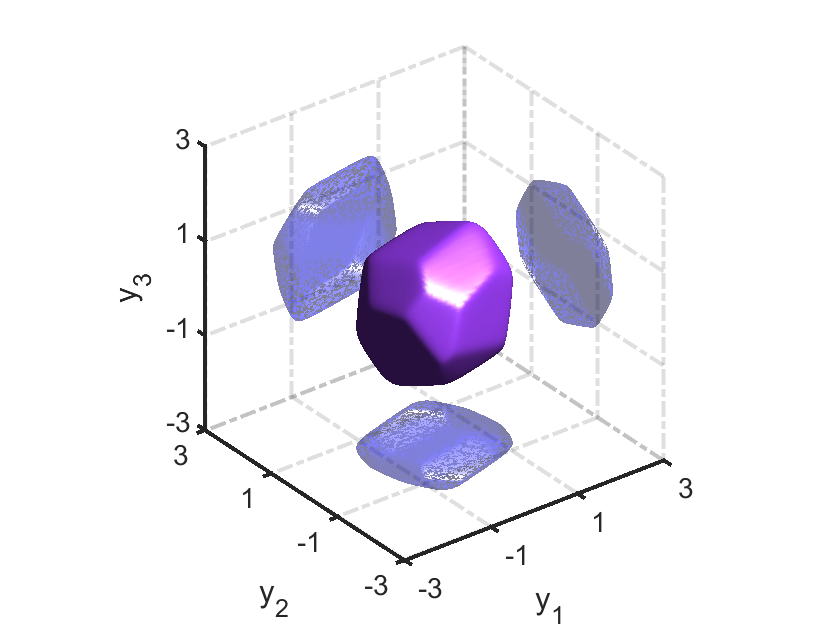}

}
\subfigure[$M=15$ ]{
\includegraphics[scale=0.3]{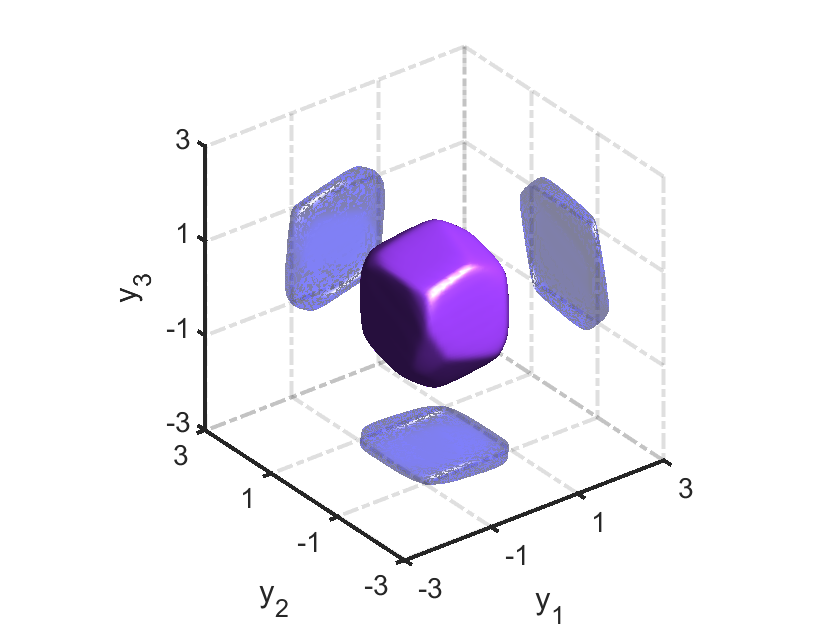}

}
%
\caption{Reconstructions for a cubic support using multi-frequency far-field data from a finite number of pairs of opposite observation directions with $M=3,7, 10,15$.
} \label{fig:2-6}
\end{figure}

\begin{figure}[H]
\centering
\subfigure[$y_1=0$]{
\includegraphics[scale=0.22]{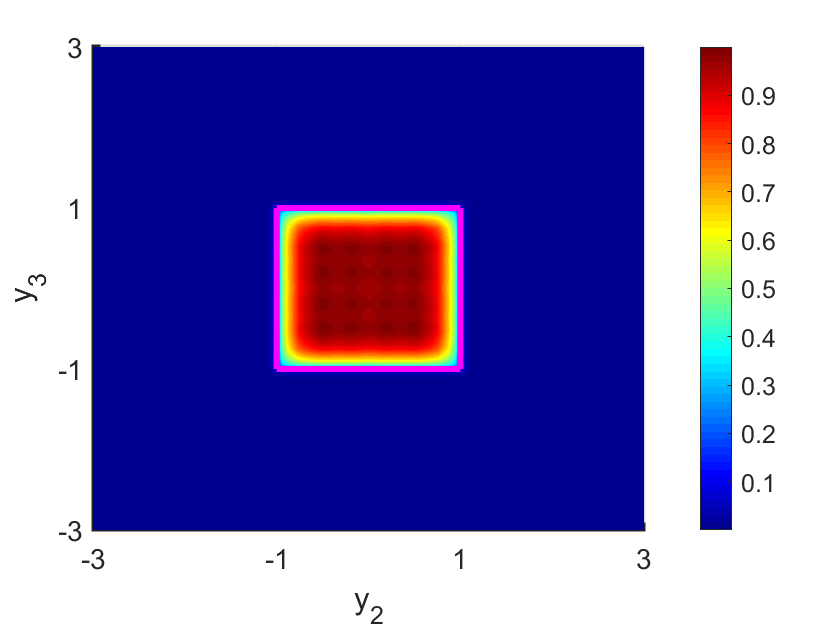}
}
\subfigure[$y_2=0$ ]{
\includegraphics[scale=0.22]{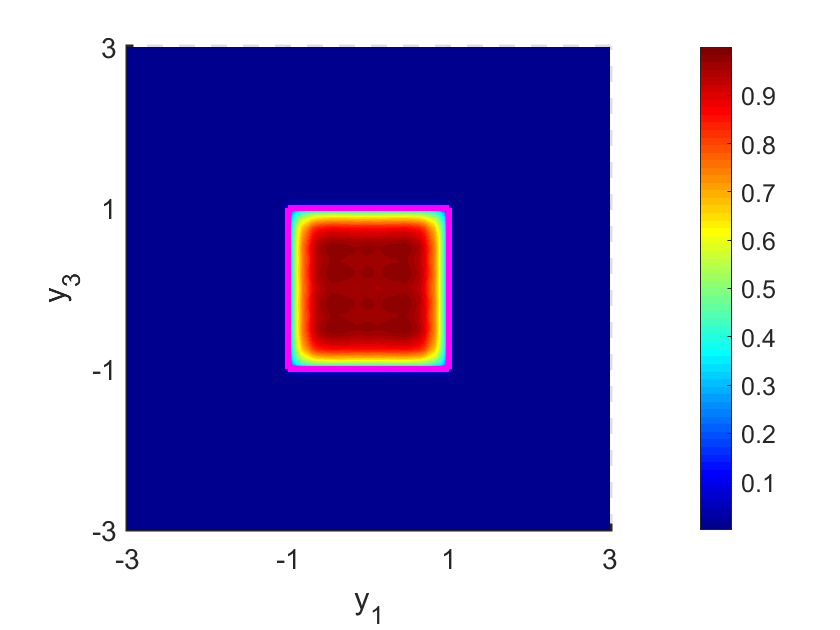}

}
\subfigure[$y_3=0$]{
\includegraphics[scale=0.22]{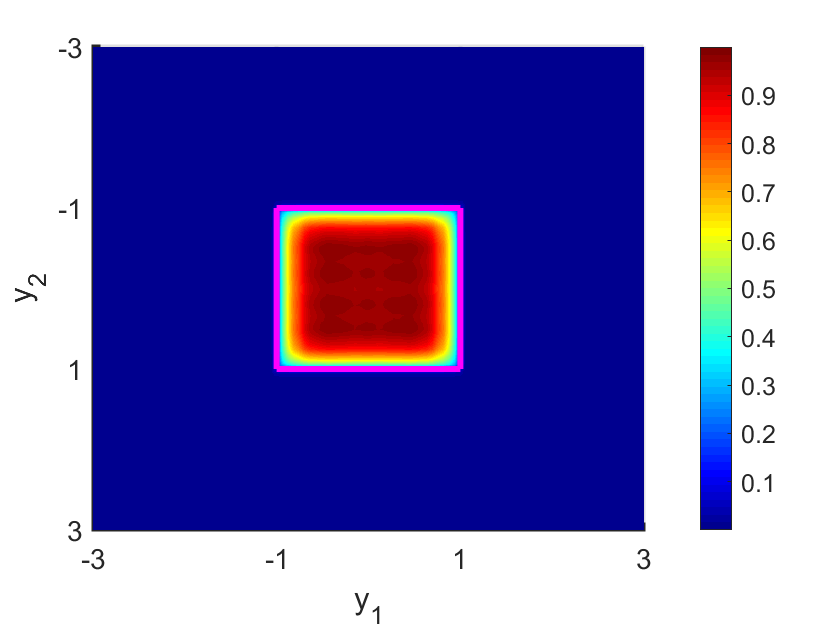}

}
\caption{Slices of reconstruction for a  cubic support using multi-frequency far-field data from $3$ pairs of  properly-selected opposite observation directions.
} \label{fig:2-7}
\end{figure}

Next, we turn our attention to the reconstructions for a cubic support, as illustrated in Figure \ref{fig:2-5}(a)). Figure \ref{fig:2-5}(b) and (c) depict a slice at $y_2=0$  and an iso-surface by plotting the indicator function $W(y)$ with data from  a single pair of opposite observation direction $(\pm1,0,0)$. It is evident that  the smallest strip $K_D^{(\hat x)}$ is well-represented in Figure \ref{fig:2-5}(b) and the cube is precisely enclosed between two hyperplanes/iso-surfaces in Figure \ref{fig:2-5}(c).
Given the exact geometry of the source support being cubic, theoretically, the geometry can be precisely recovered using far-field data from three properly-selected pairs of opposite observation directions. We showcase the reconstruction using far-field data from $3$ pairs of directions $\hat x=\{(\pm1,0,0),(0,\pm1,0),(0,0,\pm1)\}$ in Figure \ref{fig:2-6}(a), where the projections onto the three coordinate planes are also included. The results demonstrate that the cubic source support is perfectly reconstructed. Slices of the reconstruction are further presented in Figure \ref{fig:2-7}. Additionally, in Figure \ref{fig:2-6}(b)-(d), we utilize multi-frequency data from $7, 10, 15$ pairs of opposite observation directions to reconstruct the cubic source support, respectively, with  directions  uniformly distributed on the sphere with a radius of $1$. A common observation is that both the location and shapes are accurately captured as the number of observation directions increases.

\subsection{Numerical examples for near-field case in $\R^3$}
\begin{figure}[H]
\centering
\subfigure[A slice $y_2=0$ for $1/\tilde I_{\eta,\epsilon}^{(-x)}$]{
\includegraphics[scale=0.22]{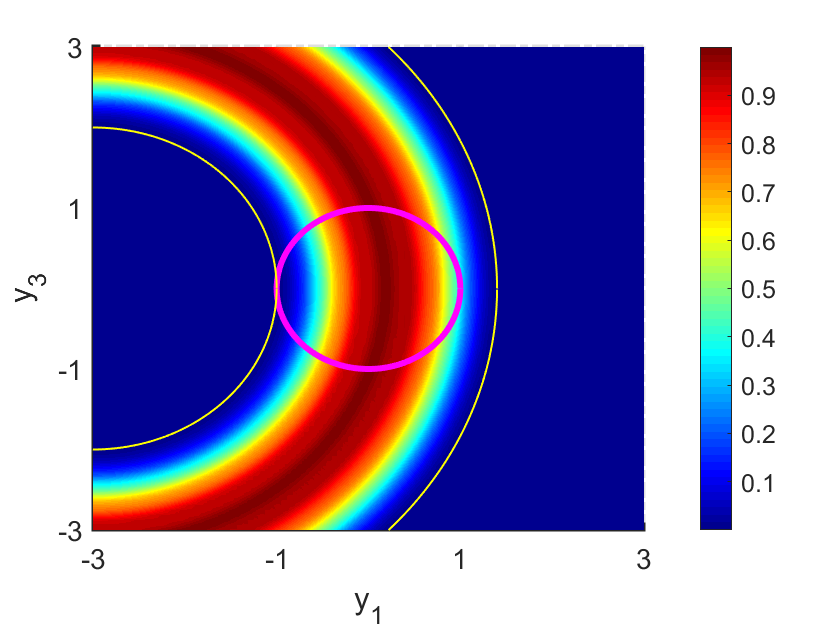}
}
\subfigure[A slice $y_2=0$ for $1/\tilde I_{\eta,\epsilon}^{(x)}$]{
\includegraphics[scale=0.22]{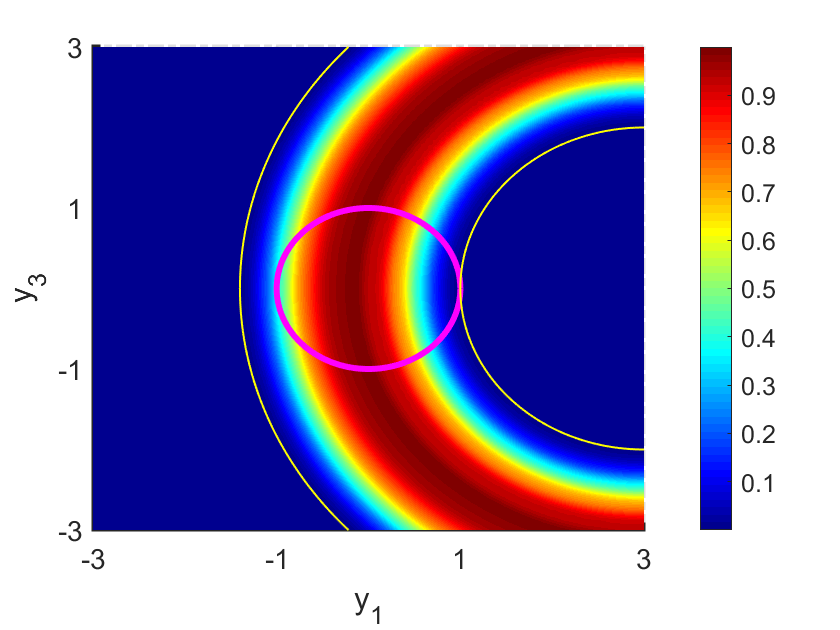}

}
\subfigure[A slice $y_2=0$ for $\widetilde W^{(x)}$  ]{
\includegraphics[scale=0.22]{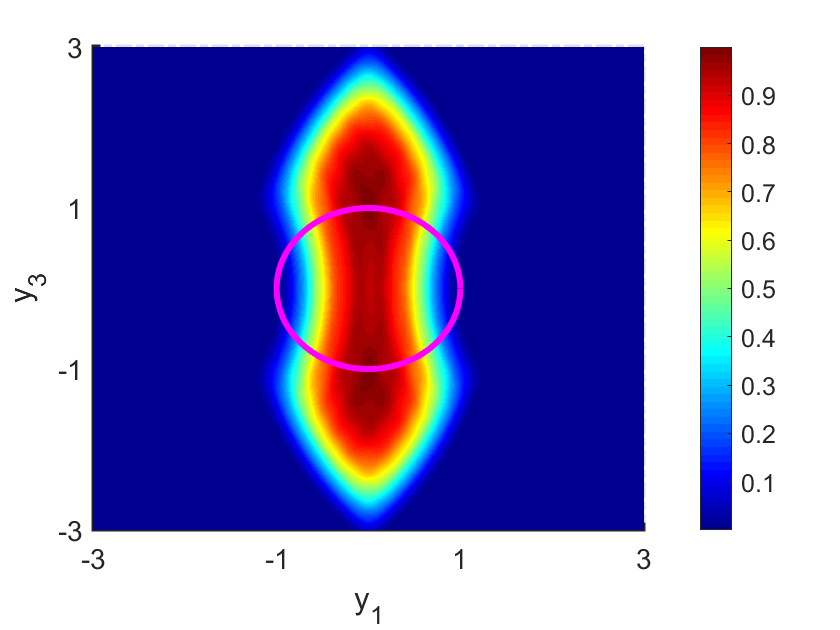}

}
\subfigure[Iso-surface for $1/\tilde I_{\eta,\epsilon}^{(-x)}$]{
\includegraphics[scale=0.22]{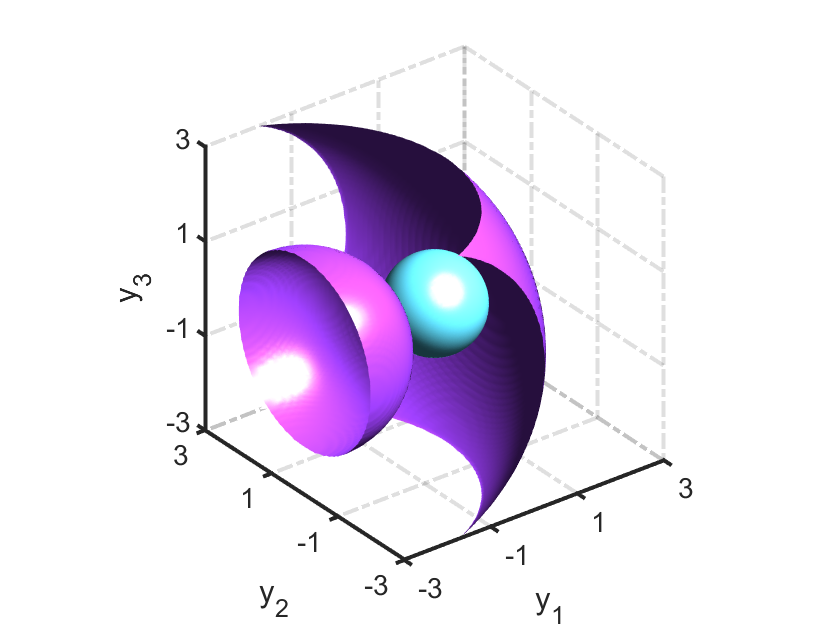}
}
\subfigure[Iso-surface for $1/\tilde I_{\eta,\epsilon}^{(x)}$]{
\includegraphics[scale=0.22]{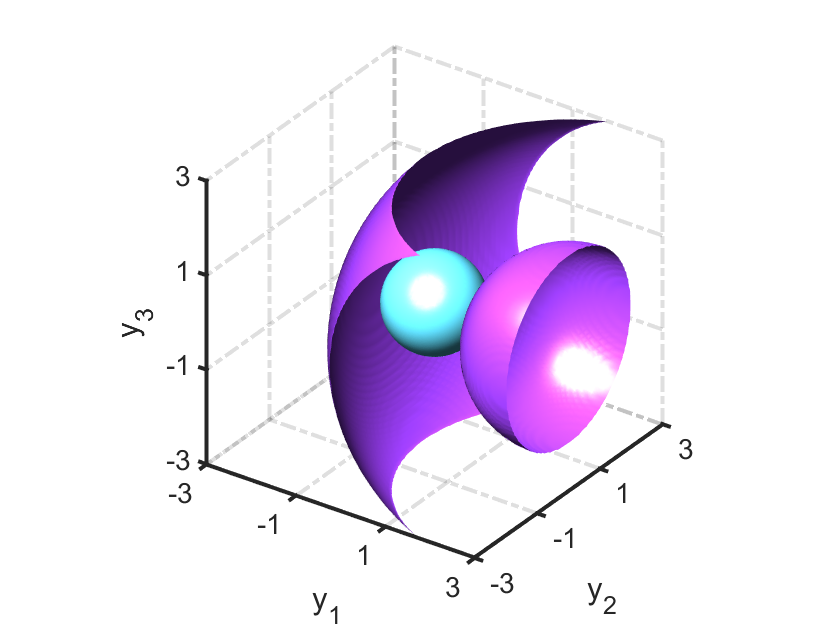}

}
\subfigure[Iso-surface for $\widetilde W^{(x)}$]{
\includegraphics[scale=0.22]{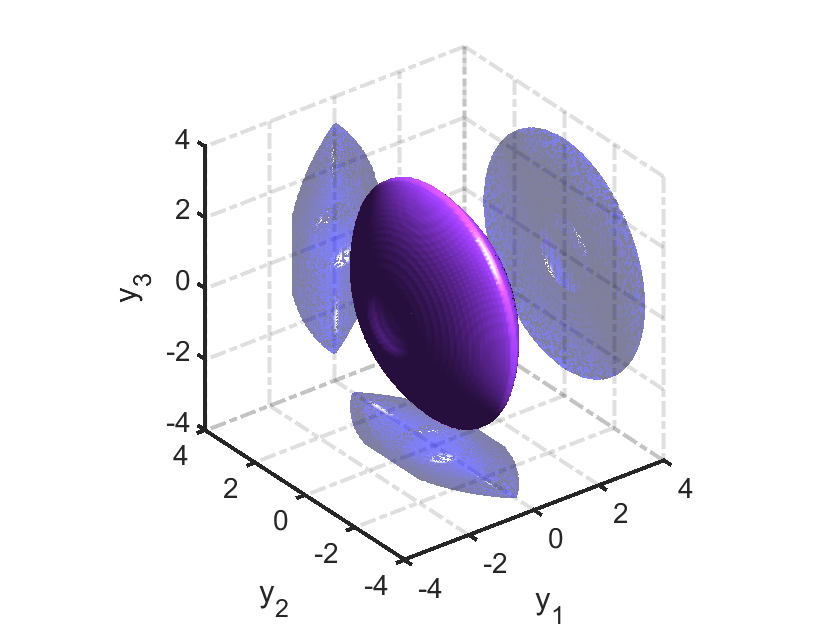}

}
\caption{Reconstructions for a spherical support using  multi-frequency near-field data from a single observation point $x=(3,0, 0)$ and $-x$ and
from a pair of opposite observation points $(\pm 3,0, 0)$.
} \label{fig:3-1}
\end{figure}

In this subsection, we  continue to conduct  the numerical examples in Section \ref{sec:far} for reconstructing the spherical-shaped support  and cubic-shaped source support but with near-field measurements.
Being different from the far-field case, we need to reconstruct the annulus $ A^{(x)}_{D,\eta}$ which contains the source support $D$ and centered at the observation point $x$ by using the auxiliary indicator function $\tilde I_{\eta,\epsilon}^{(x)}(y)$ defined in(\ref{indic-near}).
In Figure \ref{fig:3-1}(a) and (d), (b) and (e), we present a slice $y_2=0$  and iso-surfaces of reconstructions for a spherical support from one observation point $x=(3,0,0)$ by plotting the indicator function $1/\tilde I_{\eta,\epsilon}^{(x)}(y)$ and $1/\tilde I_{\eta,\epsilon}^{(-x)}(y)$ . To enhance the resolution, the boundaries  of the sphere  and the annulus  $ A^{(x)}_{D,\eta}$ at  slice  $y_2=0$ are highlighted in pink and yellow solid lines in Figure \ref{fig:3-1}(a) and (b), respectively.  It is evident that $\{y\in \R^3:|x-y|=\inf_{z\in D}|x-z|\}$ and $\{y\in \R^3:|x- y|=\sup_{z\in D}|x-z|+t_{\max}-\eta\}$  (yellow lines), which is  with a shift $t_{\max}-\eta$ to $\{y\in \R^3:|x-y|=\sup_{z\in D}|x-z|\}$,  are all precisely reconstructed.
The spherical support is enclosed between the iso-surfaces and  closely adjacent to the inner ball surface in Figures \ref{fig:3-1} (d) and (e). We further illustrate a slice $y_2=0$ and iso-surfaces of reconstruction from a pair of  opposite observation points $(\pm3,0,0)$ by plotting the indicator function $\widetilde W^{(x)}(y)$ (\ref{indic-near-1})  in Figures \ref{fig:3-1} (c) and (e). The spherical support is enclosed within the disc-shaped geometric object and the location has also been roughly determined in Figure \ref{fig:3-1}(f). It is worth noting that we can not reconstruct the $A^{(x)}_D$ defined in (\ref{annulus-a}) even using a pair of opposite observation points $\pm x$, which differs from the far-field case.

\begin{figure}[H]
\centering
\subfigure[$M=3$]{
\includegraphics[scale=0.3]{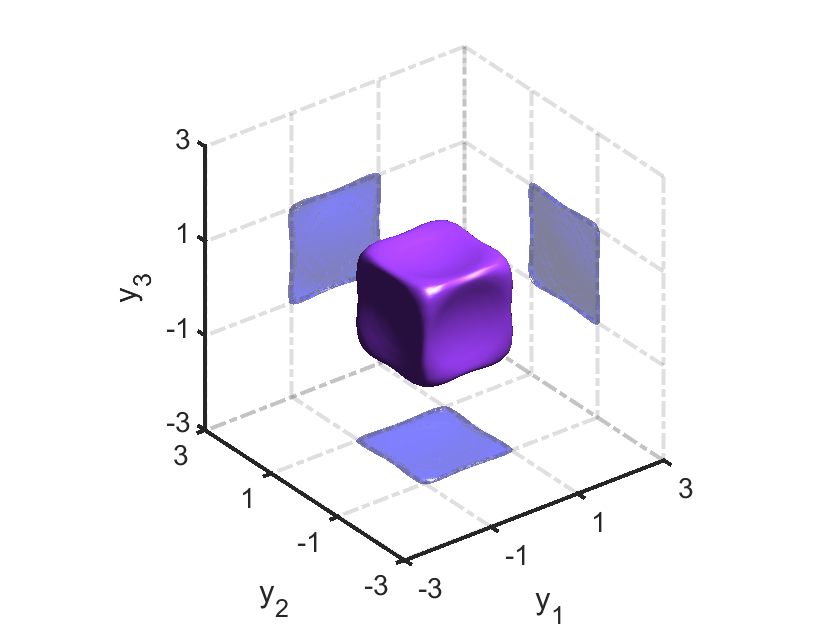}
}
\subfigure[$M=7$ ]{
\includegraphics[scale=0.3]{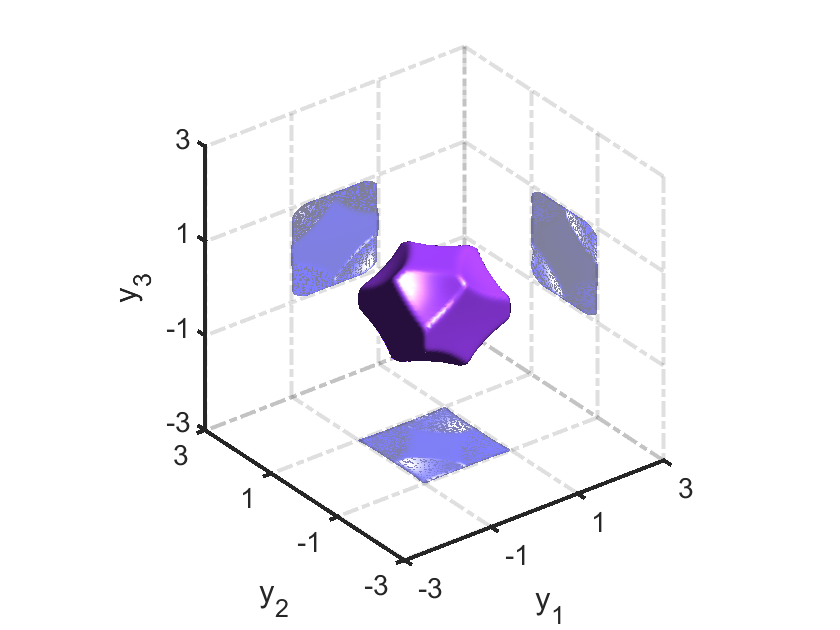}

}
\subfigure[$M=10$]{
\includegraphics[scale=0.3]{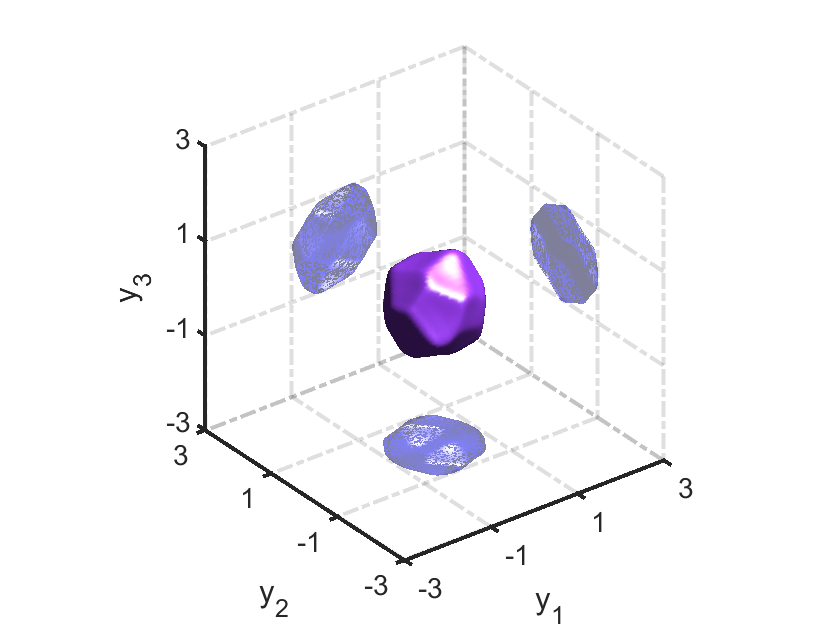}

}
\subfigure[$M=15$ ]{
\includegraphics[scale=0.3]{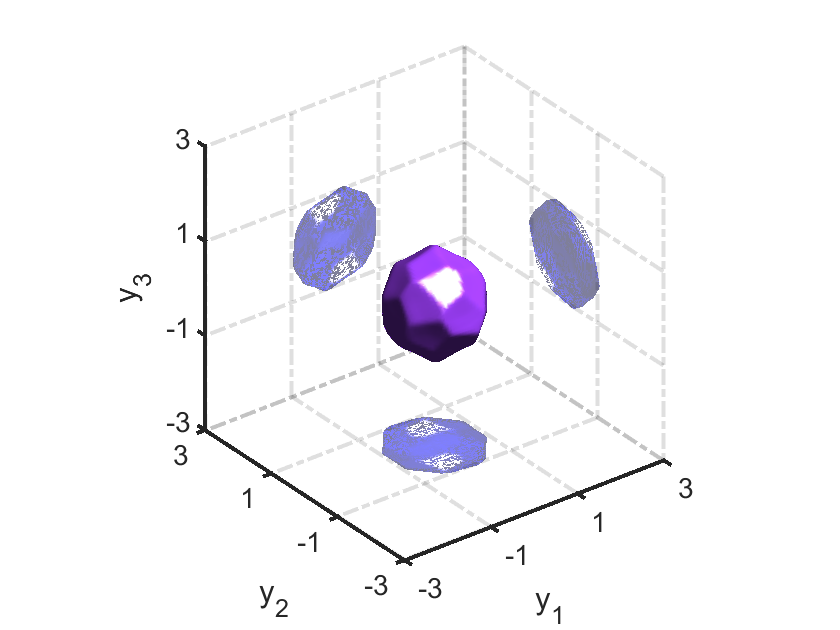}

}
\caption{Reconstructions  for  a spherical support using  multi-frequency near-field data at sparse observation points. Here we choose the number of observation points $M=3,7,10,15$.
} \label{fig:3-2}
\end{figure}

\begin{figure}
\centering
\includegraphics[scale=0.4]{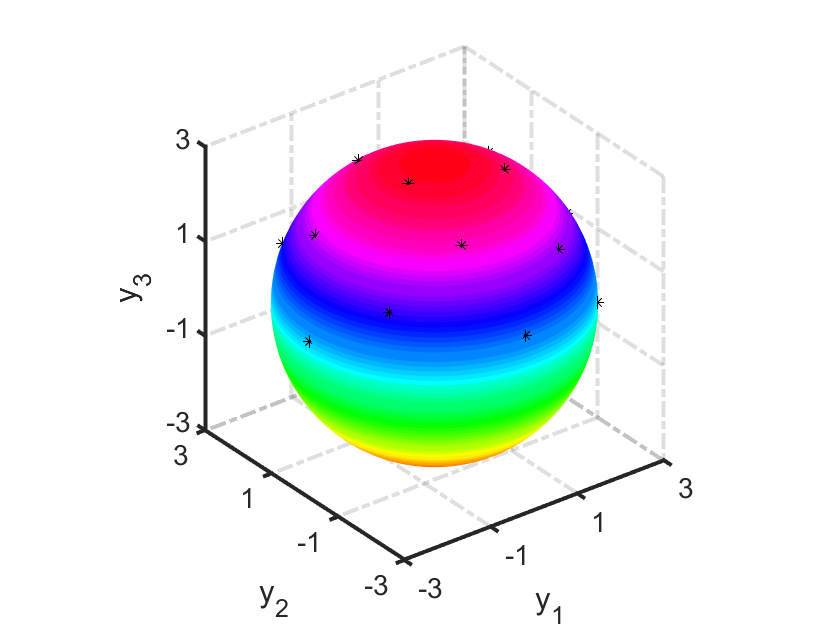}
\caption{15 points uniformly distributed on the upper hemisphere with a radius of $3$.}\label{fig:3-22}
\end{figure}

\begin{figure}[H]
\centering
\subfigure[$y_1=0$]{
\includegraphics[scale=0.22]{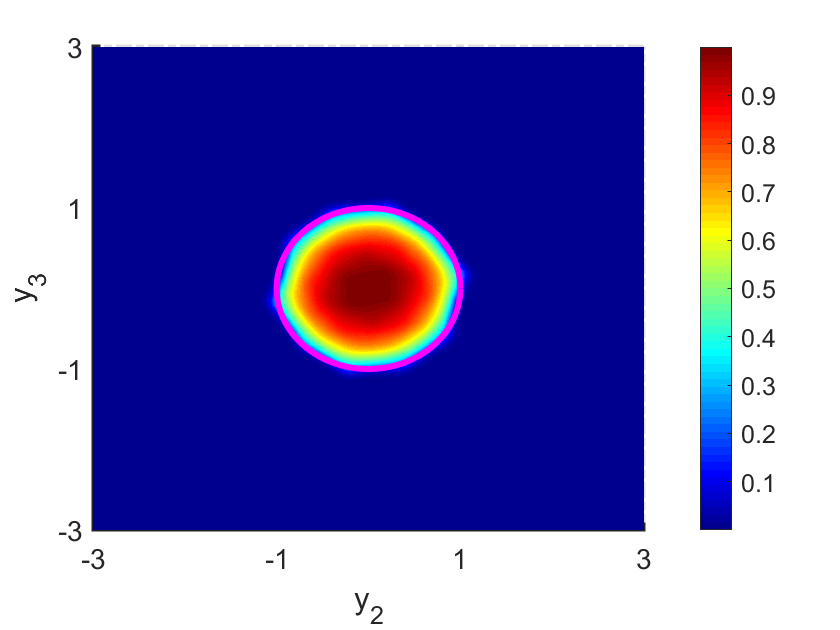}
}
\subfigure[$y_2=0$ ]{
\includegraphics[scale=0.22]{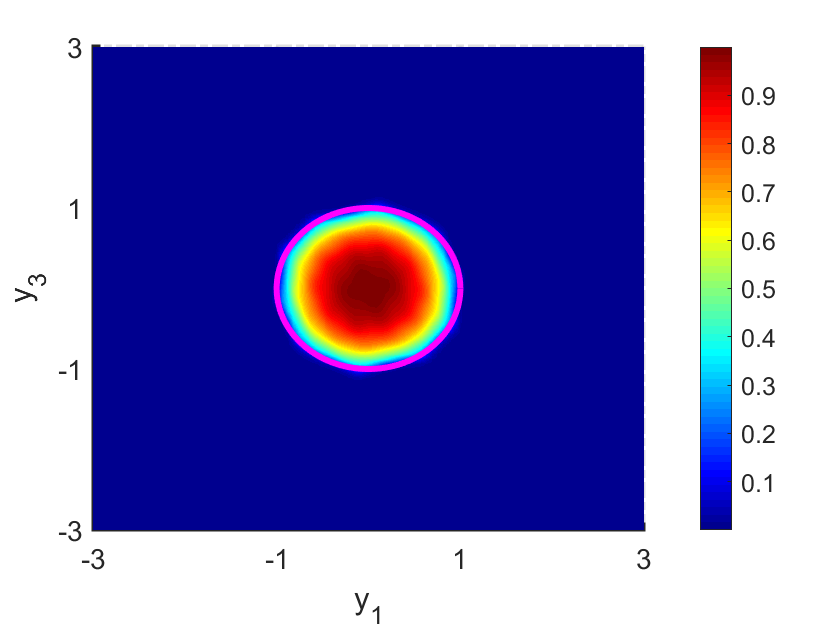}

}
\subfigure[$y_3=0$]{
\includegraphics[scale=0.22]{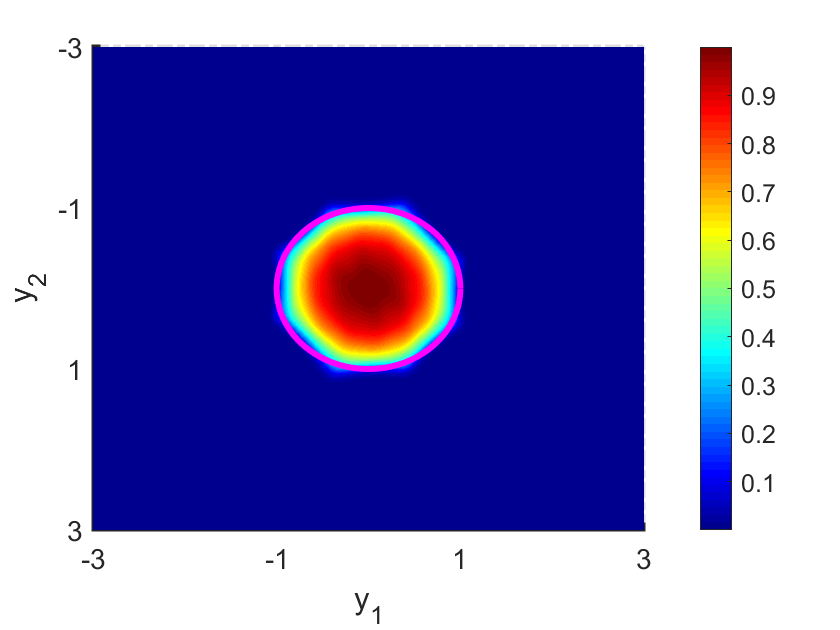}

}
\caption{Slices of reconstructions  for  a spherical support using  multi-frequency near-field data at $15$ pairs of opposite observation points.
} \label{fig:3-3}
\end{figure}

Figure \ref{fig:3-2} illustrates the reconstructions of a spherical source support from sparse observation points.  In Figure \ref{fig:3-2}(a), the precise location is obtained using near-field data from 3 pairs of observation points $x=\{(\pm 3,0,0), (0, \pm 3,0),(0,0,\pm 3)\}$. In Figure \ref{fig:3-2}(b), (c) and (d) we choose $7, 14,15$ observation points uniformly distributed on the upper hemisphere with a radius of $3$ together with their symmetric observation points, respectively. It is evident that with  increasing number of observation points, the reconstruction results  progressively approach the shape of a sphere.  The projections onto the three coordinate planes further verify the accuracy of our algorithm. Figure \ref{fig:3-22} depicts  15 points uniformly distributed on the upper hemisphere with a radius of 3. Additionally, Figure \ref{fig:3-3} shows slices of the reconstruction at the planes $y_1 = 0$, $y_2 = 0$ and $y_3 = 0$ using  data from 15 pairs of opposite observation points (see Figure \ref{fig:3-2}(d)).  For comparison purpose, the boundary of the support's slice is also demonstrated with the pink solid line. These slices provide further confirmation of the accuracy of our algorithm.

\begin{figure}[H]
\centering
\subfigure[$M=3$]{
\includegraphics[scale=0.3]{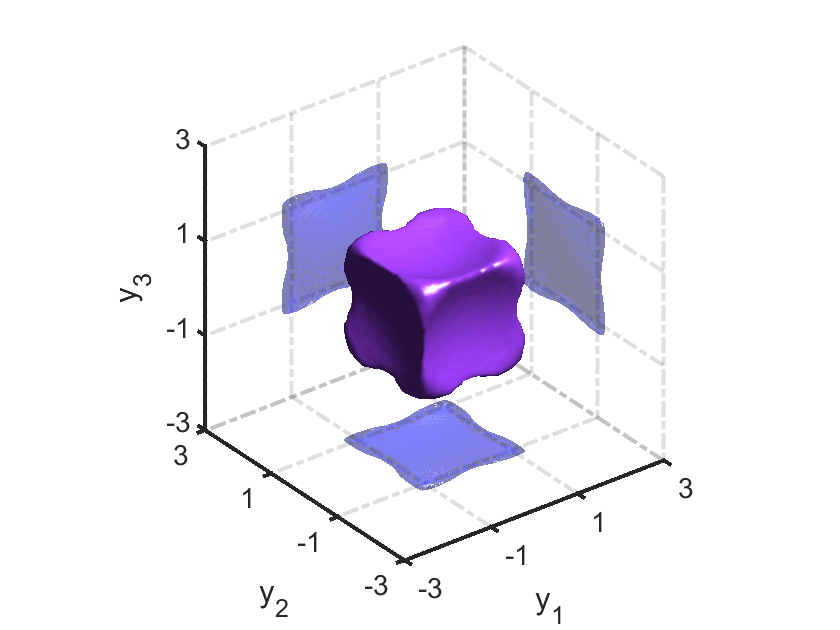}
}
\subfigure[$M=7$ ]{
\includegraphics[scale=0.3]{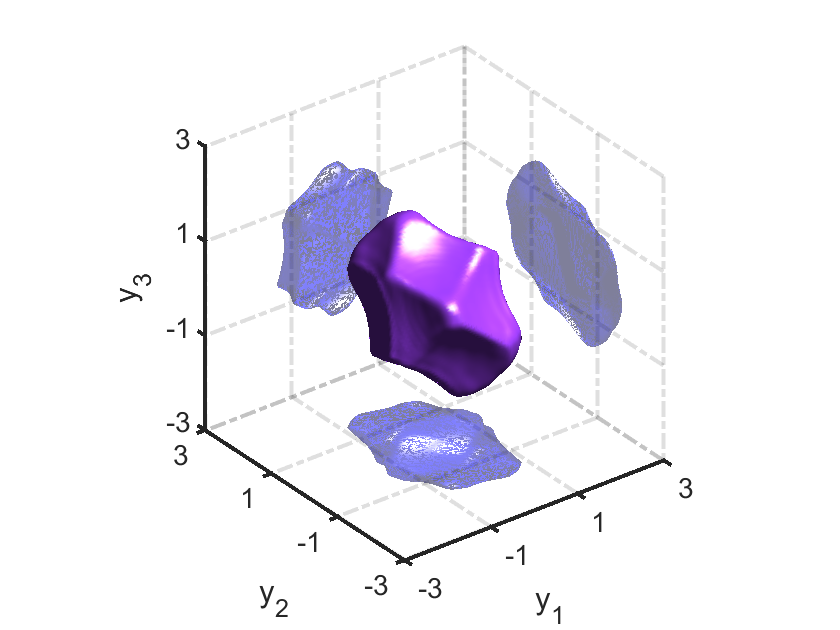}

}
\subfigure[$M=10$]{
\includegraphics[scale=0.3]{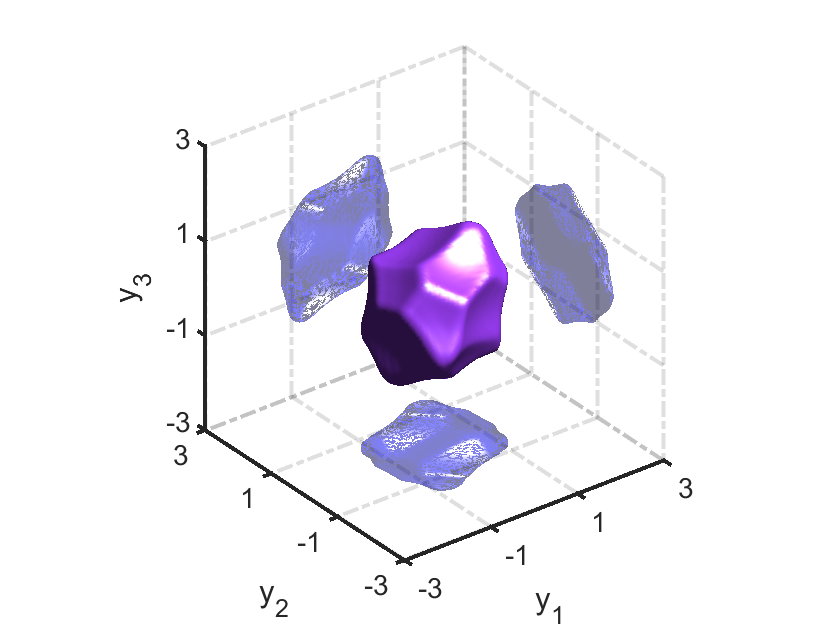}

}
\subfigure[$M=15$ ]{
\includegraphics[scale=0.3]{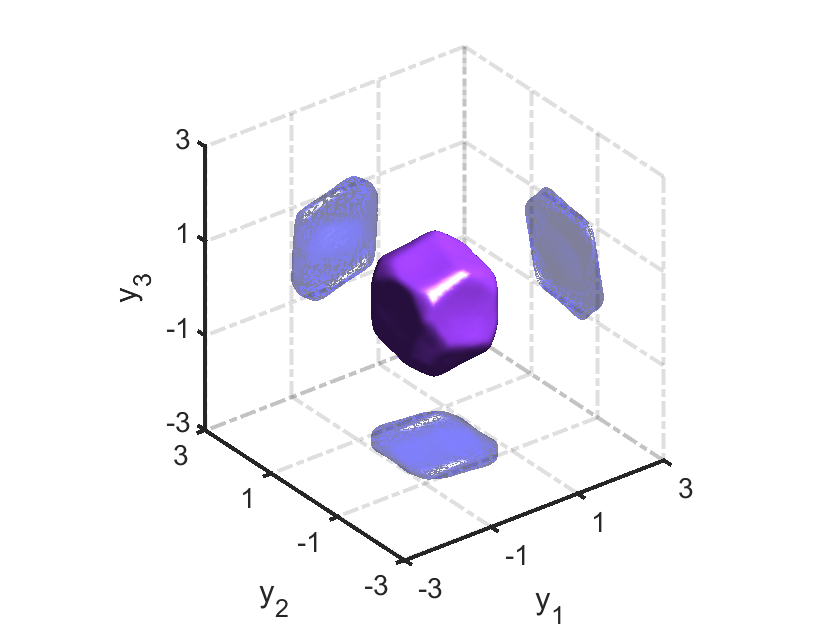}

}
\caption{Reconstructions  for  a cubic support using  multi-frequency near-field data at sparse observation points. Here we choose the number of pairs of opposite observation points $M=3, 7, 10, 15$.
} \label{fig:3-4}
\end{figure}

\begin{figure}[H]
\centering
\subfigure[$y_1=0$]{
\includegraphics[scale=0.22]{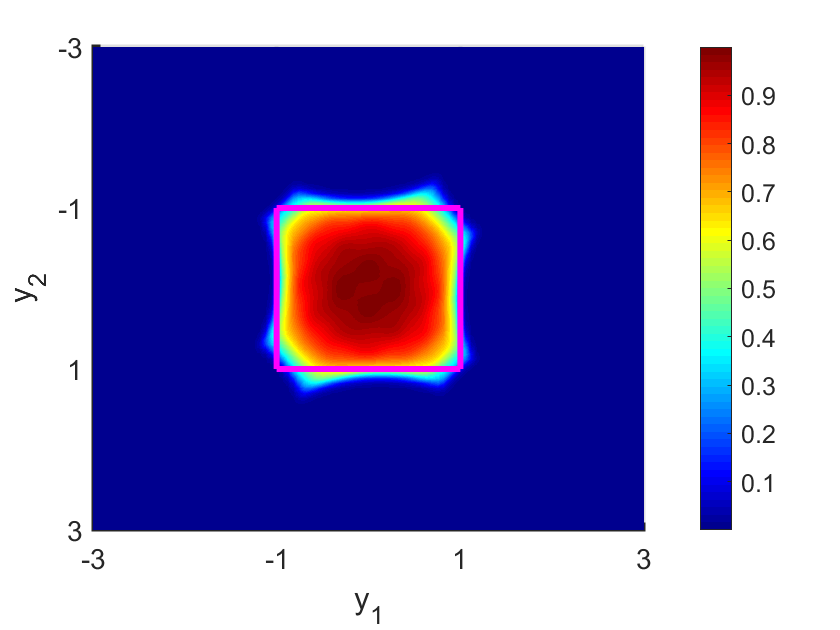}
}
\subfigure[$y_2=0$ ]{
\includegraphics[scale=0.22]{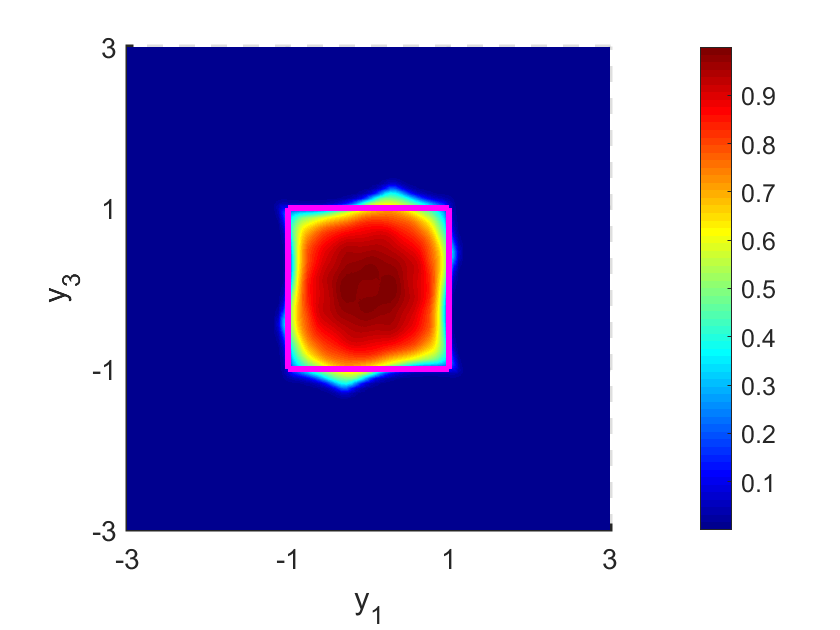}

}
\subfigure[$y_3=0$]{
\includegraphics[scale=0.22]{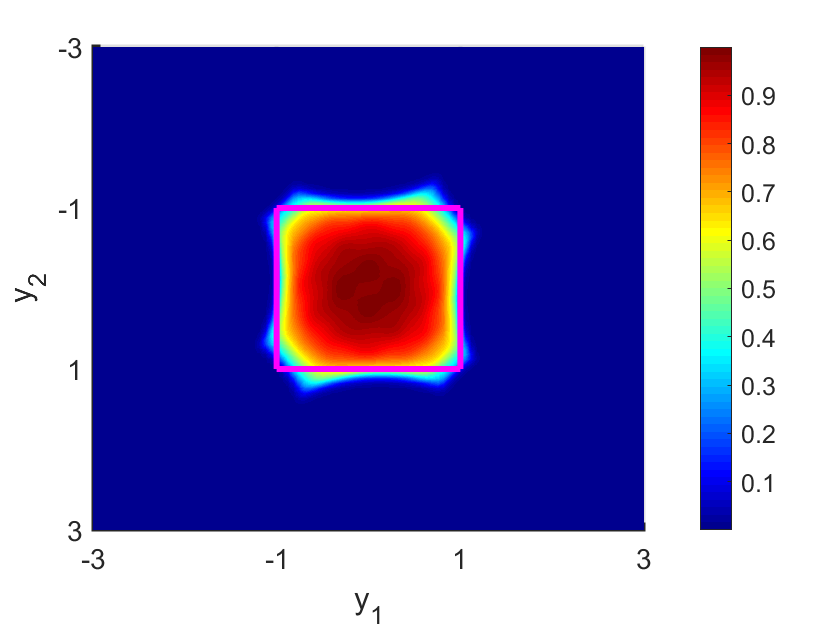}

}
\caption{Slices of reconstructions  for  a cubic support using  multi-frequency near-field data at $15$ pairs of opposite observation points.
} \label{fig:3-5}
\end{figure}

To further validate the effectiveness of our algorithm in the near-field case, we reconstruct the cubic support using  sparse observation points in Figure \ref{fig:3-4}.  We consider different numbers of pairs of opposite observation points $M=3, 7, 10, 15$. The observation points are the same as those in Figure \ref{fig:3-4}. Figure \ref{fig:3-4}(a) demonstrates that the location and shape are effectively reconstructed with only $3$ pairs of properly-selected observation points. Improved reconstruction results are achieved with $15$ pairs of opposite observation points uniformly distributed on the sphere  in Figure \ref{fig:3-4}(d). In Figure \ref{fig:3-5}, slices at $y_1=0$, $y_2=0$ and $y_3=0$ further illustrate the accuracy of the algorithm in the near-field case.

\section{Appendix}
We prove that the far-field pattern $w(\hat{x},k)$ in the frequency domain (see \eqref{u-infty}) is the inverse Fourier transform of the time-dependent far-field pattern with respect to the time variable.
Suppose that the whole space $\mathbb{R}^3$ is filled with a homogeneous and isotropic medium with a unit mass density. We designate the sound speed of the background medium as the constant $c > 0$.
The source function $S$ is supposed to be sufficiently smooth with a compact supported on $\overline{D}$ in the spatial variables. Moreover,  the wave signals are supposed to be radiated at the initial time point $t_{\min}$ and terminated at the time point $t_{\max}$. Therefore,
\begin{equation*}
 \mbox{supp }S(x,t) = \overline{D}\times [t_{\min},t_{max}] \subset B_R \times \R_+\subset \mathbb{R}^3\times \mathbb{R}_+.
\end{equation*}
The propagation of the radiated wave fields $U(x,t)$ is governed by the initial value problem
	\begin{equation}\label{TU}
		\left\{
		\begin{aligned}
			&c^{-2}\frac{\partial^2 U}{\partial t^2} = \Delta U + S(x,t), \quad &&(x, t) \in \mathbb{R}^3 \times \mathbb{R}_+, \\
			&U(x,0)=\partial_t U(x,0) = 0, &&x\in \mathbb{R}^3.
		\end{aligned}
		\right.
	\end{equation}
The solution $U(x,t)$ can be written explicitly as the convolution of the fundamental solution $G(x,t)$ to the wave equation with the source term,
	\begin{equation}\label{U}
		\begin{aligned}
		U(x,t) = G(x;t) * S(x,t) &:= \int_{\mathbb{R}^3} \int_{\mathbb{R}_+} G(x-y; t-s)S(y,s)\,dsdy\\
		&=\int_{\mathbb{R}^3} \int_{\mathbb{R}_+} \frac{\delta(t-s-c^{-1}|x-y|)}{4\pi |x-y|}S(y,s)\,dsdy,
		\end{aligned}
	\end{equation}
	where
	\begin{equation*}
		G(x;t) = \frac{\delta(t-c^{-1}|x|)}{4\pi |x|}, \quad (x,t)\in \mathbb{R}^3\backslash\{0\} \times \mathbb{R}_+.
	\end{equation*}
In the time domain,	 the far-field pattern $U^{\infty}(\hat{x},t)$ of $U(x,t)$  is defined as (see e.g., \cite{F1962})
	\begin{equation*}
	U^{\infty}(\hat{x},t):=	\lim\limits_{|x|\to +\infty} |x|\; U(x,\textcolor{black}{c^{-1}|x|}+t)\quad \mbox{ for } \hat{x}\in \s^2 \mbox{ and } t\in \R.
	\end{equation*}
By \eqref{U}, we can express
 $U^{\infty}(\hat{x},t)$  as
	\begin{equation}\label{Uf}
		\begin{aligned}
			4\pi U^{\infty}(\hat{x},t) &= 4\pi \lim\limits_{|x|\to +\infty} |x| \int_{\mathbb{R}^3} \int_{\mathbb{R}} \frac{\delta(t+c^{-1}|x|-s-c^{-1}|x-y|)}{4\pi |x-y|}S(y,s)\,dsdy\\
			&= \lim\limits_{|x|\to +\infty} \int_{B_R}\frac{|x|}{ |x-y|} \int_{\mathbb{R}} \delta(t-s+c^{-1}(|x|-|x-y|))S(y,s)\,dsdy\\
			&= \int_{B_R} \int_{\mathbb{R}} \delta(t-s+ c^{-1} \hat{x}\cdot y)S(y,s)\,dsdy\\
			&= \int_{\mathbb{R}^3} \int_{\mathbb{R}} \delta(t-s+ c^{-1} \hat{x}\cdot y)S(y,s)\,dsdy,
		\end{aligned}
	\end{equation}
which can be justified rigorously based on the definition of distribution.
Taking the inverse Fourier transform on the wave equation yields the inhomogeneous Helmholtz equations
 	\begin{equation}\label{Helmholtz-c}
 		\Delta u(x,\omega) + (\omega/c)^2 u(x,\omega) = -f(x,\omega), \qquad x\in \R^3, \;\omega>0,
 	\end{equation}
 	where the source function $f(x,\omega)$ is given by
	\begin{equation}
		f(x,\omega)=\frac{1}{\sqrt{2\pi}}\int_{\R} S(x,t)e^{i \omega t}dt.
	\end{equation}
	The solution $u(x,\omega)$ of the Helmholtz equation \eqref{Helmholtz-c} can be rephrased as
\begin{equation}
	\begin{aligned}
 		u(x,\omega)=(\mathcal{F}^{-1}U)(x,\omega) &= \sqrt{2\pi} \int_{\mathbb{R}^3} (\mathcal{F}^{-1}G)(x-y;\omega) (\mathcal{F}^{-1}S)(y,\omega)\,dy\\
 		&= \int_{\mathbb{R}^3} \Phi(x-y;\omega/c) f(y,\omega)\,dy\\
 		&= \int_{\mathbb{R}^3} \frac{e^{i\omega c^{-1}|x-y|}}{4 \pi |x-y|} f(y,\omega)\,dy.
 	\end{aligned}
\end{equation}
Here, $\Phi(x;k)$ is the outgoing fundamental solution to the Helmholtz equation $(\Delta + k^2)u = 0$, given by
\begin{equation*}
	\Phi(x;k) = \frac{e^{ik|x|}}{4 \pi |x|},\quad x\in \mathbb{R}^3,\, |x| \neq 0.
\end{equation*}
By definition of the far-field pattern in the frequency domain (see \eqref{far-field}), one can express $u^{\infty}(\hat{x},\omega)$  as
\begin{equation}
	\begin{aligned}
		u^{\infty}(\hat{x},\omega) &= \int_{\mathbb{R}^3} e^{-i\omega c^{-1}\hat{x}\cdot y} f(y,\omega)\,dy\\
		&= \int_{\mathbb{R}^3} e^{-i\omega c^{-1}\hat{x}\cdot y} \frac{1}{\sqrt{2\pi}}\int_{\R} S(y,t)e^{i \omega t}dt\,dy\\
		&= \frac{1}{\sqrt{2\pi}}  \int_{\mathbb{R}^3} \int_{\R} S(y,t) e^{i\omega (t-c^{-1} \hat{x}\cdot y)} \,dtdy.
	\end{aligned}
\end{equation}
On the other hand, the inverse Fourier transform of the time-dependent far-field pattern \eqref{Uf} can be simplified to be
\begin{equation}
	\begin{aligned}
		4\pi \left(\mathcal{F}^{-1}U^{\infty}\right) (\hat{x},\omega) &= \frac{1}{\sqrt{2\pi}} \int_{\R} \left(\int_{\mathbb{R}^3} \int_{\mathbb{R}} \delta(t-s+ c^{-1} \hat{x}\cdot y)S(y,s)\,dsdy\right) e^{i\omega t} dt \\
		&= \frac{1}{\sqrt{2\pi}} \int_{\mathbb{R}^3} \int_{\mathbb{R}} S(y,s) \int_{\mathbb{R}} \delta(t-s+ c^{-1} \hat{x}\cdot y)e^{i\omega t} \,dt\,dsdy\\
		&\xlongequal{t=s-c^{-1} \hat{x}\cdot y} \frac{1}{\sqrt{2\pi}} \int_{\mathbb{R}^3} \int_{\mathbb{R}} S(y,s) e^{i\omega (s-c^{-1} \hat{x}\cdot y)}\,dsdy\\
		&\xlongequal{s=t} \frac{1}{\sqrt{2\pi}} \int_{\mathbb{R}^3} \int_{\mathbb{R}} S(y,t) e^{i\omega (t-c^{-1} \hat{x}\cdot y)}\,dtdy\\
		&= u^{\infty}(\hat{x},\omega).
	\end{aligned}
\end{equation}
Hence, the far-field pattern $u^{\infty}(\hat{x},\omega)$ coincides with the inverse Fourier transform of the time-dependent far-field pattern $U^{\infty}(\hat{x},t)$ with respect to the time variable up to the factor $4\pi$.

\section*{Acknowledgements}
The work of G. Hu is partially supported by the National Natural Science Foundation of China (No. 12071236) and the Fundamental Research Funds for Central Universities in China (No. 63213025).


\begin{thebibliography}{00}



\bibitem{AHLS} A. Alzaalig, G. Hu, X. Liu and J. Sun, Fast acoustic source imaging using multi-frequency sparse data, Inverse Problems, 36 (2020): 025009.


\bibitem{BLLT} G. Bao, P. Li,  J. Lin and F. Triki, Inverse scattering problems with multi-frequencies, Inverse Problems, 31 (2015): 09300.

\bibitem{BLT10} G. Bao, J. Lin, and F. Triki, A multi-frequency inverse source problem, J. Differential Equations, 249 (2010): 3443--3465.

\bibitem{BLRX} G. Bao, S. Lu, W. Rundell and B. Xu, A recursive algorithm for multi-frequency acoustic inverse source
problems, SIAM J. Numer. Anal., 53 (2015): 1608-1628.



\bibitem{CHL2019} F. Cakoni, H. Haddar and A. Lechleiter,  On the factorization method for a far field inverse scattering problem in the time domain, SIAM J. Math. Anal., 51 (2019): 854-872.





\bibitem{CIL} J. Cheng, V. Isakov and S. Lu,  Increasing stability in the inverse source problem with many frequencies, J. Differential Equations, 260 (2016): 4786-4804.







\bibitem{EV09} M. Eller and N. Valdivia, Acoustic source identification using multiple frequency information, Inverse Problems, 25 (2009): 115005.


\bibitem{F1962} F. G. Friedlander,  On the radiation field of pulse solutions of the wave equation, Proc. R. Soc. Lond. A, 269 (1962): 53-65.



\bibitem{GS} R. Griesmaire and C. Schmiedecke, A Factorization method for multifrequency inverse source problem with sparse far-field measurements, SIAM J. Imaging Sciences, 10 (2017): 2119-2139.





\bibitem{GHS}R. Griesmaier, M. Hanke and J. Sylvester, Far field splitting for the Helmholtz equation, SIAM Journal on Numerical Analysis,
52 (2014): 343-362.


\bibitem{GS17} R. Griesmaier and C. Schmiedecke, A multifrequency MUSIC algorithm for locating small inhomogeneities in inverse scattering, Inverse Problems, 33 (2017): 035015

\bibitem{HL2020} H. Haddar and X. Liu,  A time domain factorization method for obstacles with impedance boundary conditions, Inverse Problems, 36 (2020): 105011.





\bibitem{GGH2022} H. Guo and G. Hu, Inverse wave-number-dependent source problems for the Helmholtz equation, arXiv:2305.07459.


\bibitem{GHZ}H. Guo, G. Hu and M. Zhao, Direct sampling method to inverse wave-number-dependent source problems: determination of the support of a stationary source,  Inverse Problems, 39 (2023): 105008.

\bibitem{GHM}H. Guo, G. Hu and G. Ma, Imaging a moving point source from multi-frequency data measured at one and sparse observation directions (part I): far-field case, SIAM Journal of Imaging Sciences, 16 (2023): 1535-1571.

\bibitem{HL2020} G. Hu and J. Li, Uniqueness to inverse source problems in an inhomogeneous medium with a single far-field pattern, SIAM J. Math. Anal., 52 (2020): 5213-5231.










\bibitem{LMZ}
X. Liu, S. Meng and B. Zhang, Modified sampling method with near field measurements, SIAM J. Appl. Math., 82 (2022): 244-266.

\bibitem{JLZ2019} X.  Ji, X. Liu and B. Zhang,  Phaseless inverse source scattering problem: phase retrieval, uniqueness and direct sampling methods, J. Comput. Phys. X, 1 (2019): 100003.

\bibitem{K98} A. Kirsch, Characterization of the shape of a scattering obstacle using the spectral data of the far field
operator, Inverse Problems, 14 (1998):  1489-1512.

\bibitem{KG08} A. Kirsch and N. Grinberg, {\rm The Factorization Method for Inverse Problems}, Oxford University Press,
Oxford, UK, 2008.


\bibitem{LY} P.  Li and G. Yuan, Increasing stability for the inverse source scattering problem with multi-frequencies, Inverse Problems and Imaging, 11 (2017): 745-759.







\bibitem{ZG} D. Zhang and Y. Guo,  Fourier method for solving the multi-frequency inverse source problem for the Helmholtz equation, Inverse Problems, 31 (2015): 035007.







\end{thebibliography}
\end{document}